\documentclass[11pt]{article}
\usepackage{amsmath}
\usepackage{amssymb}
\usepackage{amscd}
\usepackage{xspace}
\usepackage{verbatim}
\newcommand{\mysection}[1]{\section{#1}\setcounter{equation}{0}}
\title{\bf Boundary trace of positive solutions of semilinear elliptic equations in Lipschitz domains}
\author{ {\bf Moshe Marcus}\\
{\small Department of Mathematics,}\\
 {\small  Technion, Haifa, ISRAEL}
\and {\bf Laurent Veron}\\
{\small Department of Mathematics,}\\
 {\small  Univ. of Tours,  FRANCE}
}

\date{}
\begin{document}

\numberwithin{equation}{section}
\maketitle

\newcommand{\txt}[1]{\;\text{ #1 }\;}
\newcommand{\tbf}{\textbf}
\newcommand{\tit}{\textit}
\newcommand{\tsc}{\textsc}
\newcommand{\trm}{\textrm}
\newcommand{\mbf}{\mathbf}
\newcommand{\mrm}{\mathrm}
\newcommand{\bsym}{\boldsymbol}
\newcommand{\scs}{\scriptstyle}
\newcommand{\sss}{\scriptscriptstyle}
\newcommand{\txts}{\textstyle}
\newcommand{\dsps}{\displaystyle}
\newcommand{\fnz}{\footnotesize}
\newcommand{\scz}{\scriptsize}
\newcommand{\be}{
\begin{equation}
}
\newcommand{\bel}[1]{
\begin{equation}
\label{#1}}
\newcommand{\ee}{
\end{equation}
}
\newcommand{\eqnl}[2]{
\begin{equation}
\label{#1}{#2}
\end{equation}
}
\newtheorem{subn}{\name}
\renewcommand{\thesubn}{}
\newcommand{\bsn}[1]{\def\name{#1}
\begin{subn}}
\newcommand{\esn}{
\end{subn}}
\newtheorem{sub}{\name}[section]
\newcommand{\dn}[1]{\def\name{#1}}   
\newcommand{\bs}{
\begin{sub}}
\newcommand{\es}{
\end{sub}}
\newcommand{\bsl}[1]{
\begin{sub}\label{#1}}
\newcommand{\bth}[1]{\def\name{Theorem}
\begin{sub}\label{t:#1}}
\newcommand{\blemma}[1]{\def\name{Lemma}
\begin{sub}\label{l:#1}}
\newcommand{\bcor}[1]{\def\name{Corollary}
\begin{sub}\label{c:#1}}
\newcommand{\bdef}[1]{\def\name{Definition}
\begin{sub}\label{d:#1}}
\newcommand{\bprop}[1]{\def\name{Proposition}
\begin{sub}\label{p:#1}}
\newcommand{\R}{\eqref}
\newcommand{\rth}[1]{Theorem~\ref{t:#1}}
\newcommand{\rlemma}[1]{Lemma~\ref{l:#1}}
\newcommand{\rcor}[1]{Corollary~\ref{c:#1}}
\newcommand{\rdef}[1]{Definition~\ref{d:#1}}
\newcommand{\rprop}[1]{Proposition~\ref{p:#1}}
\newcommand{\BA}{
\begin{array}}
\newcommand{\EA}{
\end{array}}
\newcommand{\BAN}{\renewcommand{\arraystretch}{1.2}
\setlength{\arraycolsep}{2pt}
\begin{array}}
\newcommand{\BAV}[2]{\renewcommand{\arraystretch}{#1}
\setlength{\arraycolsep}{#2}
\begin{array}}
\newcommand{\BSA}{
\begin{subarray}}
\newcommand{\ESA}{\end{subarray}}
\newcommand{\BAL}{\begin{aligned}}
\newcommand{\EAL}{\end{aligned}}
\newcommand{\BALG}{\begin{alignat}}
\newcommand{\EALG}{\end{alignat}}
\newcommand{\BALGN}{\begin{alignat*}}
\newcommand{\EALGN}{\end{alignat*}}
\newcommand{\qeda}{\hspace{10mm}\hfill $\square$}
\newcommand{\qed}{\\
${}$ \hfill $\square$}
\newcommand{\Remark}{\note{Remark}}
\newcommand{\forevery}{\quad \forall}
\newcommand{\set}[1]{\{#1\}}
\newcommand{\lra}{\longrightarrow}
\newcommand{\sgn}{\rm{sgn}}
\newcommand{\lla}{\longleftarrow}
\newcommand{\llra}{\longleftrightarrow}
\newcommand{\Lra}{\Longrightarrow}
\newcommand{\Lla}{\Longleftarrow}
\newcommand{\Llra}{\Longleftrightarrow}
\newcommand{\warrow}{\rightharpoonup}
\newcommand{
\paran}[1]{\left (#1 \right )}
\newcommand{\sqbr}[1]{\left [#1 \right ]}
\newcommand{\curlybr}[1]{\left \{#1 \right \}}
\newcommand{\abs}[1]{\left |#1\right |}
\newcommand{\norm}[1]{\left \|#1\right \|}
\newcommand{\paranb}[1]{\big (#1 \big )}
\newcommand{\lsqbrb}[1]{\big [#1 \big ]}
\newcommand{\lcurlybrb}[1]{\big \{#1 \big \}}
\newcommand{\absb}[1]{\big |#1\big |}
\newcommand{\normb}[1]{\big \|#1\big \|}
\newcommand{
\paranB}[1]{\Big (#1 \Big )}
\newcommand{\absB}[1]{\Big |#1\Big |}
\newcommand{\normB}[1]{\Big \|#1\Big \|}

\newcommand{\thkl}{\rule[-.5mm]{.3mm}{3mm}}
\newcommand{\thknorm}[1]{\thkl #1 \thkl\,}
\newcommand{\trinorm}[1]{|\!|\!| #1 |\!|\!|\,}
\newcommand{\bang}[1]{\langle #1 \rangle}
\def\angb<#1>{\langle #1 \rangle}
\newcommand{\vstrut}[1]{\rule{0mm}{#1}}
\newcommand{\rec}[1]{\frac{1}{#1}}
\newcommand{\opname}[1]{\mbox{\rm #1}\,}
\newcommand{\supp}{\opname{supp}}
\newcommand{\dist}{\opname{dist}}
\newcommand{\myfrac}[2]{{\displaystyle \frac{#1}{#2} }}
\newcommand{\myint}[2]{{\displaystyle \int_{#1}^{#2}}}
\newcommand{\mysum}[2]{{\displaystyle \sum_{#1}^{#2}}}
\newcommand {\dint}{{\displaystyle \int\!\!\int}}
\newcommand{\q}{\quad}
\newcommand{\qq}{\qquad}
\newcommand{\hsp}[1]{\hspace{#1mm}}
\newcommand{\vsp}[1]{\vspace{#1mm}}
\newcommand{\prt}{\partial}
\newcommand{\sms}{\setminus}
\newcommand{\ems}{\emptyset}
\newcommand{\ti}{\times}
\newcommand{\nind}{\noindent}
\newcommand{\pr}{^\prime}
\newcommand{\ppr}{^{\prime\prime}}
\newcommand{\tl}{\tilde}
\newcommand{\wtl}{\widetilde}
\newcommand{\sbs}{\subset}
\newcommand{\sbeq}{\subseteq}
\newcommand{\indx}[1]{_{\scriptscriptstyle #1}}
\newcommand{\ovl}{\overline}
\newcommand{\unl}{\underline}
\newcommand{\nin}{\not\in}
\newcommand{\pfrac}[2]{\genfrac{(}{)}{}{}{#1}{#2}}

\def\ga{\alpha}     \def\gb{\beta}       \def\gg{\gamma}
\def\gc{\chi}       \def\gd{\delta}      \def\ge{\epsilon}
\def\gth{\theta}                         \def\vge{\varepsilon}
\def\gf{\phi}       \def\vgf{\varphi}    \def\gh{\eta}
\def\gi{\iota}      \def\gk{\kappa}      \def\gl{\lambda}
\def\gm{\mu}        \def\gn{\nu}         \def\gp{\pi}
\def\vgp{\varpi}    \def\gr{\rho}        \def\vgr{\varrho}
\def\gs{\sigma}     \def\vgs{\varsigma}  \def\gt{\tau}
\def\gu{\upsilon}   \def\gv{\vartheta}   \def\gw{\omega}
\def\gx{\xi}        \def\gy{\psi}        \def\gz{\zeta}
\def\Gg{\Gamma}     \def\Gd{\Delta}      \def\Gf{\Phi}
\def\Gth{\Theta}
\def\Gl{\Lambda}    \def\Gs{\Sigma}      \def\Gp{\Pi}
\def\Gw{\Omega}     \def\Gx{\Xi}         \def\Gy{\Psi}

\def\CS{{\mathcal S}}   \def\CM{{\mathcal M}}   \def\CN{{\mathcal N}}
\def\CR{{\mathcal R}}   \def\CO{{\mathcal O}}   \def\CP{{\mathcal P}}
\def\CA{{\mathcal A}}   \def\CB{{\mathcal B}}   \def\CC{{\mathcal C}}
\def\CD{{\mathcal D}}   \def\CE{{\mathcal E}}   \def\CF{{\mathcal F}}
\def\CG{{\mathcal G}}   \def\CH{{\mathcal H}}   \def\CI{{\mathcal I}}
\def\CJ{{\mathcal J}}   \def\CK{{\mathcal K}}   \def\CL{{\mathcal L}}
\def\CT{{\mathcal T}}   \def\CU{{\mathcal U}}   \def\CV{{\mathcal V}}
\def\CZ{{\mathcal Z}}   \def\CX{{\mathcal X}}   \def\CY{{\mathcal Y}}
\def\CW{{\mathcal W}} \def\CQ{{\mathcal Q}}
\def\BBA {\mathbb A}   \def\BBb {\mathbb B}    \def\BBC {\mathbb C}
\def\BBD {\mathbb D}   \def\BBE {\mathbb E}    \def\BBF {\mathbb F}
\def\BBG {\mathbb G}   \def\BBH {\mathbb H}    \def\BBI {\mathbb I}
\def\BBJ {\mathbb J}   \def\BBK {\mathbb K}    \def\BBL {\mathbb L}
\def\BBM {\mathbb M}   \def\BBN {\mathbb N}    \def\BBO {\mathbb O}
\def\BBP {\mathbb P}   \def\BBR {\mathbb R}    \def\BBS {\mathbb S}
\def\BBT {\mathbb T}   \def\BBU {\mathbb U}    \def\BBV {\mathbb V}
\def\BBW {\mathbb W}   \def\BBX {\mathbb X}    \def\BBY {\mathbb Y}
\def\BBZ {\mathbb Z}

\def\GTA {\mathfrak A}   \def\GTB {\mathfrak B}    \def\GTC {\mathfrak C}
\def\GTD {\mathfrak D}   \def\GTE {\mathfrak E}    \def\GTF {\mathfrak F}
\def\GTG {\mathfrak G}   \def\GTH {\mathfrak H}    \def\GTI {\mathfrak I}
\def\GTJ {\mathfrak J}   \def\GTK {\mathfrak K}    \def\GTL {\mathfrak L}
\def\GTM {\mathfrak M}   \def\GTN {\mathfrak N}    \def\GTO {\mathfrak O}
\def\GTP {\mathfrak P}   \def\GTR {\mathfrak R}    \def\GTS {\mathfrak S}
\def\GTT {\mathfrak T}   \def\GTU {\mathfrak U}    \def\GTV {\mathfrak V}
\def\GTW {\mathfrak W}   \def\GTX {\mathfrak X}    \def\GTY {\mathfrak Y}
\def\GTZ {\mathfrak Z}   \def\GTQ {\mathfrak Q}

\font\Sym= msam10 
\def\SYM#1{\hbox{\Sym #1}}
\newcommand{\tin}{\to\infty}
\newcommand{\ssub}[1]{_{_{\! #1}}}
\newcommand{\chr}[1]{\chi\indx{#1}}
\newcommand{\rest}[1]{\big |\indx{#1}}
\newcommand{\bdw}{\prt\Gw\xspace}
\newcommand{\wkc}{weak convergence\xspace}
\newcommand{\wrto}{with respect to\xspace}
\newcommand{\cons}{consequence\xspace}
\newcommand{\consy}{consequently\xspace}
\newcommand{\Consy}{Consequently\xspace}
\newcommand{\Essy}{Essentially\xspace}
\newcommand{\essy}{essentially\xspace}
\newcommand{\mnz}{minimizer\xspace}
\newcommand{\sth}{such that\xspace}
\newcommand{\ngh}{neighborhood\xspace}
\newcommand{\nghs}{neighborhoods\xspace}
\newcommand{\seq}{sequence\xspace}
\newcommand{\seqs}{sequences\xspace}
\newcommand{\sseq}{subsequence\xspace}
\newcommand{\ifif}{if and only if\xspace}
\newcommand{\suff}{sufficiently\xspace}
\newcommand{\abc}{absolutely continuous\xspace}
\newcommand{\sol}{solution\xspace}
\newcommand{\subss}{sub-solutions\xspace}
\newcommand{\subs}{sub-solution\xspace}
\newcommand{\supers}{super-solution\xspace}
\newcommand{\superss}{super-solutions  \xspace}
\newcommand{\Wlg}{Without loss of generality\xspace}
\newcommand{\wlg}{without loss of generality\xspace}
\newcommand{\locun}{locally uniformly\xspace}
\newcommand{\bvp}{boundary value problem\xspace}
\newcommand{\bvps}{boundary value problems\xspace}
\def\RN{\BBR^N}
\def\({{\rm (}}
\def\){{\rm )}}
\def\loc{\indx{\rm loc}}
\def\bmn{\mathbf{n}}
\def\bma{\mathbf{a}}
\def\prtn{\prt_{\bmn}}
\def\1{\\[1mm]}
\def\2{\\[2mm]}
\def\Lip{Lipschitz\xspace}
\def\BHP{boundary Harnack principle \xspace}
\def\Note{\nind\textit{Note.}\hskip 2mm}
\def\Remark{\nind\textit{Remark.}\hskip 2mm}
\def\Notation{\nind\textit{Notation.}\hskip 2mm}
\def\Proof{\nind\textit{Proof.}\hskip 2mm}
\newcommand{\tr}[1]{\text{\rm tr}\indx{#1}}
\medskip
\begin{center}
{\bf\small Abstract}

\vspace{3mm} \hspace{.05in}\parbox{4.5in} {{\small We study the
generalized boundary value problem for nonnegative solutions of of
$-\Delta u+g(u)=0$ in a bounded Lipschitz domain $\Gw$, when $g$ is
continuous and nondecreasing. Using the harmonic measure of $\Gw$,
we define a trace in the class of outer regular Borel measures. We
amphasize the case where $g(u)=|u|^{q-1}u$, $q>1$. When $\Gw$ is
(locally) a cone with vertex $y$, we prove sharp results of
removability and characterization of singular behavior. In the
general case, assuming that $\Gw$ possesses  a tangent cone at every
boundary point and $q$ is subcritical, we prove an existence and
uniqueness result for positive solutions with arbitrary  boundary
trace.  We obtain sharp results  involving  Besov spaces with
negative index  on k-dimensional edges and apply our results to the
characterization of removable sets and good measures on the boundary
of a polyhedron.}}
\end{center}


\newpage
\tableofcontents

\mysection {Introduction} In this article we study boundary value
problems with measure data on the boundary, for equations of the
form
\begin{equation}\label{M1}
-\Gd u+g(u)=0 \txt{in}\Gw
\end {equation}
 where $\Gw$ is a bounded {\em \Lip domain} in $\BBR^N$ and $g$ is a continuous nondecreasing
 function vanishing at $0$. A function $u$ is a solution of the equation
 if $u$ and $g(u)$ belong to $L^1\loc(\Gw)$ and the equation holds
 in the distribution sense. The definition of a solution satisfying
 a prescribed boundary condition is more complex and will be
 described later on.

Boundary value problems for \eqref{M1} with measure boundary data in
smooth domains (or, more precisely, in $C^2$ domains) have been
studied intensively in the last 20 years. Much of this work
concentrated on the case of power nonlinearities, namely,
$g(u)=|u|^{q-1}u$ with $q>1$. For details we address the reader to
the following papers and the references therein: Le Gall [1-2],
Dynkin and Kuznetsov [1-3], Mselati [1] (employing in an essential
way probabilistic tools) and Marcus and Veron [1-4] (employing
purely analytic methods).

The study of the corresponding linear \bvp in \Lip domains is
classical. This study shows that, with a proper interpretation,
 the basic results known for smooth domains remain valid in the \Lip
 case. Of course there are important differences too: in the Poisson integral formula
 the Poisson kernel must be replaced by the Martin kernel and, when the boundary
 data is given by a function in $L^1$, the standard surface measure must be
 replaced by the harmonic measure.
 The Hopf principle does
 not hold anymore, but it is partially replaced by the Carleson lemma
 and the \BHP due to Dahlberg \cite{Da}. A summary
of the basic results for the linear case, to the extent needed in
the present work, is presented in Section 2.

One might expect that in the nonlinear case  the results valid
for smooth domains extend to \Lip domains in a similar way.
This is indeed the case as long as the boundary data is in
$L^1$. However, in problems with measure boundary data, we
encounter essentially new phenomena.

Following is an overview of our main results on boundary value
problems for \eqref{M1}.\2

\nind{\bf A.}\hskip 2mm {\em General nonlinearity and finite measure
data.}

\medskip

\nind\/ We start with the weak $L^1$ formulation of the \bvp
\begin{equation}\label{pref-bvp}
  -\Gd u+ g(u)=0\txt{in $\Gw$,} u=\mu \txt{on $\bdw$},
\end{equation}
where $\mu\in \GTM(\bdw)$.

Let $x_0$ be a point in $\Gw$, to be kept fixed, and let
$\gr=\gr\indx{\Gw}$ denote the first eigenfunction of $-\Gd$ in
$\Gw$ normalized by $\gr(x_0)=1$. It turns out that the family of
test functions appropriate for the \bvp is
\begin{equation}\label {pref-X}
X(\Gw)=\left\{\eta\in W^{1,2}_{0}(\Gw):\gr^{-1}\Gd \eta\in
L^{\infty}(\Gw)\right\}.
\end{equation}
If $\eta\in\Gw$ then $\sup|\eta|/\gr<\infty$.

 Let $\BBK[\mu]$ denote the harmonic function in $\Gw$ with
boundary trace $\mu$. Then $u$ is an {\em $L^1$-weak solution\/} of
\eqref{pref-bvp} if
\begin{equation}\label{pref-weak'}
u\in L^1_\gr(\Gw),\q g(u)\in L^1_\gr(\Gw)
\end{equation}
and
\begin{equation}\label{pref-weak}
\int_{\Gw}\left(-u\Gd\eta+g(u)\eta\right)dx=-\int_{\Gw}\left(\BBK[\mu]\Gd\eta\right)dx\forevery
\eta\in X(\Gw). \end{equation}

Note that in \eqref{pref-weak} the boundary data appears only in an
implicit form. In the next result we present a more explicit link
between the solution and its boundary trace.

A \seq of domains $\set{\Gw_n}$ is called  a {\em \Lip exhaustion}
of $\Gw$ if,  for every $n$, $\Gw_n$ is \Lip and
\begin{equation}\label{pref-exhaustion}
  \Gw_n\sbs \bar \Gw_n\sbs \Gw_{n+1},\q \Gw=\cup\Gw_n,\q
  \BBH_{N-1}(\bdw_n)\to \BBH_{N-1}(\bdw).
\end{equation}

 \bprop{pref-tr-conv} Let $\set{\Gw_n}$ be an exhaustion of $\Gw$,
 let $x_0\in \Gw_1$ and denote by $\gw_n$
(respectively $\gw$) the harmonic measure on $\bdw_n$ (respectively
$\bdw$) relative to $x_0$. If $u$ is an $L^1$-weak solution of
\eqref{pref-bvp} then, for every $Z\in C(\bar \Gw)$,
\begin{equation}\label{pref-tr-conv}
   \lim_{n\to\infty}\int_{\bdw_n}Zu\,d\gw_n=\int_{\bdw}Z\,d\mu.
\end{equation}
\es

We note that any solution of \eqref{M1} is in $W^{1,p}\loc(\Gw)$ for
some $p>1$ and \consy possesses an integrable trace on $\bdw_n$.

In general problem \eqref{pref-bvp} does not possess a solution for
every $\mu$. We denote by $\GTM^g(\bdw)$ the set of measures $\mu\in
\GTM(\bdw)$ for which such a solution exists. The following
statements are established in the same way as in the case of smooth
domains:\2 (i) If a solution exists it is unique. Furthermore the
solution
depends monotonically on the boundary data.\\
(ii) If $u$ is an $L^1$-weak solution of \eqref{pref-bvp} then $|u|$
(resp. $u_+$) is a subsolution of this problem with $\mu$ replaced
by $|\mu|$ (resp. $\mu_+$).

A measure $\mu\in\GTM(\bdw)$ is {\em $g$-admissible} if
$g(\BBK[|\mu|])\in L^1_\gr(\Gw)$. When there is no risk of
confusion we shall simply write 'admissible' instead of
'$g$-admissible'. The following provides a sufficient condition
for existence.

\bth{pref-admissible} If $\mu$ is $g$-admissible then problem
\eqref{pref-bvp} possesses a unique solution. \es

\nind{\bf B.} {\em The boundary trace of positive solutions of
\eqref{M1}; general nonlinearity.}

\medskip
\nind\/ We say that $u\in L^1\loc(\Gw)$ is a {\em regular solution}
of the equation \eqref{M1} if $g(u)\in L^1_\gr(\Gw)$.

 \bprop{pref-reg-sol} Let $u$ be a positive
solution of the equation \eqref{M1}. If $u$ is regular then
$u\in L^1_\gr(\Gw)$ and it possesses a boundary trace
$\mu\in\GTM(\bdw)$. Thus $u$ is the solution of the \bvp
\eqref{pref-bvp} with this measure $\mu$. \es

As in the case of smooth domains, a positive solution possesses a
boundary trace even if the solution is not regular. The boundary
trace may be defined in several ways; in every case it is expressed
by an {\em unbounded measure.} A definition of trace is 'good' if
the trace uniquely determines the solution. A discussion of the
various definitions of boundary trace, for boundary value problems
in $C^2$ domains, with power nonlinearities, can be found in
\cite{MV6}, \cite{Dy1} and the references therein. In \cite{MV1} the
authors introduced a definition of trace -- later referred to  as
the 'rough trace' by Dynkin \cite{Dy1} -- which proved to be 'good'
in the subcritical case, but not in the supercritical case (see
\cite{MV2}). Mselati \cite{Ms} obtained a 'good' definition of trace
for the problem with $g(u)=u^2$ and $N\geq 4$, in which case this
non-linearity is supercritical. His approach employed probabilistic
methods developed by Le Gall in a series of papers. For a
presentation of these methods  we refer the reader to his book
\cite{LeG-book}. Following this work the authors introduced in
\cite{MV6} a notion of trace, called 'the precise trace', defined in
the framework of the fine topology associated with the capacity
$C_{2/q,q'}$ on $\bdw$. This definition of trace turned out to be
'good' for all power nonlinearities $g(u)=u^q$, $q>1$, at least in
the class of $\gs$-moderate solutions.  In the subcritical case, the
precise trace reduces to the rough trace. At the same time Dynkin
\cite{Dy2} extended Mselati's result to the case $(N+1)/(N-1)\leq
q\leq 2$.

In the present paper we confine ourselves to boundary value problems
with rough trace data. (See the definition below.) However we
develop a framework for the study of existence and uniqueness (see
\rth{pref-reg+sing} below) which can be applied to a large class of
nonlinearities  and can be adapted to other notions of trace as
well. In particular, it can be adapted to the 'precise trace' for
power nonlinearities (in smooth domains) and to a related notion of
trace for \Lip domains. This issue will be addressed in a subsequent
paper.

Here are the main results in this part of the paper, including the
relevant definitions.

 \bdef{pref-trace} Let $u$ be a positive
supersolution, respectively subsolution, of \eqref{M1}. A point
$y\in \bdw$ is a {\em regular boundary point} relative to $u$ if
there exists an open \ngh $D$ of $y$ \sth $g\circ u\in
L^1_\gr(\Gw\cap D)$. If no such \ngh exists we say that $y$ is a
{\em singular boundary point} relative to $u$.

The set of regular boundary points of $u$ is denoted by $\CR(u)$;
its complement on the boundary is denoted by $\CS(u)$. Evidently
$\CR(u)$ is relatively open.
 \es

\bth{pref-reg-trace} Let $u$ be a positive solution of \eqref{M1} in
$\Gw$. Then $u$ possesses a trace  on $\CR(u)$, given by a Radon
measure $\nu$.

Furthermore, for every compact set $F\sbs \CR(u)$,
\begin{equation}\label{pref-local}
\int_{\Gw}\left(-u\Gd\eta+g(u)\eta\right)dx=-\int_{\Gw}\left(\BBK[\gn\chi\indx{F}]\Gd\eta\right)dx
\end{equation}
for every $\eta\in X(\Gw)$ \sth $\supp \eta\cap\bdw\sbs F$ and
$\nu\chi\indx{F}\in \GTM^g(\bdw)$. \es

\bdef{pref-g-trace} Let $g\in \CG$. Let $u$ be a positive solution
of \eqref{M1} with regular boundary set $\CR(u)$ and singular
boundary set $\CS(u)$. The Radon measure $\nu$ in $\CR(u)$
associated with $u$  as in \rth{pref-reg-trace} is called the {\em
regular part of the trace of $u$}. The couple $(\nu,\CS(u))$ is
called the {\em boundary trace} of $u$ on $\bdw$. This trace is also
represented by the (possibly unbounded) Borel measure $\bar\nu$
given by
\begin{equation}\label{pref-barnu}
 \bar\nu(E)=\begin{cases} \nu(E), &\text{if $E\sbs
\CR(u)$}\\ \infty,&\text{otherwise.}\end{cases}
\end{equation}

The boundary trace of $u$ in the sense of this definition will be
denoted by $\tr{\bdw}u$.

Let
\begin{equation}\label{pref-Vnu}
V_\nu:=\sup\set{u\indx{\nu\chi_F}: \;F\sbs \CR(u),\;F\text{
compact}}
\end{equation}
where $u\indx{\nu\chi_F}$ denotes the solution of \eqref{pref-bvp}
with $\mu=\nu\chi_F$. Then $V_\nu$ is called the {\em semi-regular
component\/} of $u$. \es

\bdef{pref-removable} A compact set $F\sbs \bdw$ is {\em removable}
relative to \eqref{M1} if the only non-negative solution $u\in
C(\bar\Gw\sms F)$ which vanishes on $\bar\Gw\sms F$ is the trivial
solution $u=0$. \es

\blemma{pref-maxsol} Let $g\in\CG$ and assume that $g$ satisfies the
Keller-Osserman condition. Let $F\sbs\bdw$ be a compact set and
denote by $\CU_F$ the class of solutions $u$ of \eqref{M1} which
satisfy the condition,
\begin{equation}\label{pref-Fc-vanish}
   u\in C(\bar\Gw\sms F),\q u=0\txt{on $\bdw\sms F$.}
\end{equation}
Then there exists a function $U_F\in \CU_F$ \sth
$$u\leq U_F \forevery u\in \CU_F.$$

\nind Furthermore, $\CS(U_F)=:F'\sbs F$; $F'$ need not be equal to
$F$. \es

\bdef{pref-maxsol} $U_F$ is called the {\em maximal solution}
associated with $F$. The set $F'=\CS(U_F)$ is  called the $g$-kernel
of $F$ and denoted by $k_g(F)$. \es

\bth{pref-reg+sing} Let $g\in\CG$ and assume that $g$ is convex and
satisfies the Keller-Osserman condition.

\nind{\sc Existence.}  The following set of conditions is  necessary
and sufficient  for existence of a solution $u$ of the generalized
\bvp
\begin{equation}\label{gen-bvp}
  -\Gd u+ g(u)=0 \txt{in}\Gw,\q \tr{\bdw}u=(\nu,F),
\end{equation}
where $F\sbs\bdw$ is a compact set and $\nu$ is a Radon measure on
$\bdw\sms F$.

(i)   For every compact set $E\sbs \bdw\sms F$, $\nu\chi\indx{E}\in
\GTM^g(\bdw)$.

(ii) If $k_g(F)=F'$, then $F\sms F'\sbs\CS(V_\nu)$.

\smallskip

\nind When this holds,
\begin{equation}\label{pref u<V+U}
V_\nu\leq  u\leq V_\nu + U_F.
\end{equation}
Furthermore if $F$ is a removable set then \eqref{pref-bvp}
possesses exactly one solution.

\medskip
\nind{\sc Uniqueness.} Given a compact set $F\sbs \bdw$, assume that
\begin{equation}\label{pref-uniq1}
  U_{E} \text{ is the unique solution with trace
 $(0,k_g(E))$}
\end{equation}
for every compact $E\sbs F$. Under this assumption:

\nind{\rm (a)} If $u$ is a solution of \eqref{gen-bvp} then
\begin{equation}\label{pref u>max(V,U)}
\max(V_\nu,U_F)\leq u\leq V_\nu + U_F.
\end{equation}

\nind{\rm (b)} Equation \eqref{M1} possesses at most one solution
satisfying \eqref{pref u>max(V,U)}.

\nind{\rm (c)} Condition \eqref{pref-uniq1} is necessary and
sufficient in order that \eqref{gen-bvp} possess at most one
solution.

\medskip

\nind{\sc Monotonicity.}

 \nind{\rm (d)} Let $u_1,u_2$
be two positive solutions of \eqref{M1} with boundary traces
$(\nu_1,F_1)$ and $(\nu_2,F_2)$ respectively. Suppose that $F_1\sbs
F_2$ and that $\nu_1\leq \nu_2\chi\indx{F_1}=:\nu'_2$.
If \eqref{pref-uniq1} holds for $F=F_2$ then $u_1\leq u_2$. \es

In the remaining part of this paper we consider equation \eqref{M1}
with power nonlinearity:
\begin{equation}\label{Mq}
  -\Gd u+|u|^{q-1}u=0
\end{equation}
with $q>1$.

\medskip
\nind{\bf C.} {\em Classification of positive solutions in a conical
domain possessing an isolated singularity at the vertex.}

\medskip

Let $C_{_{S}}$ be a cone with vertex  $0$ and opening $S\subset
S^{N-1}$, where $S$ is a \Lip domain. Put $\Gw=C_S\cap B_1(0)$.
Denote by $\gl\indx{S}$ the first eigenvalue and by $\gf\indx{S}$
the first eigenfunction of $-\Gd'$ in $W^{1,2}_{0}(S)$ normalized by
$\max \gf_{_{S}}=1.$  Put

$$\ga_S=\rec{2}(N-2+\sqrt{(N-2)^2+4\gl\indx{S}}$$
and
$$\Gf_1=\rec{\gg} x^{-\ga_{_S}}\gf_{_{S}}(x/\abs x)$$
where  $\gg\indx{S}$ is a positive number. $\Gf_1$ is a harmonic
function in $C_S$ vanishing on $\prt C_S\sms \{0\}$ and  $\gg$ is
chosen so that the boundary trace of $\Gf_1$ is $\gd_0$ (=Dirac
measure on $\prt C_S$ with mass $1$ at the origin). Further denote
$\Gw_S=C_S\cap B_{1}(0)$.

It was shown in \cite{FV} that, if $q\geq 1+\frac{2}{\ga\indx{S}}$
there is no solution of \eqref{Mq} in $\Gw$ with isolated
singularity at $0$. We obtain the following result.

\bth{pref-cone1} Assume that $1<q<1+\frac{2}{\ga\indx{S}}$. Then
$\gd_0$ is admissible for $\Gw$ and \consy, for every real $k$,
there exists a unique solution of this equation in $\Gw$ with
boundary trace $k\gd_0$. This solution, denoted by $u_k$ satisfies
\begin{equation}\label{pref-lim-k}
u_{k}(x)=k\Gf_1(x)(1+o(1))\quad\text{as }x\to 0.
\end{equation}

The function $$u_\infty=\lim_{k\tin} u_k$$  is the unique positive
solution of \eqref{eq-q} in $\Gw_S$ which vanishes on $\prt \Gw\sms
\{0\}$ and is strongly singular at $0$, i.e.,
\begin{equation}\label{ssing}
 \int_\Gw u_\infty^q \gr\, dx=\infty
\end{equation}
where $\gr$ is the first eigenfunction of $-\Gd$ in $\Gw$ normalized
by $\gr(x_0)=1$ for some (fixed) $x_0\in \Gw$. Its asymptotic
behavior at $0$ is given by,
\begin{equation}\label{pref-lim-infi}
u_{\infty}(x)= |x|^{-\frac{2}{q-1}}\gw_{S}(x/|x|)(1+o(1))\quad
\txt{as}x\to0
\end{equation}
where  $\gw$ is the (unique) positive solution of
\begin{equation}\label{pref-NEVP}
-\Gd'\gw-\gl_{_{N,q}}\gw+\abs\gw^{q-1}\gw=0
\end{equation}
on $S^{N-1}$ with
\begin{equation}\label{pref-NEVP1}
\gl_{_{N,q}}=\myfrac{2}{q-1}\left(\myfrac{2q}{q-1}-N\right).
\end{equation}
\es

As a consequence one can state the following classification result.

\bth{pref-cone2} Assume that $1<q<q_{_S}=1+2/\ga_{_S}$ and denote
$$\tl\ga\indx{S}=\rec{2}\big(2-N+\sqrt{(N-2)^2+4\gl\indx{S}}\,\big).$$
If $u\in C(\bar\Gw_{S}\setminus\{0\})$ is a  positive solution of
(\ref{Mq}) vanishing on $(\prt C_{_{S}}\cap
B_{r_{0}}(0))\sms\{0\}$, the following alternative holds:

\smallskip
Either
  $$\limsup_{x\to 0}\abs
x^{-\tilde\ga_{_{S}}}u(x)<\infty,
$$

\smallskip
or $$\text{there exists $k>0$ such that (\ref{pref-lim-k}) holds,}$$

\smallskip
or
 $$\text{(\ref{pref-lim-infi}) holds.}$$
 \es

In the first case $u\in C(\bar\Gw$; in the second, $u$ possesses a
{\em weak singularity} at the vertex while in the last case $u$ has
a {\em strong singularity} there.

\medskip
\nind{\bf D.} {\em Criticality in \Lip domains.}

\medskip
Let $\Gw$ be a \Lip domain and let $\gx\in \bdw$. We say that
$q_\gx$ is the {\em critical value for \eqref{Mq} at $\gx$\/} if,
for $1<q<q_\gx$, the equation possesses a solution with boundary
trace $\gd_\gx$ while, for $q> q_\gx$ no such solution exists. We
say that $q^\sharp_\gx$ is the {\em secondary critical value at
$\gx$\/} if for $1<q<q^\sharp_\gx$ there exists a non-trivial
solution of \eqref{Mq} which vanishes on $\bdw\sms\{\gx\}$ but for
$q> q^\sharp_\gx$ no such solution exists.

In the case of smooth domains, $q_\gx=q^\sharp_\gx$ and
$q_\gx=(N+1)/(N-1)$ for every boundary point $\gx$. Furthermore, if
$q=q_\gx$ there is no solution with isolated singularity at $\gx$,
i.e.,  an isolated singularity at $\gx$ is  removable.

In \Lip domains {\em the critical value depends on the point}.
Clearly $q_\gx\leq q^\sharp_\gx$, but the question whether, in
general, $q_\gx=q^\sharp_\gx$ remains open. However we prove that,
if $\Gw$ is a polyhedron, $q_\gx=q^\sharp_\gx$ at every point and
the function $\gx\to q_\gx$ obtains only a finite number of values.
In fact it is constant on each open face and each open edge, of any
dimension. In addition, if $q=q_\gx$, an isolated singularity at
$\gx$ is removable. The same holds true in a   piecewise $C^2$
domain $\Gw$ except that $\gx\to q_\gx$ is not constant on edges but
it is continuous on every relatively open edge.

For general \Lip domains, we can provide only a partial answer to
the question posed above.

We say that $\Gw$ possesses a {\em tangent cone} at a point $\gx\in
\bdw$ if the limiting inner cone with vertex at $\gx$ is the same as
the limiting outer cone at $\gx$.

\bth{pref-q=q'} Suppose that  $\Gw$ possesses a tangent cone
$C_\gx^\Gw$ at a point $\gx\in\bdw$ and denote by $q_{c,\gx}$  the
critical value for this cone at the vertex $\gx$.  Then
$$q_\gx=q^\sharp_\gx= q_{c,\gx}.$$
Furthermore,  if $1<q<q_\gx$ then  $\gd_\gx$ is admissible,
i.e.,
$$M_\gx:=\int_\Gw K(x,\gx)^q\gr(x)dx<\infty.$$\es

We do not know if, under the assumptions of this theorem, an
isolated singularity at $\gx$ is removable when $q=q_{c,\gx}$. It
would be useful to resolve this question.

\medskip

\noindent{\bf E.} {\em The generalized \bvp in \Lip domains: the
subcritical case.}

\medskip
In the case of smooth domains, a \bvp for equation \eqref{Mq} is
either subcritical or supercritical. This is no longer the case when
the domain is merely \Lip since the criticality varies from point to
point. In this part of the paper we discuss the generalized \bvp in
the strictly subcritical case. Later we discuss the mixed case
(partly subcritical and partly supercritical) when $\Gw$ is a
polyhedron and the boundary data is given by a bounded measure.

Under the conditions of \rth{pref-q=q'} we know that, if
$\gx\in \bdw$ and $1<q <q_\gx$ then $K(\cdot,\gx)\in
L^1_\gr(\Gw)$. In the next result, we derive, under an
additional restriction on $q$, {\em uniform} estimates of the
norm $\norm{K(\cdot,\gx)}\indx{L^1_\gr(\Gw)}$. Such estimates
are needed in the study of existence and uniqueness. For its
statement we need the following notation:

If $z\in\bdw$, we denote by $S_{z,r}$ the opening of the largest
cone $C_S$ with vertex at $z$ \sth $C_S\cap B_r(z)\sbs\Gw\cup\{z\}$.
If $E$ is a compact subset of $\bdw$ we denote:
$$q^*_E=\lim_{r\to 0}\inf\left\{q\indx{S_{z,r}}:z\in\prt\Gw,
\;\dist(z,E)<r\right\}.$$ We observe that
$$q^*_E\leq \inf\{q_{c,z}:z\in E\}$$
but this number also measures, in a sense, the rate of convergence
of interior cones to the limiting cones. If $\Gw$ is convex then
$q^*_E\leq(N+1)/(N-1)$ for every non-empty set $E$. On the other
hand if $\Gw$ is the {\em complement} of a bounded convex set then
$q^*_E=(N+1)/(N-1)$.

 \bth{pref-admi} If $E$ is a compact subset of $\bdw$ and
 $1<q<q^*_E$ then, there exists $M>0$ such that,
\begin{equation}\label{pref-adm1}
\myint{\Gw}{}K^q(x,y)\gr(x)dx\leq M \forevery y\in E.
\end {equation}
\es

Using this theorem we obtain,

\bth{pref-sub} Assume that $\Gw$ is a bounded Lipschitz domain and
$u$ is a positive solution of \eqref{Mq}. If $y\in\CS(u)$ (i.e.
$y\in \bdw$ is a singular point of $u$) and $1<q<q^*_{\{y\}}$ then,
for every $k>0$, the measure $k\gd_y$ is admissible and
\begin{equation}\label{pref-u>uk}
  u\geq u_{k\gd_{y}}=\text{ solution with boundary trace $k\gd_y$.}
\end{equation}
\es

\Remark It can be shown that, if $q>q^*_{\{y\}}$,
\eqref{pref-u>uk} may not hold. For instance, such solutions
exist if $\Gw$ is a smooth, obtuse cone and $y$ is the vertex
of the cone. Therefore the condition $q<q^*_{\{y\}}$ for every
$y\in \bdw$ is, in some sense necessary for uniqueness in the
subcritical case.

As a consequence we first obtain the existence and uniqueness result
in the context of bounded measures.

 \bth{pref-ex-sub} Let $E\sbs \bdw$ be a closed set and assume that $1<q<q^*_E$. Then,
for every $\mu\in \GTM(\Gw)$ \sth $\supp\mu\sbs E$ there exists a
(unique) solution $u_\mu$ of \eqref{eq-q} in $\Gw$ with boundary
trace $\mu$.
\es

Further, using Theorems \ref{t:pref-reg+sing}, \ref{t:pref-cone1}
and \ref{t:pref-admi}, we establish the existence and uniqueness
result for generalized \bvp{}s.

\bth{pref-wellposed} Let $\Gw$ be a bounded Lipschitz domain which
possesses a tangent cone at every boundary point. If
$$1<q<q^*_{\prt\Gw}$$
then, for every positive, outer regular Borel measure $\bar\gn$ on
$\prt\Gw$, there exists a unique solution $u$ of (\ref{Mq}) such
that $tr_{_{\prt\Gw}}(u)=\bar\gn$. \es

\medskip
{\bf F.} {\em On the action of Poisson type kernels with fractional
dimension.}

In preparation for the study of supercritical \bvp{}s
 we establish an
harmonic analytic result, extending a well known result on the
action of Poisson kernels on Besov spaces with negative index
(see \cite[1.14.4.]{Tri} and \cite{MV5}). We first quote the
classical result for comparison purposes.

\bprop {pref-repr}Let $1<q<\infty$ and $s>0$. Then, for any bounded
Radon measure $\gm$ in $\BBR^{n-1}$,
\begin{equation}\label{pref-P4}
I(\gm)=\int_
{\BBR_+^{n}}\abs{\BBK_{n}[\gm](y)}^qe^{-y_{1}}y_{1}^{sq-1}dy\approx
\norm{\gm}^q_{B^{-s,q}(\BBR^{n-1})}.
\end{equation}
\es

Here $\BBK_n[\mu]$ denotes the Poisson potential of $\mu$ in
$\BBR^n=\BBR_+\ti\BBR^{n-1}$, namely,
\begin{equation}\label{pref-BBKn}
  \BBK_n[\mu](y)=\gg_{n}y_{1}\int_{\BBR^{n-1}}\myfrac{d\gm(z)}{\left(y_1^2+|\gz-z|^2\right)^{n/2}}\forevery y=(y_1,\gz)\in
\BBR^n_{+}
\end{equation}
where $\gg_n$ is a constant depending only on $n$.
 The notation $I\approx J$ means that $c^{-1}I\leq J\leq cI$ for some $c>0$.

 In this paper we prove,

 \bth{pref-general-nu} Let $1<q$,  $m$ a positive integer and $\nu\in \BBR$
 \sth $m+1\leq \nu$.
 For every $s\in (0,m/q')$  there exists a
positive constant $c$ \sth, for every bounded positive measure
$\mu$ supported in $\BBR^m\cap B_{R/2}(0)$, $R>1$,
\begin{equation}\label{pref-Mns<}\BAL
\rec{c}\norm{\gm}^q_{B^{-s,q}(\BBR^{m})}&\leq
\int_0^R\Big(\int_{|\gz|<R}\abs{\BBK_{\nu,m}[\mu](\tau,\gz)}^q
d\gz\Big) \tau^{ sq-1}\,d\tau\\
& \leq c R^{(s+\nu-m) q+1}\norm{\gm}^q_{B^{-s,q}(\BBR^{m})}.\EAL
\end{equation}
Here
\begin{equation}\label{pref-Knum}
 \BBK_{\nu,m}[\mu](\tau,\gz)=\int_{\BBR^{m}}
\myfrac{\tau^{\nu-m}d\gm(z)}{(\tau^2+|\gz-z|^2)^{\gn/2}}\forevery
\tau\in [0,\infty).\end{equation}

 This also holds when $s=m/q'$, provided that the diameter of
$\supp \mu$ is sufficiently small. \es

This is proved in Section 7 (see \rth{general-nu}) using a slightly
different notation. Note that
$$\BBK_n[\mu]=\gg_n\BBK_{n,n-1}[\mu].$$

\medskip
\nind{\bf G.} {\em The admissibility condition and the critical
value in a k-wedge.}

\medskip

The next step towards the study of \bvp{}s in a polyhedron is the
treatment of such problems in  a k-wedge (or k-dihedron) i.e., the
domain defined by the intersection of $k$ hyperplanes in $\BBR^N$,
$1<k<N$. The edge is an $(N-k)$ dimensional space. We note that the
case $k=N$ (which corresponds to a cone) has been studied previously
in this paper while the case $N=1$ (i.e. a half space) is classical.

We denote by $D_A$ a k-wedge \sth, its edge $d_A$ is identified with
$\BBR^{N-k}$ and the 'opening' of the wedge is $A=D_A\cap S^{k-1}$.
If  $S_A$ denotes the spherical domain
\begin{equation}\label{pref-dom}
S_A=\{x\in \BBR^N:\abs x=1,\,x\in A\ti\prod_{j=k}^{N-1}[0,\gp]\}\sbs
S^{N-1}\}
\end{equation}
then
$$D_{A}=\{x=(r,\gs):r>0,\gs\in S_A\},\q D_{A,R}=D_A\cap \Gg_R$$
where $$\Gg_R=\{x=(x',x'')\in
\BBR^k\ti\BBR^{N-k}:\,|x'|<r,\;|x''|<R\}.$$

Let $\gl_{A}$ be the first eigenvalue of $-\Gd_{_{S^{N-1}}}$ in
$W^{1,2}_{0}(S_A)$
 and denote
   \begin{equation}\label{pref-kappa2}\BAL
   \gk_{+}&=\myfrac{1}{2}\left(2-N+\sqrt{(N-2)^2+4\gl_{A}}\right)\\
   \gk_{-}&=\myfrac{1}{2}\left(2-N-\sqrt{(N-2)^2+4\gl_{A}}\right).
\EAL\end{equation}

One can show that the Martin kernel $K_A$ in $D_A$ relative to
points $z\in d_A$ is given by

\begin{equation}\label{pref-k4}
K_{A}(x,z)=c_{_{A}}\myfrac{|x'|^{\gk_{+}}\gw^{\{N-k+1\}}(\gs_{N-k+1})}{|x-z|^{(N-2+2\gk_{+})}},
\end{equation} where $\gw^{\{N-k+1\}}$ is a related eigenfunction
in $A$ and $x=(x',x'')\in \BBR^k\ti\BBR^{N-k}$. Using this formula
we obtain the admissibility condition for a measure $\mu\in
\GTM(d_A)$ \sth $\supp \mu\sbs B_R(0)$:

\begin{equation}\label{pref-k8}
\int_{\Gg_R} \Big(\int_{\BBR^{N-k}}
\frac{|x'|^{\gk_+}d|\gm|(z)}{(|x'|^2+|x''|^2)^{(N-2+2\gk_+)/2}}\Big)^q
|x'|^{\gk_{+}}dx<\infty
\end{equation}
where $\Gg_R:=\{x=(x',x'')\in
\BBR^k\ti\BBR^{N-k}:|x'|<R,\,|x''|<R\}.$

Using this expression we show that the condition
\begin{equation}\label{pref-q-critk}
1<q<q_{c}:=\myfrac{\gk_{+}+N}{\gk_{+}+N-2}.
\end{equation}
is necessary and sufficient in order that the Dirac measure
$\mu=\gd_P$, supported at a point $P\in d_A$, be admissible.

In addition we show that  the condition
\begin{equation}\label {pref-qk}
1<q<
q^*_c:=1+\myfrac{2-k+\sqrt{(k-2)^2+4\gl_{A}-4(N-k)\gk_{+}}}{\gl_{A}-(N-k)\gk_{+}}
\end{equation}
is necessary and sufficient for the existence of a non-trivial
solution $u$ of (\ref{Mq}) in $D_{A}$ which vanishes on $\prt
D_{A}\sms d_A$. Furthermore, when this condition holds, there exist
admissible non-trivial positive {\em bounded measures} $\mu$ on
$d_A$, i.e., measures \sth $\BBK[\gm]\in L^q_{\gr}(\Gg_{R}\cap
D_{A})$.

Finally we have the following removability result:

\bth{pref-edge-bvp} Assume that $q_c\leq q<q^*_c$. A measure $\mu\in
\GTM(\prt D_A)$, with compact support contained in $d_A$, is  good
relative to \eqref{Mq} in $D_A$ if and only if $\mu$ vanishes on
every Borel set $E\sbs d_A$ \sth $C_{s,q'}^{N-k}(E)=0$, where
$s=2-\frac{k+\kappa_+}{q'}$ and $C_{s,q'}^{N-k}$ is the Bessel
capacity with the indicated indices in $\BBR^{N-k}$  (which we
identify with the edge $d_A$). \es

Recall that $\mu$ is 'good' if the specified equation possesses a
solution with boundary data $\mu$. The above result implies in
particular that sets with $C_{s,q'}^{N-k}$-capacity zero are
conditionally removable. However we obtain a much stronger result
later on.


\medskip
\nind{\bf H.} {\em Boundary value problems in a polyhedron: the
supercritical case.}

\medskip

In the final part of the paper (Sections 8) we study \bvp{}s in the
supercritical case in polyhedrons, with trace given by bounded
measures. For such domains $\Gw$ we provide a complete
characterization of 'good measures', i.e., measures $\mu$ on $\bdw$
\sth \eqref{Mq} possesses a (unique) solution with boundary trace
$\mu$. We also provide a complete characterization of {\em removable
sets}. These results, with rather obvious modifications, also apply
to piecewise $C^2$ domains. The case of general \Lip domains and
boundary trace given by unbounded Borel measures will be treated in
a subsequent paper.

\bth{pref-poly-good} Let $\Gw$ be an $N$-dimensional polyhedron. Let
$L$ denote one of the faces, or edges, or vertices of $\Gw$ and let
$Q_L$ denote the half space with boundary $L$, or the wedge with
edge $L$, or the cone with vertex $L$ \sth $\Gw\sbs D_L$ and $\prt
Q_L$ is determined by the faces of $\Gw$ adjacent to $L$. Thus
$\bdw$ is the union of the sets $\prt\Gw\cap\prt Q_L$.

Denote by $A_L$ the opening of $Q_L$ so that, in the notation of
{\bf G}, $Q_L=D_{A_L}$ and denote by $\gk_+(L)$, $q_c(L)$ etc. the
various notations introduced in {\bf G} relative to $A_L$. In
particular let $k(L)$ denote the co-dimension of the linear space
spanned by $L$ and put
$$s(L)=2-\frac{k(L)+\kappa_+(L)}{q'}.$$

Let $\mu$ be a bounded measure on $\prt \Gw$, \(possibly a signed
measure\). Then  $\mu$ is a good measure relative to \eqref{Mq} in
$\Gw$, if and only if, for every $L$ as above and every Borel set
$E\sbs L$ the following condition holds.

If $1\leq k=\mathrm{codim}L<N$ then
\begin{equation}\label{pref-poly-good1}\BAL
  C^{N-k}_{s(L),q'}(E)=0 &\Lra
 \mu(E)=0 &&\txt{if} q_c(L)\leq q<q^*_c(L)\\
 q\geq q_c^*(L) &\Lra
 \mu(L)=0 &&\txt{if} q\geq q^*_c(L)
\EAL\end{equation}  and if $k=N$ (i.e., $L$ is a vertex)
\begin{equation}\label{pref-poly-good2}
q\geq
q_c(L)=\frac{N+2+\sqrt{(N-2)^2+4\gl_{A_L}}}{N-2+\sqrt{(N-2)^2+4\gl_{A_L}}}\Lra
\mu(L)=0.
\end{equation}

In all cases, if $1<q<q_c(L)$ then there is no restriction on
$\mu\chi\indx{L}$.
\es

\medskip

\nind{\bf I.} {\em Characterization of removable sets.}

Let $\Gw$ be an N-dimensional polyhedron.

\rth{pref-poly-good} provides a necessary and sufficient
condition for the removability of a singular set $E$ relative
to the family of solutions $u$ \sth
$$\int_\Gw |u|^q\gr\,dx<\infty.$$
 The next
result provides a necessary and sufficient condition for
removability in the full sense, as  defined in
\rdef{pref-maxsol}.

\bth{pref-rem} Let $\Gw$ be an N-dimensional polyhedron and let
$E$ be a compact subset of $\bdw$. A nonempty compact set
$E\sbs\bdw$ is removable if and only if, for every $L$ as in
{\bf G} \sth $1\leq
k=\mathrm{codim}L<N$ the following condition holds:\\
either

$$q_c(L)\leq q<q^*_c(L)\txt{and} C^{N-k}_{s(L),q'}(E)=0$$
or $q\geq q^*_c(L)$. In the case $k=N$ the condition is $q\geq
q^*_c(L)=q_c(L)$. \es

\medskip

\nind\/{\bf Aknowledgements}. The authors are grateful to the High
Council for  Scientific and Technological Cooperation between France
and Israel for its financial support.

\mysection{Boundary value problems }
\subsection{Classical harmonic analysis in Lipschitz domains}
A bounded domain $\Gw\subset \BBR^N$ is called a Lipschitz domain if there exist positive numbers $r_0,\gl_0$
and a cylinder
\begin{equation}\label{Or0}
 O_{r_0}=\set{\gx=(\gx_1,\gx')\in\RN: |\gx'|<r_0,\;|\gx_1|<r_0}
\end{equation}
such that, for every $y\in \bdw$ there exist:\1
(i) A Lipschitz function $\psi^y$ on the $(N-1)$-dimensional ball $B'_{r_0}(0)$ with Lipschitz constant
$\geq \gl_0$;\1
(ii) An isometry $T^y$ of $\RN$ \sth
\begin{equation}\label{Lipdom}
\BAL
&T^y(y)=0,\q (T^y)^{-1}(O_{r_0}):=O^y_{r_0},\\
&T^y(\bdw\cap O^y_{r_0})=\set{(\psi^y(\gx'),\gx'):\, \gx'\in B'_{r_0}(0)}\\
&T^y(\Gw\cap O^y_{r_0})=\set{(\gx_1,\gx'):\, \gx'\in B'_{r_0}(0),\;-r_0<\gx_1<\psi^y(\gx')}
\EAL
\end{equation}

The constant $r_0$ is called a {\em localization constant} of $\Gw$; $\gl_0$ is called a
{\em Lipschitz constant} of $\Gw$.
The pair $(r_0,\gl_0)$ is called a {\em Lipshitz character} (or, briefly, L-character) of $\Gw$. Note that,
if $\Gw$ has L-character $(r_0,\gl_0)$ and $r'\in(0,r_0)$, $\gl'\in (\gl_0,\infty)$ then $(r',\gl')$
is also an L-character of $\Gw$.

By the Rademacher theorem, the outward normal  unit vector  exists $\CH^{N-1}-$a.e. on $\prt\Gw$,
where $\CH^{N-1}$ is the N-1 dimensional Hausdorff measure. The unit normal at a point $y\in \bdw$ will
be denoted by $\bmn_y$.

\medskip
 We list below some facts concerning the Dirichlet problem in Lipschitz domains.

 \medskip
\noindent{\bf A.1}- Let $x_0\in \Gw$, $h\in C(\prt\Gw)$ and denote $L_{x_0}(h):=v_h(x_0)$ where $v_h$ is
 the solution of the Dirichlet problem
\begin{equation}\label{D1}\left\{\BA {l}
-\Gd v=0\quad\in \Gw\\
\phantom{-\Gd}
v=h\quad\text{on }\prt\Gw.
\EA\right.
\end{equation}
Then $L_{x_0}$ is a continuous linear functional on $C(\bdw)$. Therefore there exists a unique Borel
measure on $\prt\Gw$, called the harmonic measure in $\Gw$,  denoted
by $\gw^{x_0}_\Gw$ \sth
\begin{equation}\label{D2}
 v_h(x_0)=\int_{\bdw}hd\gw^{x_0}_\Gw \forevery h\in C(\prt\Gw).
\end{equation}
When there is no danger of confusion, the subscript $\Gw$ will be dropped. Because of Harnack's
inequality the measures $\gw^{x_0}$ and $\gw^{x}$, $x_0,\,x\in \Gw$ are mutually absolutely continuous.
For every fixed $x\in \Gw$ denote the Radon-Nikodym derivative by
\begin{equation}\label{HM1}
K(x,y):=\myfrac{d\gw^{x}}{d\gw^{x_0}}(y) \q\text{for $\gw^{x_0}$-a.e. } y\in\prt\Gw.
\end{equation}
Then, for every $\bar x\in \Gw$,  the function $y\mapsto K(\bar x,y)$ is positive and continuous on $\prt\Gw$ and,
for every  $\bar y\in \bdw$, the function $x\mapsto K(x,\bar y)$ is harmonic in $\Gw$ and satisfies
$$\lim_{x\to y}K(x,\bar y)=0\forevery y\in \prt\Gw\setminus\{\bar y\}.$$

By \cite{HW}
\begin{equation}\label{mart}\lim_{z\to y}\myfrac{G(x,z)}{G(x_0,z)}=K(x,y)\forevery y\in\prt\Gw
\end{equation}
Thus the kernel $K$ defined above is the {\em Martin kernel}.

The following is an equivalent definition of the harmonic measure \cite{HW}:\\
 For any closed set $E\sbs \prt\Gw$
\begin{equation}\label{w-measure}\BAL
  &\gw^{x_0}(E):=\\
  &\inf\{\gf(x_0):\,\gf\in C(\Gw)_+\,\text{ superharmonic in $\Gw$,\;}
  \liminf_{x\to E}\gf(x)\geq 1\}.
\EAL\end{equation}
The extension to open sets and then to arbitrary Borel sets is standard.

By \eqref{D2}, \eqref{HM1} and \eqref{w-measure}, the unique solution $v$ of (\ref{D1}) is given by
\begin{equation}\label{HM}\BAL
&v(x)=\myint{\prt\Gw}{}K(x,y)h(y)d\gw^{x_0}(y)=\\&\inf\{\gf\in C(\Gw): \gf\text{ superharmonic,\;}
 \liminf_{x\to y}\gf(x)\geq h(y),\;\forall y\in \bdw\}.&
\EAL\end{equation}
For details see \cite{HW}.\2
\noindent {\bf A.2}- Let $(x_0,y_0)\in\Gw\ti\prt\Gw$. A function $v$ defined
in $\Gw$ is called a kernel function at $y_0$ if it is positive and
harmonic in $\Gw$ and verifies $v(x_0)=1$ and $\lim_{x\to y}v(x)=0$
for any $y\in \prt\Gw\setminus\{y_0\}$. It is proved in \cite[Sec 3]{HW} that the kernel function
at $y_0$ is unique. Clearly this unique function is $K(\cdot, y_0)$. \2
\noindent {\bf A.3}- We denote by $G(x,y)$ the Green kernel for the Laplacian in $\Gw\ti\Gw$.
This means that the solution of the Dirichlet problem
\begin{equation}\label{dir}\left\{\BA {l}
-\Gd u=f\quad\text{in }\Gw,\\
\phantom{-\Gd} u=0\quad\text{on }\prt\Gw,
\EA\right.\end{equation}
with $f\in C^2(\overline\Gw)$, is expressed by
\begin{equation}\label{dir1}
u(x)=\myint{\Gw}{}G(x,y)f(y)dy\forevery y\in\overline\Gw.
\end{equation}
We shall write (\ref{dir1}) as $u=\BBG[f]$.\2
\noindent {\bf A.4}- Let $\Gl$ be the  first eigenvalue of $-\Gd$ in $W^{1,2}_{0}(\Gw)$ and denote by $\gr$ the corresponding eigenfunction
normalized by $\max_\Gw \gr =1$.

 Let $0<\gd<\dist(x_{0},\Gw) $ and put
 $$C_{x_{0},\gd}:=\max_{|x-x_0|=\gd}G(x,x_0)/\gr(x).$$ Since $C_{x_{0},\gd}\,\gr-G(\cdot,x_0)$ is
 superharmonic,  the maximum principle implies that
\begin{equation}\label{GR'}
0\leq G(x,x_{0})\leq C_{x_{0},\gd}\,\gr(x)\quad\forall x\in \Gw\setminus B_{\gd}(x_{0}).
\end{equation}

On the other hand, by
\cite[Lemma 3.4]{KP}: for any $x_{0}\in\Gw$ there exists a constant
$C_{x_{0}}>0$ such that
\begin{equation}\label{GR}
0\leq \gr(x)\leq C_{x_{0}}G(x,x_{0})\quad\forall x\in \Gw.
\end{equation}

\noindent {\bf A.5}- For every
bounded regular Borel measure $\gm$ on $\prt\Gw$ the function
\begin{equation}\label{HM2}
v(x)=\myint{\prt\Gw}{}K(x,y)d\gm(y)\forevery x\in\Gw,
\end{equation}
is harmonic in $\Gw$. We denote this relation by $v=\BBK[\gm]$.\2
\noindent {\bf A.6}- Conversely, for every {\em positive} harmonic function $v$ in $\Gw$
there exists a unique positive
bounded regular Borel measure $\gm$ on $\prt\Gw$ \sth (\ref{HM2}) holds.
The measure $\gm$ is constructed as follows \cite[Th 4.3]{HW}.

Let $SP(\Gw)$ denote the set of continuous, non-negative
superharmonic functions in $\Gw$. Let $v$ be a positive harmonic
function in $\Gw$.

If $E$ denotes a relatively closed subset of
$\Gw$, denote by $R^E_v$ the function defined in $\Gw$ by
$$R^E_v(x)=\inf\set{\gf(x):\,\gf\in SP(\Gw),\; \gf\geq v \text{ in }E}.$$
Then $R^E_v$ is superharmonic in $\Gw$, $R^E_v$ decreases as $E$ decreases and,
if $F$ is another relatively closed subset of $\Gw$, then
$$R^{E\cup F}_v\leq R^E_v+R^F_v.$$
Now, relative to a point $x\in \Gw$, the measure $\gm$ is defined by,
\begin{equation}\label{constr}
\gm^x_v(F)=\inf\{R_v^{E}(x):E=\bar D\cap \Gw,\;D\text { open in }\BBR^N,\;D\supset F\},
\end{equation}
for every compact set $F\sbs \bdw$. From here it is extended to open sets and then to
arbitrary Borel sets in the usual way.

It is easy to see that, if $D$ contains $\bdw$ then $R_v^{\bar D\cap \Gw}=v$. Therefore
\begin{equation}\label{gm(bndOm)}
   \gm^x_v(\bdw)=v(x).
\end{equation}
In addition, if $F$ is a compact subset of the boundary, the function $x\mapsto \gm^x_v(F)$
 is harmonic in $\Gw$ and vanishes on $\bdw\sms F$.\2
\noindent {\bf A.7}- If $x,x_0$ are two points in $\Gw$,  the Harnack inequality  implies that $\gm^x_v$ is absolutely continuous
with respect to $\gm^{x_0}_v$. Therefore, for $\mu^{x_0}_v$-a.e.
 point $y\in \bdw$, the density function $d\gm^x_v/d\gm^{x_0}_v(y)$ is a kernel function at $y$.
By the uniqueness of the kernel
function it follows that
\begin{equation}\label{densityK}
\frac{d\gm^x_v}{d\gm^{x_0}_v}(y)=K(x,y),\quad \gm^{x_0}_v\text{-a.e. } y\in\prt\Gw.
\end{equation}
Therefore, using \eqref{gm(bndOm)},
\begin{equation}\label{constr1}\BAL &(a) \q\gm^x_v(F)=\myint{F}{}K(x,y)d\gm^{x_0}(y),\\
&(b)\q  v(x)=\myint{\prt\Gw}{}K(x,y)d\gm^{x_0}(y).
\EAL\end{equation}
{\bf A.8}- By a result of Dahlberg \cite[Theorem 3]{Da}, the (interior) normal derivative of
$G(\cdot,x_0)$ exists $\CH_{N-1}$-a.e. on $\bdw$ and is positive. In addition,
for every Borel set $E\sbs\bdw$,
\begin{equation}\label{Da-Poisson}
  \gw^{x_0}(E)=\gg_N\int_E \prt G(\gx,x_0)/\prt {\bf n_\gx}\,dS_\gx,
\end{equation}
where $\gg_N (N-2)$ is the surface area of the unit ball in $\BBR^N$ and $dS$ is surface measure on
$\bdw$. Thus, for each fixed $x\in \Gw$, the harmonic measure
 $\gw^{x}$ is absolutely continuous relative to
$\CH_{N-1}\big|_{\bdw}$ with density function $P(x,\cdot)$ given by
\begin{equation}\label{P(x,y)}
    P(x,\gx)=\prt G(\gx,x)/\prt {\bf n_\gx}\txt{for a.e.} \gx\in \bdw.
\end{equation}
In view of \eqref{HM}, the unique solution $v$ of (\ref{D1}) is given by
\begin{equation}\label{HM'}
 v(x)=\int_\Gw\,P(x,\gx)h(\gx)dS_\gx
\end{equation}
for every $h\in C(\bdw)$.  Accordingly $P$ is the {\em Poisson kernel} for $\Gw$.
The expression on the right hand side of \eqref{HM'} will be denoted by $\BBP[h]$. We observe that,
\begin{equation}\label{K,P}
   \BBK[h\gw^{x_0}]=\BBP[h] \forevery h\in C(\bdw).
\end{equation}
\noindent {\bf A.9}- The {\em \BHP}, first proved in \cite{Da}, can be formulated as follows \cite{JK1}.

Let $D$ be a \Lip domain with L-character $(r_0,\gl_0)$. Let $\gx\in
\prt D$ and $\gd\in (0,r_0)$. Assume that $u,v$ are positive
harmonic functions in $D$, vanishing on $\prt D\cap B_{\gd}(\gx)$.
Then there exists a constant $C=C(N,r_0,\gl_0)$ \sth,
\begin{equation}\label{BHP}
  C^{-1}u(x)/v(x)\leq u(y)/v(y)\leq Cu(x)/v(x) \forevery x,y\in B_{\gd/2}(\gx).
\end{equation}

\noindent {\bf A.10}- Let $D,D'$ be two \Lip domains with L-character $(r_0,\gl_0)$. Assume that $D'\sbs D$ and $\prt D\cap\prt D'$
contains a relatively open set $\Gamma$. Let $x_0\in D'$ and let $\gw,\gw'$ denote the harmonic measures
of $D,D'$ respectively, relative to $x_0$. Then, for every compact set $F\sbs \Gamma$, there exists a constant
$c_F=C(F,N,r_0,\gl_0, x_0)$ \sth
\begin{equation}\label{gw-gw'}
   \gw'\lfloor_F\leq \gw\lfloor_F\leq c_F\gw'\lfloor_F.
\end{equation}

Indeed, if $G,G'$ denote the Green functions of $D,D'$ respectively then, by the \BHP,
\begin{equation}\label{dnG-dnG'}
  \prt G'(\gx,x_0)/\prt {\bf n_\gx}\leq \prt G(\gx,x_0)/\prt {\bf n_\gx}
  \leq c_F\prt G(\gx,x_0)/\prt {\bf n_\gx}\txt{for a.e.} \gx\in F.
\end{equation}
Therefore \eqref{gw-gw'} follows from \eqref{Da-Poisson}.

\noindent {\bf A.11}- By \cite[Lemma 3.3]{KP}, for every positive harmonic function $v$ in $\Gw$,
\begin{equation}\label{GR2}
\myint{\Gw}{}v(x)G(x,x_{0})dx<\infty.
\end{equation}
In view of \eqref{GR}, it follows that $v\in L^1_\gr(\Gw)$.
\subsection{The dynamic approach to boundary trace.}
Let $\Gw$ be a bounded Lipschitz domain and $\set{\Gw_n}$ be a {\em
\Lip exhaustion} of $\Gw$. This means that, for every $n$, $\Gw_n$
is \Lip and
\begin{equation}\label{exhaustion}
  \Gw_n\sbs \bar \Gw_n\sbs \Gw_{n+1},\q \Gw=\cup\Gw_n,\q
  \BBH_{N-1}(\bdw_n)\to \BBH_{N-1}(\bdw).
\end{equation}

 \blemma{hm-w-conv} Let $x_0\in \Gw_1$ and denote by $\gw_n$
(respectively $\gw$) the harmonic measure in $\Gw_n$ (respectively
$\Gw$) relative to $x_0$. Then, for every $Z\in C(\bar \Gw)$,
\begin{equation}\label{hm-w-conv}
   \lim_{n\to\infty}\int_{\bdw_n}Z\,d\gw_n=\int_{\bdw}Z\,d\gw.
\end{equation}
\es
\Proof By the definition of harmonic measure
$$\int_{\bdw_n}d\gw_n=1.$$
We extend $\gw_n$ as a Borel measure on $\bar \Gw$ by setting $\gw_n(\bar\Gw\sms \bdw_n)=0$,
and keep the notation $\gw_n$ for the extension. Since the \seq $\set{\gw_n}$ is bounded,
there exists a weakly convergent  \sseq (still denoted by $\set{\gw_n}$).
Evidently the limiting measure, say $\tl\gw$  is supported in $\bdw$ and $\tl\gw(\bdw)=1$.
It follows that for every $Z\in C(\bar\Gw)$,
$$\int_{\bdw_n}Z\,d\gw_n\to \int_{\bdw}Z\,d\tl\gw.$$
Let $\gz:=Z\mid_{\bdw}$ and $z:=\BBK^\Gw[\gz]$. Again by the definition of harmonic measure,
$$\int_{\bdw_n}z\,d\gw_n= \int_{\bdw}\gz\,d\gw=z(x_0).$$
It follows that
$$\int_{\bdw}\gz\,d\tl\gw=\int_{\bdw}\gz\,d\gw,$$
for every $\gz\in C(\bdw)$. \Consy $\tl\gw=\gw$. Since the limit does not depend on the \sseq
it follows that the whole \seq $\set{\gw_n}$ converges weakly to $\gw$. This implies \eqref{hm-w-conv}.
\qed

In the next lemma we continue to use the  notation introduced above.
\blemma{tr-w-conv}
Let $x_0\in \Gw_1$, let $\mu$ be a bounded Borel measure on $\bdw$ and put $v:=\BBK^\Gw[\mu]$.
 Then, for every $Z\in C(\bar \Gw)$,
\begin{equation}\label{tr-w-conv}
   \lim_{n\to\infty}\int_{\bdw_n}Zv\,d\gw_n=\int_{\bdw}Z\,d\mu.
\end{equation}
\es
\Proof It is sufficient to prove the result for positive $\mu$. Let $h_n:=v\mid_{\bdw_n}$.
Evidently $v=\BBK^{\Gw_n}[h_n\gw_n]$ in $\Gw_n$. Therefore
$$v(x_0)=\int_{\bdw_n}h_nd\gw_n=\mu(\bdw).$$
Let $\mu_n$ denote the extension of $h_n\gw_n$ as a measure in $\bar\Gw$ \sth
$\mu_n(\bar\Gw\sms \bdw_n)=0$. Then $\set{\mu_n}$ is bounded and \consy there exists a
weakly convergent  \sseq  $\set{\mu_{n_j}}$.
The limiting measure, say $\tl\mu$,  is supported in $\bdw$ and
\begin{equation}\label{tr-w-conv0}
\tl\mu(\bdw)=v(x_0)=\mu(\bdw).
\end{equation}
It follows that for every $Z\in C(\bar\Gw)$,
$$\int_{\bdw_{n_j}}Z\,d\mu_{n_j}\to \int_{\bdw}Z\,d\tl\mu.$$

To complete the proof, we have to show that $\tl\mu=\mu$. Let $F$ be a closed subset of
$\bdw$ and put,
$$\mu^F=\mu\chi_{F}, \q v^F=\BBK^\Gw[\mu^F].$$
Let $h_n^F:=v^F\mid_{\bdw_n}$ and let $\mu_n^F$ denote the extension of $h_n^F\gw_n$ as a
measure in $\bar\Gw$ \sth $\mu_n^F(\bar\Gw\sms \bdw_n)=0$.
As in the previous part of the proof, there exists a weakly convergent \sseq of
$\set{\mu_{n_j}^F}$.
The limiting measure $\tl \mu^F$ is supported in $F$ and
$$\tl\mu^F(F)=\tl\mu^F(\bdw)=v^F(x_0)=\mu^F(\bdw)=\mu(F).$$
As $v^F\leq v$, we have $\tl\mu^F\leq \tl\mu$. \Consy
\begin{equation}\label{tr-w-conv2}
\mu(F)\leq \tl\mu(F).
\end{equation}
Observe that $\tl{\mu}$ depends on the first \sseq
$\set{\mu_{n_j}}$, but not on the second \sseq. Therefore
\eqref{tr-w-conv2} holds for every closed set $F\sbs \bdw$, which
implies that $\mu\leq \tl\mu$. On the other hand, $\mu$ and $\tl\mu$
are positive measures which, by \eqref{tr-w-conv0}, have the same
total mass. Therefore $\mu=\tl\mu$. \qed

\blemma{K(mu)} Let $\mu\in\GTM(\bdw)$ (= space of bounded Borel measures on $\bdw$).
Then $\BBK[\mu]\in L^1_\gr(\Gw)$ and there exists a constant $C=C(\Gw)$ \sth
\begin{equation}\label{K(mu)<}
   \norm{\BBK[\mu]}_{L^1_\gr(\Gw)}\leq C\norm{\mu}_{\GTM(\bdw)}.
\end{equation}
In particular if $h\in L^1(\bdw;\gw)$ then
\begin{equation}\label{K(h)<}
   \norm{\BBP [h]}_{L^1_\gr(\Gw)}\leq C\norm{h}_{L^1(\bdw;\gw)}.
\end{equation}
\es
\Proof Let $x_0$ be a point in $\Gw$ and let $K$ be defined as in \eqref{HM1}. Put
$\gf(\cdot)=G(\cdot,x_0)$ and $d_0=\dist(x_0,\Gw)$. Let $(r_0,\gl_0)$ denote the Lipschitz character of $\Gw$.

By \cite[Theorem 1]{Bog}, there exist positive constants $c_1(N,r_0,\gl_0,d_0)$ and
$c_0(N,r_0,\gl_0,d_0)$ such that
for every $y\in \bdw$,
\begin{equation}\label{K-estimate}
 c_1^{-1} \frac{\gf(x)}{\gf^2(x')}|x-y|^{2-N}   \leq  K(x,y)\leq c_1 \frac{\gf(x)}{\gf^2(x')}|x-y|^{2-N},
\end{equation}
for all $x,x'\in \Gw$ such that
\begin{equation}\label{BHar1}
c_0|x-y|< \dist(x',\bdw)\leq |x'-y|<|x-y|<\rec{4}\min(d_0,r_0/8).
\end{equation}
Therefore, by \eqref{GR} and \eqref{GR'}, there exists a constant $c_2(N,r_0,\gl_0,d_0)$ \sth
$$ c_2^{-1}\frac{\gf^2(x)}{\gf^2(x')}|x-y|^{2-N}\leq \gr(x)K(x,y)\leq c_2 \frac{\gf^2(x)}{\gf^2(x')}|x-y|^{2-N}$$
for $x,x'$ as above.
There exists a constant $\bar c_0$, depending on $c_0,N$, such that, for every $x\in \Gw$ satisfying
$|x-y|<\rec{4}\min(d_0,r_0/8)$ there exists $x'\in \Gw$ which satisfies \eqref{BHar1} and also
$$|x-x'|\leq\bar c_0\min(\dist(x,\bdw),\dist(x',\bdw)).$$
 By the Harnack chain argument, $\gf(x)/\gf(x')$ is bounded by a constant depending  on $N,\bar c_0$.
 Therefore
\begin{equation}\label{grK}
 c_3^{-1}|x-y|^{2-N}\leq \gr(x)K(x,y)\leq  c_3|x-y|^{2-N}
\end{equation}
for some constant
$c_3(N,r_0,\gl_0,d_0)$ and all $x\in \Gw$ sufficiently close to the boundary.

Assuming that $\mu\geq 0$,
$$\int_\Gw \BBK[\mu](x)\gr(x)dx=\int_{\bdw}\int_\Gw K(x,\gx)\gr(x)dx\, d\mu(\gx)
\leq C\norm{\mu}_{\GTM(\bdw)}.$$
In the general case we apply this estimate to $\mu_+$ and $\mu_-$. This implies \eqref{K(mu)<}.
For the last statement of the theorem see \eqref{K,P}.
\qed

\bprop{IBP} Let $v$ be a positive harmonic function in $\Gw$ with
boundary trace $\gm$. Let  $Z\in C^2(\overline\Gw)$  and let $\tilde
G\in C^\infty(\Gw)$ be a function that coincides with $x\mapsto
G(x,x_{0})$ in $Q\cap \Gw$ for some neighborhood $Q$ of $\prt\Gw$
and some fixed $x_{0}\in\Gw$. In addition assume that there exists a
constant $c>0$ \sth
\begin{equation}\label{ibp0}
|\nabla Z\cdot\nabla \tl G|\leq c\gr.
\end{equation}
Under these assumptions, if $\gz:=Z\tilde G$ then
\begin{equation}\label{ibp1}
-\myint{\Gw}{}v\Gd\gz\, dx=\myint{\prt\Gw}{}Z d \gm.
\end{equation}
\es
\Remark This result is useful in a k-dimensional dihedron in the case where $\mu$ is concentrated on the edge. In such a case one can find, for every smooth function on the edge, a lifting $Z$ \sth condition \eqref{ibp0} holds. See Section 8 for such an application.

\Proof Let $\set{\Gw_n}$ be a $C^1$ exhaustion of $\Gw$. We assume
that $\bdw_n\sbs Q$ for all $n$ and $x_0\in \Gw_1$. Let $\tl G_n(x)$
be a function in $C^1(\Gw_n)$ \sth $\tl G_n$ coincides with
$G^{\Gw_n}(\cdot, x_0)$ in $Q\cap\Gw_n$, $\tl G_n(\cdot,x_0)\to \tl
G(\cdot,x_0)$ in $C^2(\Gw\sms Q)$ and $\tl G_n(\cdot,x_0)\to\tl
G(\cdot,x_0)$ in $\mathrm{Lip}\,(\Gw)$. If $\gz_n=Z\tl G_n$ we have,

$$\BAL-\int_{\Gw_n}v\Gd\gz_n\,dx&=\int_{\bdw_n} v\prt_{\bmn}\gz\, dS=\int_{\bdw_n} vZ\prt_{\bmn}\tl G_n(\gx,x_0)\,dS\\
&=\int_{\bdw_n} vZP^{\Gw_n}(x_0,\gx)\, dS= \int_{\bdw_n} vZ\,d\gw_n.
\EAL$$
By \rlemma{tr-w-conv},
$$\int_{\bdw_n} vZ\,d\gw_n\to\int_{\bdw}Z\,d\mu.$$
On the other hand, in view of \eqref{ibp0}, we have
$$\Gd \gz_n=\tl G_n\Gd Z +Z\Gd\tl G_n + 2\nabla Z\cdot \nabla \tl G_n\to \Gd Z$$
in $L^1_\gr(\Gw)$; therefore,
$$-\int_{\Gw_n}v\Gd\gz_n\,dx\to -\int_{\Gw}v\Gd\gz\,dx.$$
\qed

\bdef{unif-ex} Let $D$ be a \Lip domain and let $\set{D_n}$ be a \Lip  exhaustion of $D$. We say that
$\set{D_n}$ is a {\em uniform \Lip exhaustion} if there exist positive numbers $\bar r, \bar \gl$ \sth
$D_n$ has L-character $(\bar r, \bar \gl)$ for all $n\in\BBN$. The pair $(\bar r, \bar \gl)$ is an L-character
of the exhaustion.
\es

\blemma{D'sbsD} Assume $D,D'$ are two Lipschitz domains such that
$$\Gamma\sbs \prt D\cap\prt D'\sbs\prt(D\cup D')$$
where $\Gamma$ is a relatively open set.
 Suppose $D,D', D\cup D'$ have L-character $(r_0,\gl_0)$. Let $x_0$ be a point
in $D\cap D'$ and put  $$d_0=\min(\dist(x_0,\prt D),\dist(x_0,\prt D')).$$
 Let $u$ be a positive harmonic function in $D\cup D'$ and denote
its boundary trace on $D$ (resp. $D'$) by $\mu$ (resp. $\mu'$).
Then, for every compact set $F\sbs \Gamma$,
 there exists a constant $c_F=c(F,r_0,\gl_0,d_0,N)$ \sth
\begin{equation}\label{mu--mu'}
  c_F^{-1}\mu'\lfloor_F\leq \mu\lfloor_F\leq c_F\mu'\lfloor_F.
\end{equation}
\es

\Proof We  prove \eqref{mu--mu'} in the case that $D'\sbs D$.
This implies \eqref{mu--mu'} in the general case by comparison of the boundary trace on $\prt D$ or $\prt D'$
with the boundary trace on $\prt (D\cup D')$.

Let $Q$ be an open set \sth $Q\cap D$ is Lipschitz and
 $$F\sbs Q, \q \bar Q\cap D\sbs D',\q \bar Q\cap\prt D\sbs\Gamma.$$
Then there exist uniform  \Lip exhaustions of $D$ and $D'$, say $\set{D_n}$ and $\set{D'_n}$,
possessing the following properties:

\smallskip
(i) $\bar D'_n\cap Q=\bar D_n\cap Q$.

(ii) $x_0\in D'_1$ and $\dist(x_0,\prt D'_1)\geq \rec{4}d_0$.

(iii) There exist $r_Q>0$ and $\gl_Q>0$ \sth both exhaustions have L-character $(r_Q,\gl_Q)$.

\smallskip

Put $\Gamma_n:=\prt D_n\cap Q=\prt D'_n\cap Q$ and let $\gw_n$ (resp. $\gw'_n$) denote
the harmonic measure, relative to $x_0$,
of $D_n$ (resp. $D'_n$). By \rlemma{tr-w-conv},
$$\int_{\Gamma_n} \phi\,u(y)\,d\gw_n(y)\to \int_{\Gamma} \phi\,d\mu,$$
and
$$\int_{\Gamma_n} \phi\,u(y)d\gw'_n(y)\to \int_{\Gamma} \phi\,d\mu'$$
for every $\phi\in C_c(Q)$. By {\bf A.10} there exists a constant $c_Q=c(Q,r_Q,\gl_Q,d_0,N)$ \sth
$$\gw'_n\lfloor_{\Gamma_n}\leq \gw_n\lfloor_{\Gamma_n}\leq c_Q\gw'_n\lfloor_{\Gamma_n}.$$
This implies \eqref{mu--mu'}.
\qed

\subsection{$L^1$ data}
We denote by $X(\Gw)$ the space of test functions,
\begin{equation}\label {L1}
X(\Gw)=\left\{\eta\in W^{1,2}_{0}(\Gw):\gr^{-1}\Gd \eta\in L^{\infty}(\Gw)\right\}.
\end{equation}
Let $X_{+}(\Gw)$ denote its positive cone.

Let $f\in L^\infty(\Gw)$, and let $u$ be the weak $W^{1,2}_0$
solution of the Dirichlet problem
\begin{equation}\label{Gdu=f}
  -\Gd u=f\txt{in $\Gw$},\q u=0\txt{on $\bdw$}
\end{equation}
If $\Gw$ is a Lipschitz domain (as we assume here) then $u\in C(\bar
\Gw)$ (see \cite{Tru}). Since $\BBG[f]$ is a weak $W^{1,2}_0$
solution, it follows that the  solution
 of \eqref{Gdu=f}, which is unique in $C(\bar\Gw)$,  is given by $u=\BBG[f]$. If, in addition,
 $|f|\leq c_1\gr$ then, by the maximum principle,
 \begin{equation}\label{ugr-bound}
   |u|\leq (c_1/\Gl)\gr,
 \end{equation}
 where $\Gl$ is the first eigenvalue of $-\Gd$ in $\Gw$.

In particular, if  $\eta\in X(\Gw)$ then $\eta\in C(\bar\Gw)$ and it
satisfies
\begin{align}\label{rep2}
-\BBG[\Gd \eta]&=\eta,\\
|\eta|&\leq \Gl^{-1}\norm{\gr^{-1}\Gd\eta}_{L^\infty}\gr\label{X}.
\end{align}

If, in addition, $\Gw$ is a $C^2$ domain then the solution of
\eqref{Gdu=f} is in $C^1(\bar\Gw)$.

\blemma{equivlin}Let $\Gw$ be a Lipschitz bounded domain. Then for any $f\in L^1_{\gr}(\Gw)$
there exists a unique $u\in L^1_{\gr}(\Gw)$ such that
\begin{equation}\label {L'4}
-\myint{\Gw}{}u\Gd\eta\,dx=\myint{\Gw}{}f\eta dx\forevery\eta \in X(\Gw).
\end{equation}
Furthermore $u=\BBG[f]$.
Conversely, if $f\in L^1_{loc}(\Gw)$, $f\geq 0$  and there exists $x_0\in\Gw$ \sth $\BBG[f](x_0)<\infty$ then
$f\in L^1_{\gr}(\Gw)$.
Finally
\begin{equation}\label{ulro}
   \norm{u}_{L_\gr(\Gw)}\leq \Gl^{-1} \norm{f}_{L_\gr(\Gw)}
\end{equation}
\es
\Proof First assume that $f$ is bounded. We have already observed  that, in this case, the weak
$W^{1,2}_0$ solution $u$ of the Dirichlet problem \eqref{Gdu=f} is in $C(\bar \Gw)$ and $u=\BBG[f]$.
Furthermore, it follows from \cite {BB} that
$$\int_{\Gw}\nabla \eta\cdot\nabla u dx
=-\int_{\Gw}u\Gd\eta dx.
$$
Thus $u=\BBG[f]$ is also a weak $L^1_\gr$ solution (in the sense of
\eqref{L'4}).

 Let $\eta_0$ be the weak
$W^{1,2}_0$ solution of \eqref{Gdu=f} when $f={\rm sgn}(u)\gr$; evidently $\eta_0\in X(\Gw)$.
If $u\in L^1_{\gr}(\Gw)$ is a solution of \eqref{L'4} for some $f\in L^1_{\gr}(\Gw)$ then
\begin{equation}\label{new1}
\int_{\Gw}|u|\gr dx=\int_{\Gw}f\eta_{0}dx\leq \Gl^{-1}\int_{\Gw}|f|\gr dx.
\end{equation}
The second inequality follows from \eqref{ugr-bound}. This proves \eqref{ulro} and implies uniqueness.

Now assume that $f\in L^1_{\gr}(\Gw)$ and let $\set{f_n}$ be a \seq of bounded functions \sth $f_n\to f$
in this space. Let $u_n$ be the weak
$W^{1,2}_0$ solution of  \eqref{Gdu=f} with $f$ replaced by $f_n$. Then $u_n$ satisfies
\eqref{L'4} and $u_n=\BBG[f_n]$. By \eqref{ulro},
 $\set{u_n}$ converges in $L^1_\gr(\Gw)$, say $u_n\to u$. In view of \eqref{GR'} it follows that
 $u=\BBG[f]$ and that $u$ satisfies \eqref{L'4}.

If $f\in L^1_{loc}(\Gw)$, $f\geq 0$  and $\BBG[f](x_0)<\infty$  then, by \eqref{GR},
$f\in L^1_{\gr}(\Gw)$.
\qed
\blemma{lin} Let $\Gw$ be a Lipschitz bounded domain. If  $f\in L^1_{\gr}(\Gw)$ and
$h\in  L^1(\prt\Gw;\gw)$, there exists a unique $u\in L^1_{\gr}(\Gw)$ satisfying
\begin{equation}\label {L4}
\int_{\Gw}{}\left(-u\Gd\eta-f\eta\right)dx=-\int_{\Gw}{}\BBP[h]\Gd\eta dx\forevery\eta \in X(\Gw)
\end{equation}
or equivalently
\begin{equation}\label{u=Gf-Kh}
  u=\BBG[f]-\BBP[h].
\end{equation}
The following estimate holds
\begin{align}\label {L5}
\norm u_{L^1_{\gr}(\Gw)}&\leq c\left(\norm f_{L^1_{\gr}(\Gw)}+
\norm {\BBP[h]}_{L^1_{\gr}(\Gw)}\right)\\
&\leq c\left(\norm f_{L^1_{\gr}(\Gw)}+
\norm {h}_{L^1(\bdw,\gw)}\right).\notag
\end{align}
Furthermore, for any nonnegative element $\eta\in X(\Gw)$, we have
\begin{equation}\label {L5'}
-\myint{\Gw}{}\abs u\Gd\eta\,dx\leq -\myint{\Gw}{}\BBP[\abs {h}]\Gd\eta dx+\myint{\Gw}{}\eta f{\rm sgn}(u)\,dx,
\end{equation}
and
\begin{equation}\label {L5+}
-\myint{\Gw}{}u_{+}\Gd\eta\,dx\leq -\myint{\Gw}{}\BBP[h_{+}]\Gd\eta dx+\myint{\Gw}{}\eta f{\rm sgn}_{+}(u)\,dx.
\end{equation}
\es
\Proof {\it  Existence.}\hskip 2mm By \rlemma{K(mu)}, the assumption on $h$ implies that
$\BBP[|h|]\in L^1_{\gr}(\Gw)$. If we denote by $v$ the
unique function in $L^1_{\gr}(\Gw)$ which satifies
$$-\myint{\Gw}{}v\Gd\eta dx=-\myint{\Gw}{}f\eta dx
\forevery \eta\in X(\Gw),$$
then $u=v-\BBP[h]\in L^1_{\gr}(\Gw)$ and (\ref{L4}) holds.

By \rlemma{equivlin}, \eqref{u=Gf-Kh} is equivalent to \eqref{L4}.\1
\noindent{\it Estimate (\ref{L5})} This inequality follows from \eqref{L4} and
\eqref{ulro}.\1
\noindent {\it Estimate (\ref{L5+}).}\hskip 2mm Let $\set{\Gw_n}$ be
an exhaustion of $\Gw$ by {\em smooth domains}. If $u$ is the
solution of \eqref{L4} and $h_n:=u\big |_{\bdw_n}$ then, in $\Gw_n$,
$$u=\BBG^{\Gw_n}[f]-\BBP^{\Gw_n}[h_n]\txt{in $\Gw_n$,}$$
or equivalently,
\begin{align}\label{temp2.1}
\int_{\Gw_n}{}\left(-u\Gd\eta-f\eta\right)dx&=-\int_{\Gw_n}{}\BBP[h_n]\Gd\eta dx\\
&=-\int_{\bdw_n}(\prt\eta/\prt{\bf n})h_ndx \forevery\eta \in
X(\Gw_n). \notag
\end{align}
We recall that, since $\Gw_n$ is smooth, $\eta\in X(\Gw_n)$ implies
that $\eta\in C^1(\bar \Gw_n)$. In addition it is known that (see
e.g. \cite{Ve1}), for every non-negative $\eta \in X(\Gw_n)$,
\begin{align}\label{temp2.2}
\int_{\Gw_n}{}\left(-|u|\Gd\eta-f\eta\,{\rm sign}\, u\right)dx&\leq-\int_{\bdw_n}\prt\eta/\prt{\bf n}|h_n|dx
\end{align}
Let $\gr_n$ be the first eigenfunction of $-\Gd$ in $\Gw_n$, normalized by $\gr_n(\bar x)=1$ for some $\bar x\in \Gw_1$.
Let $\eta$ be a non-negative function in $X(\Gw)$ and let $\eta_n$ be the solution of the problem
$$\Gd z=(\Gd \eta)\gr_n/\gr \txt{in $\Gw_n$,} z=0 \txt{on $\bdw_n$.}$$
Then $\eta_n\in X(\Gw_n)$ and,
since $\gr_n\to\gr$,
$$\Gd\eta_n\to\Gd\eta,\q \eta_n\to\eta.$$
If $v:=\BBP[|h|]$ then $v\geq|u|$ so that
$$\tl h_n:=v\big|_{\bdw_n}\geq |h_n|.$$
 Therefore
\begin{align}\label{temp2.3}
&-\int_{\bdw_n}\prt\eta_n/\prt{\bf n}|h_n|dx \leq -\int_{\bdw_n}\prt\eta/\prt{\bf n}|\tl h_n|dx=\\
&-\int_{\Gw_n}{}\BBP^{\Gw_n}[\tl h_n]\Gd\eta_n dx=-\int_{\Gw_n}{}v\Gd\eta_n dx
\to-\int_{\Gw}{}v\Gd\eta dx.\notag
\end{align}
Finally, \eqref{temp2.2} and \eqref{temp2.3} imply \eqref{L5'}.\1
{\it Estimate  \eqref{L5+}} This inequality  is obtained
by adding \eqref{L4} and \eqref{L5'}.
\qed
\bdef{class} We shall say that a function $g:\BBR\mapsto\BBR$ belongs to $\CG(\BBR)$
if it is continuous, nondecreasing and $g(0)=0$.
\es
\blemma{Sol} Let $\Gw$ be a Lipschitz bounded domain and $g\in \CG(\BBR)$.
If  $f\in L^1_{\gr}(\Gw)$ and
$h\in  L^1(\prt\Gw;\gw)$, there exists a unique $u\in L^1_{\gr}(\Gw)$ \sth $g(u)\in L^1_{\gr}(\Gw)$
and
\begin{equation}\label{nln}
\myint{\Gw}{}\left(-u\Gd\eta+(g(u)-f)\eta\right)dx=-\myint{\Gw}{}
\BBP[h]\Gd \eta\,dx\forevery\eta\in X(\Gw).
\end{equation}
The correspondence $(f,h)\mapsto u$ is increasing.

If $u,u'$ are solutions of \eqref{nln} corresponding
to data $f,h$ and $f',h'$ respectively then the following estimate holds:
\begin{align}\label {NL5}
&\norm {u-u'}_{L^1_{\gr}(\Gw)}+\norm{ g(u)-g(u')}_{L^1_{\gr}(\Gw)}\\
&\leq c\left(\norm{ f-f'}_{L^1_{\gr}(\Gw)}+\norm {\BBP[h-h']}_{L^1_{\gr}(\Gw)}\right) \notag\\
&\leq c\left(\norm {f-f'}_{L^1_{\gr}(\Gw)}+\norm {h-h'}_{L^1(\bdw,\gw)}\right).\notag
\end{align}
Finally, for any nonnegative element $\eta\in X(\Gw)$, we have
\begin{equation}\label {NL5'}
-\int_{\Gw}{}\abs u\Gd\eta\,dx + \int_\Gw |g(u)|\eta\,dx \leq
-\int_{\Gw}{}\BBP[\abs {h}]\Gd\eta dx+\myint{\Gw}{}\eta f{\rm sgn}(u)\,dx,
\end{equation}
and
\begin{equation}\label {NL5+}
-\myint{\Gw}{}u_{+}\Gd\eta\,dx+ \int_\Gw g(u)_+\eta\,dx \leq -\myint{\Gw}{}\BBP[h_{+}]\Gd\eta dx+\myint{\Gw}{}\eta f{\rm sgn}_{+}(u)\,dx.
\end{equation}
\es
\Proof If $u,u'$ are two solutions as stated above then $v=u-u'$ satisfies
\begin{equation}\label{nln'}
\myint{\Gw}{}\left(-v\Gd\eta+F\eta\right)dx=-\myint{\Gw}{}
\BBP[h-h']\Gd hdx\forevery\eta\in X(\Gw)
\end{equation}
where $F=g(u)-g(u')-(f-f')\in L^1_\gr(\Gw)$. Applying \eqref{L5'} to
this equation and using the properties of $g$ described in
\rdef{class} we obtain \eqref{NL5}. Similarly we obtain \eqref{NL5'}
and \eqref{NL5+}, using \eqref{L5'} and \eqref{L5+}. These
inequalities imply uniqueness and monotone dependence on data.

In the case that $f$ and $h$ are bounded, existence is obtained by the standard variational
method. In general we approach $f$ in $L^1_\gr(\Gw)$ by functions in $C_c^\infty(\Gw)$ and $h$ in
$L^1(\bdw;\gw)$ by functions in $C(\bdw)$ and employ \eqref{NL5}.
\qed

\mysection{Measure data}
Denote by $\GTM_\gr(\Gw)$ the space of Radon measures $\nu$ in $\Gw$ \sth
$\gr|\nu|$ is a bounded measure.
\blemma{tr} Let $\Gw$ be a Lipschitz bounded domain.
Let $\nu\in \GTM_\gr(\Gw)$ and $u\in L^1_{loc}(\Gw)$ be a nonnegative solution of
$$-\Gd u=\nu\quad\text {in }\Gw.
$$
Then $u\in L^1_\gr(\Gw)$ and there exists a unique positive Radon measure $\gm$ on $\prt\Gw$
 such that
\begin{equation}\label{tr1}
u=\BBK[\gm]+\BBG[\nu].
\end{equation}
\es
\Proof Let $D$ be a smooth subdomain of $\Gw$ \sth $\bar D\sbs \Gw$.
Since $u\in W^{1,p}_{loc}(\Gw)$ for some $p>1$ it follows that $u$
possesses a trace, say $h_D$, in $W^{1-\rec{p},p}(\prt D)$. Put
$v:=u-\BBG^D[\nu]$. Then $-\Gd v=0$ in $D$ and $v\geq 0$ on $\prt D$
and therefore in $D$. If $\set{D_n}$ is an increasing \seq of such
domains, converging to $\Gw$, then $\BBG^{D_n}[\nu]\uparrow
\BBG^\Gw[\nu]$. Thus $v=u-\BBG^\Gw[\nu]$ is a non-negative harmonic
function in $\Gw$ and \consy possesses a boundary trace
$\mu\in\GTM(\bdw)$ \sth $v=\BBK[\mu]$. \qed
\blemma{M-lin} Let $\Gw$ be a Lipschitz bounded domain. If  $\nu\in\GTM_\gr(\Gw)$ and
$\mu\in  \GTM(\bdw)$, there exists a unique $u\in L^1_{\gr}(\Gw)$ satisfying
\begin{equation}\label {ML4}
\int_{\Gw}{}-u\Gd\eta\,dx=\int_\Gw \eta\,d\nu-\int_{\Gw}{}\BBK[\mu]\Gd\eta dx\forevery\eta \in X(\Gw).
\end{equation}
This is equivalent to
\begin{equation}\label {ML4eq}
u=\BBG[\nu]+\BBK[\mu].
\end{equation}
The following estimate holds
\begin{align}\label {ML5}
\norm u_{L^1_{\gr}(\Gw)}&\leq c\left(\norm {\nu}_{\GTM_{\gr}(\Gw)}+
\norm {\BBK[\mu]}_{L^1_{\gr}(\Gw)}\right)\\
&\leq c\left(\norm{\nu}_{\GTM_{\gr}(\Gw)}+
\norm {\mu}_{\GTM(\bdw)}\right).\notag
\end{align}
In addition, if $d\nu=fdx$ for some $f\in L^1_\gr(\Gw)$ then, for any nonnegative element
$\eta\in X(\Gw)$, we have
\begin{equation}\label {ML5'}
-\int_{\Gw}{}\abs u\Gd\eta\,dx  \leq
-\int_{\Gw}{}\BBK[\abs {\mu}]\Gd\eta dx+\myint{\Gw}{}\eta f{\rm sgn}(u)\,dx,
\end{equation}
and
\begin{equation}\label {ML5+}
-\myint{\Gw}{}u_{+}\Gd\eta\,dx\leq -\myint{\Gw}{}\BBK[\mu_{+}]\Gd\eta dx+\myint{\Gw}{}\eta f{\rm sgn}_{+}(u)\,dx.
\end{equation}
 \es
\Proof   We approximate $\mu$ by a \seq $\set{h_nP(x_0,\cdot)}$ and
$\nu$ by a \seq $\set{f_n}$ \sth
$$h_nP(x_0,\cdot)\in L^1(\bdw),\q h_nP(x_0,\cdot)\CH_{N-1}\to \mu \q\text{weakly in measure}$$
 and
$$f_n\in L^1_\gr(\Gw),\q f_n\to\nu \q\text{weakly relative to $C_\gr(\Gw)$,}$$
where $C_\gr$ denotes the space of  functions $\gz\in C(\Gw)$ \sth $\gr\gz\in L^\infty(\Gw)$.
Applying \rlemma{lin} to problem \eqref{L5} ($f,h$ replaced by $f_n,h_n$) and taking the limit
we obtain a solution
$u\in L^1_{\gr}(\Gw)$ of \eqref{ML4} satisfying \eqref{ML5}.

\rlemma{equivlin} implies that any solution $u$ of \eqref{ML4} satisfies \eqref{ML4eq}. Therefore the solution is unique and hence
\eqref{ML5} holds for all solutions.

Inequalities \eqref{ML5'} and \eqref{ML5+} are proved in the same way as the corresponding
 inequalities in \rlemma{lin}
\qed 
\bdef {solmeas}Let $\Gw$ be a bounded Lipschitz domain and let $g\in
\CG(\BBR)$. If $\gm\in\GTM(\prt\Gw)$, a function $u\in
L_{\gr}^1(\Gw)$ is a weak solution of
\begin{equation}\label{meas1}\left\{\BA {l}
-\Gd u+g(u)=0\quad\text {in }\Gw\\
\phantom{-\Gd u+g()} u=\gm\quad\text {in }\prt\Gw
\EA\right.\end{equation}
if $g(u)\in L^1_\gr(\Gw)$ and
\begin{equation}\label{meas2}
u+\BBG[g(u)]=\BBK[\gm]
\end{equation}
a.e. in $\Gw$. Equivalently
\begin{equation}\label{meas2bis}
\int_{\Gw}{}\left(-u\Gd\eta+g(u)\eta\right)dx=-\int_{\Gw}{}\left(\BBK[\gm]\Gd\eta\right)dx\forevery
\eta\in X(\Gw).
\end{equation}
The measure $\mu$ is called the {\em boundary trace} of $u$ on
$\bdw$.

Similarly a function $u\in L_{\gr}^1(\Gw)$ is a weak supersolution,
respectively subsolution, of \eqref{meas1} if $g(u)\in L^1_\gr(\Gw)$
and
\begin{equation}\label{meas-ssol}
  u+\BBG[g(u)]\geq \BBK[\gm] \txt{\rm respectively} u+\BBG[g(u)]\leq
  \BBK[\gm].
\end{equation}
This is equivalent to \eqref{meas2bis}, with $=$  replaced by $\geq$
or $\leq$, holding for every positive $\eta\in X(\Gw)$.
\es 
\Remark It follows from this definition and \rlemma{Sol} that,  if
$$\gm_{n}\rightharpoonup\gm  \txt{weakly in $\GTM(\prt\Gw)$,}
u_{n}\to u, \q  g(u_{n})\to g(u) \txt{in $L^1_{\gr}(\Gw)$,}$$
 and if
$$u_{n}=\BBK[\gm_{n}]-\BBG[g(u_{n})],$$
then $u=\BBK[\gm]-\BBG[g(u)]$.
\blemma{M-nln} Let $\Gw$ be a Lipschitz bounded domain and let $g\in \CG$. Suppose that
$\mu\in  \GTM(\bdw)$ and that
there exists a solution of problem \eqref{meas1}. Then the solution is unique.

If $\mu,\mu'$ are two measures in $\GTM(\bdw)$, for which problem
\eqref{meas1} possesses solutions $u,u'$ respectively, then the
following estimate holds:
\begin{align}\label {MNL5}
\norm{u- u'}_{L^1_{\gr}(\Gw)}+\norm{g(u)-g(u')}_{L^1_\gr(\Gw)}&\leq
\norm {\BBK[\mu-\mu']}_{L^1_{\gr}(\Gw)})\\
&\leq \norm {\mu-\mu'}_{\GTM(\bdw)}.\notag
\end{align}
If $\mu\leq \mu'$ then  $u\leq u'$.

In addition, for any nonnegative element
$\eta\in X(\Gw)$, we have
\begin{equation}\label {MNL5'}
-\int_{\Gw}(\abs u\Gd\eta-|g(u)|\eta)\,dx  \leq
-\int_{\Gw}{}\BBK[\abs {\mu}]\Gd\eta dx
\end{equation}
and
\begin{equation}\label {MNL5+}
-\int_{\Gw}(u_{+}\Gd\eta- g(u)_+\eta)\,dx \leq
-\int_{\Gw}{}\BBK[\mu_{+}]\Gd\eta dx.
\end{equation}
 \es
\Proof This follows from \rlemma{M-lin} in the same way that
\rlemma{Sol} follows from \rlemma{lin}. \qeda

 \bdef{ext-tr} Assume that $u\in W_{loc}^{1,p}(\Gw)$ for some $p>1$.
 We say that $u$ possesses a boundary
trace $\mu\in\GTM(\bdw)$ if, for every \Lip exhaustion $\set{\Gw_n}$
of $\Gw$,
\begin{equation}\label{tr-w-conv1}
   \lim_{n\to\infty}\int_{\bdw_n}Zu\,d\gw_n=\int_{\bdw}Z\,d\mu,
\end{equation}
 holds for every $Z\in C(\bar\Gw)$.

Similarly we say that $u$ possesses a trace $\mu$ on a relatively
open set $A\sbs \bdw$ if \eqref{tr-w-conv1} holds for every $Z\in
C(\bar\Gw)$ \sth $\supp Z\sbs \Gw\cup A$. \es

\Remark If $u\in W_{loc}^{1,p}(\Gw)$ for some $p>1$ then, by
Sobolev's trace theorem, for every relatively open $(N-1)$-
dimensional \Lip surface $\Gs$, $u$ possesses a trace in
$W^{1-\rec{p},p}(\Gs)$. In particular the trace is in $L^1(\Gs)$. In
fact there exists an element of the Lebesgue equivalence class of
$u$ \sth the trace on $\Gs$ is precisely the restriction of $u$ to
$\Gs$. When it is relevant, as in \eqref{tr-w-conv1}, we assume that
$u$ is represented by such an element.

If $u\in W^{1,p}(\Gw)$ then, by the same token, $u$ possesses a
trace in $W^{1-\rec{p},p}(\bdw)$. If $\set{\Gw_n}$ is a uniform \Lip
exhaustion and $h_n$ (resp. $h$) denotes the trace of $u$ on
$\bdw_n$ (resp. $\bdw$) then
$$\norm{h_n}\indx{W^{1-\rec{p},p}(\bdw_n)}\to
\norm{h}\indx{W^{1-\rec{p},p}(\bdw)}.$$

\nind{}
This follows from the continuity of the imbedding
$$W^{1,p}(\Gw)\hookrightarrow W^{1-\rec{p},p}(\bdw)$$
and the fact that $C^1(\bar \Gw)$ is dense in $W^{1,p}(\Gw)$.

Similarly, if $\set{\Gw_n}$ is a  \Lip exhaustion (not necessarily
uniform, but  satisfies \eqref{exhaustion}) then
$$\norm{h_n}\indx{L^1(\bdw_n)}\to
\norm{h}\indx{L^1(\bdw)}.$$

In particular, if $u\in W^{1,p}_0(\Gw)$ then its boundary trace is
zero, in the sense of the above definition.

\bprop{trace=limit} Let $u$ be a weak solution of \eqref{meas1}. If
$\{\Gw_n\}$ is a \Lip exhaustion of $\Gw$ then, for every $Z\in
C(\bar \Gw)$,
\begin{equation}\label{tr-w-convN}
   \lim_{n\to\infty}\int_{\bdw_n}Zu\,d\gw_n=\int_{\bdw}Z\,d\mu,
\end{equation}
where $\gw_n$ is the harmonic measure of $\Gw_n$ \(relative to a
point $x_0\in\Gw_1$\). \es

\Proof If $v:=\BBG[ g\circ u]$ then $v\in L^1_\gr(\Gw)$ and $u+v$ is
a harmonic function.  By \eqref{meas2}, $u+v=\BBK^\Gw[\mu]$.
Therefore, by \rlemma{tr-w-conv},
\begin{equation}\label{tr-w-conv1'}
  \lim_{n\to\infty}\int_{\bdw_n}Z(u+v)\,d\gw_n=\int_{\bdw}Z\,d\mu
\end{equation}
for every $Z\in C(\bar \Gw)$. As $v\in W^{1,p}_0(\Gw)$ for some
$p>1$ its boundary trace is zero. Therefore \eqref{tr-w-conv1'}
implies \eqref{tr-w-convN}.
 \qeda 

\bdef{admissible} A measure $\mu\in \GTM(\bdw)$ is called
$g$-admissible if $g(\BBK[|\mu|])\in L^1_\gr(\Gw)$. \es
\bth{admissible} If $\mu$ is $g$-admissible then problem
\eqref{meas1} possesses a unique solution. \es

\Proof  First assume that $\mu>0$. Under the admissibility
assumption, $U=\BBK[\mu]$ is a supersolution of \eqref{meas1}. Let
$\set{D_n}$ be an increasing \seq of smooth domains \sth $\bar
D_n\sbs D_{n+1}\sbs \Gw$ and $D_n\uparrow \Gw$. Let $u_n$ be the
solution of problem \eqref{meas1} in $D_n$ with boundary data
$h_n=U\big|_{\prt D_n}$. Then $\set{u_n}$ decreases and the limit
$u=\lim u_n$ satisfies \eqref{meas1}.

In the general case we define $\bar U=\BBK[|\mu|]$ and $U$, $u_n$ as
before. By assumption $g(\bar U)\in L^1_\gr(\Gw)$ and $\bar U$
dominates $|u_n|$ for all $n$. Let $\eta$ be a non-negative function
in $X(\Gw)$ and let $\gz_n$ be the solution of the problem
$$\Gd \gz=(\Gd \eta)\gr_n/\gr \txt{in $D_n$,} \gz=0 \txt{on $\prt D_n$.}$$
Then $\gz_n\in X(D_n)$ and, since $\gr_n\to\gr$,
$$(\Gd\gz_n)\to(\Gd\eta),\q \gz_n\to\eta.$$
In addition, $(\Gd\gz_n)/\gr_n=(\Gd\eta)/\gr$ is bounded and, by \eqref{ugr-bound},
 the sequence
$\{\gz_n/\gr_n\}$ is uniformly bounded.

The solutions $u_n$ satisfy,
\begin{equation}\label{meas-n}
\int_{D_n}{}\left(-u_n\Gd\gz_n+g(u_n)\gz_n\right)dx=-\int_{D_n}\BBP^{D_n}[h_n]\Gd\gz_ndx.
\end{equation}
 The \seq
$\set{u_k:k>n}$ is bounded in $W^{1,p}(D_n)$ for every $n$. \Consy
there exists a \sseq (still denoted by $\set{u_n}$) which converges
pointwise a.e. in $\Gw$. We denote its limit by $u$. Since
$\set{u_n}$ is dominated by $\bar U$ it follows that
$$\lim_{n\to\infty}\int_{D_n}\left(-u_n\Gd\gz_n+g(u_n)\gz_n\right)dx=
\int_{\Gw}\left(-u\Gd\eta+g(u)\eta\right)dx.$$
Furthermore,
$$\int_{D_n}\BBP^{D_n}[h_n]\Gd\gz_ndx=\int_{D_n} U\Gd\eta(\gr_n/\gr)\,dx\to \int_\Gw U\Gd\eta dx=
\int_\Gw \BBK[\mu]\Gd\eta\,dx.$$
Thus $u$ is the solution of \eqref{meas1}.
\qed

\Remark If we do not assume that $g(0)=0$ the admissibility
condition becomes,
\begin{equation}\label{admiss}
g(\BBK[\gm_{+}]+\gr (g(0))_{+})\in L^1_\gr(\Gw)\quad\text {and }\;
g(-\BBK[\gm_{-}]-\gr (g(0))_{-})\in L^1_\gr(\Gw).
 \end{equation}
\mysection{The boundary trace of positive solutions} \hskip 5mm As
before we assume that $\Gw$ is a bounded Lipschitz domain and $g\in
\CG$. We denote by $\gr$ the first eigenfunction of $-\Gd$ in $\Gw$
normalized by $\gr(x_0)=1$ at some (fixed) point $x_0\in \Gw$.

A function
$u\in L^1_{loc}(\Gw)$ is a solution of the equation
\begin{equation}\label{eq1}
-\Gd u+g(u)=0\quad \text{in }\Gw,
\end {equation}
if $g\circ u\in L^1_{loc}(\Gw)$ and $u$ satisfies the equation in
the distribution sense.

A function $u\in L^1_{loc}(\Gw)$ is a supersolution (resp.
subsolution) of the equation \eqref{eq1} if $g\circ u\in
L^1_{loc}(\Gw)$ and
$$-\Gd u+g\circ u \geq 0 \q\text{(resp. $\leq 0$)}$$
in the distribution sense.

 \bprop{reg-sol} Let $u$ be a positive
solution of \eqref{eq1}. If  $g\circ u\in L^1_\gr(\Gw)$ then $u\in
L^1_\gr(\Gw)$ and it possesses a boundary trace $\mu\in\GTM(\bdw)$,
i.e., $u$ is the solution of the \bvp \eqref{meas1} with this
measure $\mu$.

 \es

\Proof If $v:=\BBG[ g\circ u]$ then $v\in L^1_\gr(\Gw)$ and $u+v$ is
a positive harmonic function. Hence $u+v\in L^1_\gr(\Gw)$ and there
exists a non-negative measure $\mu\in\GTM(\bdw)$ such that
$u+v=\BBK[\mu]$. In view of \eqref{meas2}, this implies our
assertion.
 \qeda 

 \blemma{g-basic1} If $u$ is a non-negative solution of
\eqref{eq1} then $u\in C^1(\Gw)$.

Let $\set{u_n}$ be a \seq of non-negative solutions of \eqref{eq1}
which is uniformly
bounded in every compact subset of $\Gw$.
 Then there exists a \sseq
$\set{u_{n_j}}$ which converges in $C^1(\bar\Gw')$ for every $\Gw'\Subset\Gw$
to a solution $u$ of \eqref{eq1}. \es

\noindent\Proof  Since $g\circ u\in L^1\loc(\Gw)$ it follows that
$u\in W^{1,p}\loc (\Gw)$ for some $p\in [1,N/(N-1))$. Let $\Gw'$ be
a smooth domain \sth $\Gw'\Subset\Gw$. By the trace imbedding
theorem, $u$ possesses a trace  $h\in L^1(\bdw')$. If $U$ is the
harmonic function in $\Gw'$ with boundary trace $h$ then $u<U$. Thus
 $u$ (and hence $g\circ u$) is bounded in every
compact subset of $\Gw$. By elliptic p.d.e. estimates, $u\in
C^1(\Gw)$.

The second assertion of the lemma follows from the first by a
standard  argument. \qeda

\bth{ssol-1} \(i\) Let $u$ be a non-negative supersolution \(resp.
subsolution\) of \eqref{eq1}. Then $u\in W^{1,p}\loc (\Gw)$ for some
$p\in [1,N/(N-1))$. In particular, if $\Gw'$ is a $C^1$ domain \sth
$\Gw'\Subset\Gw$ then $u$ possesses a trace  $h\in L^1(\bdw')$.

\medskip
\nind\(ii\) If $u$ is a positive supersolution, there exists a
non-negative solution $\unl u\leq u$ which is the largest among all
solutions dominated by $u$.

 If $u$ is a  positive subsolution and $u$ is dominated by a
solution $w$ of \eqref{eq1} then there exists a minimal solution
$\bar u$ \sth $u\leq \bar u$. In particular, if $g\in \CG$ satisfies
the Keller-Osserman condition then such a solution exists.


\medskip
\nind\(iii\) Under the assumptions of \(ii\),  if $g\circ \unl u\in
L^1_\gr(\Gw)$ \(resp. $g\circ \bar u\in L^1_\gr(\Gw)$\) then
 the boundary trace of $\unl u$ \(resp. $\bar
u$\)  is also the boundary trace of $u$ in the sense of
\rdef{ext-tr}.

 \es

\nind\Proof First consider the case of a supersolution. Since $-\Gd
u+g(u)\geq 0$ there exists a positive Radon measure $\tau$ in $\Gw$
\sth
$$-\Gd u+g(u)=\tau\txt{in $\Gw$.}$$
Therefore $u\in W^{1,p}\loc (\Gw)$ and \consy $u$ possesses an $L^1$
trace on $\bdw'$ for every $\Gw'$ as above.

Next, let $\set{\Gw_n}$ be a $C^1$ exhaustion of $\Gw$ which is also
uniformly \Lip. Let $v_n$ be the solution of the \bvp
\begin{equation}\label{vn=u}
  -\Gd v+g(v)=0\txt{in $\Gw_n$,}v=u\txt{on $\bdw_n$.}
\end{equation}
Since $u$ possesses a trace in $L^1(\bdw_n)$ this \bvp possesses a
(unique) solution. By the comparison  principle $0\leq v_n\leq u$ in
$\Gw_n$. Therefore the \seq $\set{v_n}$ decreases and \consy it
converges to a solution $\unl u$ of \eqref{eq1}. Evidently this is
the largest solution dominated by $u$.

Now suppose  that $g\circ \unl u\in L^1_\gr(\Gw)$ (but not
necessarily $g\circ u\in L^1_\gr(\Gw)$). By \rprop{reg-sol}, $\unl
u\in L^1_\gr(\Gw)$ and  $\unl u$ possesses a boundary trace $\mu$.
By the definition of $v_n$,
$$\BAL \int_{\bdw_n}ud\gw_n=\int_{\bdw_n}
P^{\Gw_n}(x_0,y)u(y)dS &=v_n(x_0)+\int_{\Gw_n}
G^{\Gw_n}(x,x_0)g(v_n(x))dx\\
&\to \unl u(x_0)+\int_\Gw G^{\Gw}(x,x_0)g(\unl u(x))dx. \EAL $$
Hence, taking a \sseq if necessary, we may assume that
$$u\chi\indx{\bdw_n}\gw_n\rightharpoonup \mu'$$
where $\mu'$ is a measure on $\bdw$ \sth
$$\mu'(\bdw)= \unl u(x_0)+\int_\Gw G^{\Gw}(x,x_0)g(\unl u(x))dx.$$
 On the other hand, as $\mu$ is the
boundary trace of $\unl u$,
$$\unl u(x_0)+\int_\Gw G^{\Gw}(x,x_0)g(\unl u(x))dx=\mu(\bdw).$$
Thus $\mu(\bdw)=\mu'(\bdw)$. However, as $\unl u\leq u$, we have
$\mu\leq \mu'$. This implies that $\mu=\mu'$.

 Next we treat the case of a subsolution. The proof of (i) is the same
 as before. We turn to (ii). In the present case,
 the corresponding \seq
$\set{v_n}$ is increasing and, in general, may not converge. But, as
we assume that $u$ is dominated by a solution $w$, the \seq
converges to a solution $\bar u$ which is clearly the smallest
solution above $u$. In particular, if $g$ satisfies the
Keller-Osserman condition then $\set{v_n}$ is uniformly bounded in
every compact subset of $\Gw$ and \consy converges to a solution.

The proof of (iii) for subsolutions is again the same as in the
 case of supersolutions. \qed

\bcor{ssol-2} {\rm I.}  Let $u$ be a non-negative supersolution of
\eqref{eq1}. Let $A$ be a relatively open subset of $\bdw$. Suppose
that, for every \Lip domain $\Gw'$ \sth
\begin{equation}\label{Gw'}
   \Gw'\sbs \Gw, \q \bdw'\cap\bdw\sbs A,
\end{equation}
we have
\begin{equation}\label{g(u)}
  g\circ u\in L^1_\gr(\Gw').  
\end{equation}
 Then both $u$ and
$\unl u$  possess    traces on $A$ and the two traces are equal.

\medskip
\nind{\rm II.}  Let $u$ be a non-negative subsolution of
\eqref{eq1}. Let $A$ be a relatively open subset of $\bdw$. Suppose
that  for every \Lip domain $\Gw'$ satisfying \eqref{Gw'} we have
\begin{equation}\label{g(bar u)}
  g\circ \bar u\in L^1_\gr(\Gw').
  \end{equation}
  Then both $u$ and
$\bar u$  posses    traces on $A$ and the two traces are equal. \es

\Proof  Let $u$ be a supersolution and let $\Gw'$ be a domain as
above. Denote by $\gr'$ the first eigenfunction of $-\Gd$ in $\Gw'$
normalized by $\gr'(x_0)=1$ for some $x_0\in \Gw'$. Since $\gr'\leq
c\gr$, \eqref{Gw'} implies that $g\circ  u\in L^1_{\gr'}(\Gw')$. Let
$\unl u'$ denote the largest solution of \eqref{eq1} in $\Gw'$
dominated by $u$.  Then $g\circ  \unl u'\in L^1_{\gr'}(\Gw')$ and,
by \rth{ssol-1}, $\unl u'\in L^1_\gr(\Gw')$ and $\unl u'$ has a
trace $\nu'$ on $\bdw'$ which is also the boundary trace of $u$ on
$\bdw$.

Let $\set{\Gw_n}$ be an increasing uniformly \Lip \seq of domains
\sth $\bdw_n\cap\Gw$ is a $C^1$ surface, $D_n:=\Gw\sms \Gw_n$ is
\Lip and
$$F_n:=\bdw_n\sms \Gw\sbs F_{n+1}^0\sbs A,\q \cup\Gw_n=\Gw,\q \cup F^0_n=A,$$
where $F^0_n$ is the relative interior of $F_n$. Denote by $\unl
u_n$ the largest solution dominated by $u$ in $\Gw_n$ and observe
that $\set{\unl u_n}$ is decreasing and converges to a solution.
Obviously this is the largest solution dominated by $u$, namely,
$\unl u$.

Let $\tau_n$ be the trace of $\unl u_n$ on $\bdw_n$. Put
$\nu_n=\tau_n\chi\indx{F_n}$. Recall that $\tau_n$ is also the trace
of $u$ so that
$$\nu'_n=\tau_n-\nu_n=u\chi\indx{\bdw_n\sms F_n}dS.$$



\nind {\em Assertion A.}\hskip 2mm {\em  There exists a Radon
measure $\nu$ on $A$ \sth $\nu_n\rightharpoonup \nu$
 and $\nu$ is the trace of $u$, as well as of $\unl u$, on $A$.}

\medskip
Let $E$ be a compact subset of $A$ and denote,
$$n(E):=\inf\set{m\in \BBN: \, E\sbs F_m^0}.$$
 In view of the fact that, for $n\geq n(E)$,
$\nu_n$ is the trace of $u$, relative to $\Gw_n$, on a set
$F_{n(E)}^0$ in which $E$ is strongly contained and the fact that
$\set{\Gw_n}$ is \Lip, \rlemma{D'sbsD} implies that the set
$\set{\nu_n(E):n\geq n(E)}$ is bounded. 
By taking a \seq if necessary we may assume that
$$\nu_n\lfloor\indx{E}\rightharpoonup \nu_E.$$
Applying this procedure to  $E=F_m$ for each $m\in \BBN$ and then
using the  diagonalization method  we obtain a \sseq, again denoted
by $\set{\nu_n}$, \sth
$$\nu_n\rightharpoonup \nu$$
where $\nu$ is a Radon measure on $A$ (not necessarily bounded).

 Next we
wish to show that $\nu$ is the trace of $u$ on $A$ relative to
$\Gw$. To this purpose
 we construct a $C^1$ exhaustion  of $\Gw$, say $\set{D_n}$, \sth $D_n\Subset \Gw_n$ and
 $\prt D_n=\Gg_n\cup \Gg'_n$ where
 $$\BAL
 \Gg'_n&=\prt \Gw_n\cap \set{y\in \Gw:\dist(y,F_{n})\geq \ge_n}\\
 \Gg_n&\sbs \set{y\in \Gw_n:\dist(y,F_{n})< \ge_n},
 \EAL$$
 where $0<\ge_n<\rec{2}\dist(F_{n},\bdw\sms A)$ is chosen so that
 $$\BBH\indx{N-1}\chi\indx{\Gg_n}\rightharpoonup\BBH\indx{N-1}\chi\indx{A}
\txt{and} u\chi\indx{\Gg_n}d\gw^n\rightharpoonup\nu.$$ Here $d\gw^n$
is the harmonic measure in $D_n$.
This is possible because, if $\Gg_n$ is sufficiently close to
$\bdw_n$, then
$$u\chi\indx{\Gg_n}d\gw^n-\nu_n\chi\indx{F_n}\rightharpoonup 0.$$
(As usual in this paper, $\nu_n\chi\indx{F_n}$ denotes the Borel
measure in $\BBR^N$ that is equal to $\nu_n$ on $F_n$ and zero
elsewhere.) This implies that $\nu$
 is the trace of $u$ on $A$.

Since $\nu_n$ is also the trace of $\unl u_n$ on $F_n$ it follows
that, if $\Gg_n$ is sufficiently close to $\bdw_n$,
$$\unl u_n\chi\indx{\Gg_n}d\gw^n-\nu_n\chi\indx{F_n}\rightharpoonup 0.$$
 As $\unl
 u_n\downarrow \unl u$ we deduce that $\nu$ is also the trace of $\unl u$ on
 $A$.

  If $u$ is a subsolution the argument is essentially the same. Let $\bar u_n$ be
  the smallest solution that dominates $u$ in $\Gw_n$.
  Then the \seq $\set{\bar u_n}$ is increasing, but it is
  dominated by a solution $w$. Therefore it converges to a solution
  and this is the smallest solution dominating $u$, namely,
  $\bar u$. By \rth{ssol-1}, $\unl u_n$ and $u\lfloor\index{\Gw_n}$
  possess the same trace on $\bdw_n$.  Let $\tau_n$ be the trace of $\unl u_n$ on
  $\bdw_n$ and put $\nu_n=\tau_n\chi\indx{F_n}$. The rest of the
  proof is as before.
  \qed

 \bdef{trace} Let $u$ be a positive
supersolution, respectively subsolution, of \eqref{eq1}. A point
$y\in \bdw$ is a {\em regular boundary point} relative to $u$ if
there exists an open \ngh $D$ of $y$ \sth $g\circ u\in
L^1_\gr(\Gw\cap D)$. If no such \ngh exists we say that $y$ is a
{\em singular boundary point} relative to $u$.

The set of regular boundary points of $u$ is denoted by $\CR(u)$;
its complement on the boundary is denoted by $\CS(u)$. Evidently
$\CR(u)$ is relatively open.
 \es

\bth{reg-trace} Let $u$ be a positive solution of \eqref{eq1} in
$\Gw$. Then $u$ possesses a trace  on $\CR(u)$, given by a Radon measure $\nu$.

Furthermore, for every compact set $F\sbs \CR(u)$, 
\begin{equation}\label{meas2b}
\int_{\Gw}\left(-u\Gd\eta+g(u)\eta\right)dx=-\int_{\Gw}\left(\BBK[\gn\chi\indx{F}]\Gd\eta\right)dx
\end{equation}
for every $\eta\in X(\Gw)$ \sth $\supp \eta\cap\bdw\sbs F$. \es
\Proof The first assertion is an immediate consequence of
\rcor{ssol-2}.

We turn to the proof of the second assertion. Let $F$ be a compact
subset of $\CR(u)$ and let $\eta\in X(\Gw)$ be a function \sth the
following conditions hold for some open set $E_\eta$:
$$\supp\eta\sbs
\bar\Gw\cap E_\eta,\q F\sbs E_\eta\cap\bdw,\q \bar
E_\eta\cap\CS(u)=\ems,\q x_0\in D_\eta:=\Gw\cap E_\eta.$$

\noindent By \rdef{trace}, if $D$ is a subdomain of $\Gw$ \sth $\bar
D\cap \CS(u)=\ems$ then $g\circ u\in L^1_\gr(D)$, where $\gr$ is the
first normalized eigenfunction of $\Gw$. Let $E$ be a $C^2$ domain
\sth
$$\bar E_\eta\sbs E, \q \BBH_{N-1}(\bdw\cap\prt E)=0,\q \bar
E\cap\CS(u)=\ems.$$

\noindent Put $D:=E\cap \Gw$ and note that $g\circ u\in L^1_\gr(D)$.


 If $\gf$ denotes the first normalized eigenfunction in $D$
then $\gf\leq c\gr$ for some positive constant $c$. Therefore the
fact that $g\circ u\in L^1_\gr(D)$ implies that $g\circ u\in
L^1_\gf(D)$ and the properties of $\eta$ imply that  $\eta\in X(D)$.
Hence $u$ possesses a boundary trace $\tau^D$ on $\prt D$ and
\begin{equation}\label{tauD}
  \int_{D}\left(-u\Gd\eta+g(u)\eta\right)dx=-\int_{D}\left(\BBK^D[\tau^D]\Gd\eta\right)dx.
\end{equation}

Let $\Gg= \bar E\cap\bdw$ and $\Gg'=\prt D\sms \Gg$; note that
$\Gg\cap\CS(u)=\ems$ and $\eta$ vanishes in a \ngh of $\prt E\cap
\bar \Gw$. Put $\tau_\Gg^D=\tau^D\chi\indx{\Gg}$ and
$\tau^D_{\Gg'}=\tau^D-\tau^D_\Gg$. Then $d\tau^D_{\Gg'}=udS$ on
$\Gamma'$ and, as $u\in C(\bar D\sms \Gg)$,
$$\BBK^D[\tau^D_{\Gg'}]\in C(\bar D\sms \Gamma).$$
Furthermore $\eta$ vanishes in a \ngh of $\Gamma'$ and \consy
$$\BAL\int_{D}\left(\BBK^D[\tau^D_{\Gg'}]\Gd\eta\right)dx&= \int_D\left(\int_{\prt
D\sms \Gg} P^D(x,y)u(y)dS_y\right)\Gd\eta(x)dx\\
&=\int_{\prt D\sms
\Gg}\left(\int_DP^D(x,y)\Gd\eta(x)dx\right)u(y)dS_y=0.\EAL$$
Thus
\begin{equation}\label{tauGD}
\int_{\Gw}\left(-u\Gd\eta+g(u)\eta\right)dx=-\int_{\Gw}\BBK^D[\tau_\Gg^D]\Gd\eta\,dx.
\end{equation}
(Changing the domain of integration from $D$ to $\Gw$ makes no
difference since $\eta$ vanishes in $\Gw\sms D$.)

Now, $\tau^D_\Gg$ is the trace of $u$ on $\Gg$ relative to $D$ while
$\nu\chi\indx{\Gg}$ is the trace of $u$ on $\Gg$ relative to $\Gw$.
Since $D\sbs \Gw$ it follows that
\begin{equation}\label{tdgg}
\tau^D_\Gg\leq \nu\chi\indx{\Gg}.
\end{equation}

Let $\set{E^j}$ be an increasing sequence of $C^2$ domains \sth each
domain possesses the same properties as $E$ and,
\begin{equation}\label{D^j}
\bar E^j\cap\bdw=\bar E\cap\bdw=\Gamma, \txt{and}
D^j:=E^j\cap\Gw\uparrow \Gw.
\end{equation}
For each $j\in \BBN$ and $y\in \Gamma$, the function
$K^{D^j}(\cdot,y)$ is harmonic in $D^j$, vanishes on $\prt
D^j\sms\{y\}$ and $K^{D^j}(x_0,y)=1$. Furthermore the \seq
$\set{K^{D^j}(\cdot,y)}$ is non-decreasing. Therefore it converges
 uniformly in compact subsets of $(\Gw\cup\Gamma)\sms \{y\}$.
The limit is the corresponding kernel function in $\Gw$, namely
$K^\Gw(\cdot,y)$. (Recall that the kernel function is unique.)

In view of \eqref{tdgg}, the sequence $\set{\tau_\Gg^{D^j}}$ is
bounded. Therefore there exists a \sseq, which we still denote by
$\set{\tau_\Gg^{D^j}}$, \sth
$$\tau_\Gg^{D^j}\rightharpoonup\tau_\Gg$$
 weakly relative to $C(\Gamma)$. Combining these facts we obtain,
 $$\BBK^{D_j}[\tau_\Gg^{D^j}]\to \BBK^\Gw[\tau_\Gg].$$
Hence, by \eqref{tauD},
\begin{equation}\label{tauGg}
  \int_{\Gw}\left(-u\Gd\eta+g(u)\eta\right)dx=-\int_{\Gw}\left(\BBK^\Gw[\tau_\Gg]\Gd\eta\right)dx.
\end{equation}
Finally, as $\tau^{D_j}_\Gg$ is the trace of $u$ on $\Gg$ relative
to $D_j$ then, in view of \eqref{D^j},  the limit $\tau_\Gg$ is the
trace of $u$ on $\Gg$ relative to $\Gw$, i.e.,
 $$\tau_\Gg=\nu\chi\indx{\Gg}.$$
This relation and \eqref{tauGg} imply \eqref{meas2b}. \qed

\bth{ssol-3} \hskip 2mm{\rm I.}  Let $u$ be a positive supersolution
of \eqref{eq1} in $\Gw$ and let $\unl u$ be the largest solution
dominated by $u$. Then,
\begin{equation}\label{su=su*}
   \CS(u)=\CS(\unl u),\q \CR(u)=\CR(\unl u).
\end{equation}
Both $u$ and $\unl u$ possess a trace on $\CR(u)$ and the two traces
are equal.

\nind{\rm II.}\hskip 2mm Let $u$ be a positive subsolution of
\eqref{eq1} in $\Gw$ and let $\bar u$ be the smallest solution which
dominates $u$. If $u$ is dominated by a solution $w$ of \eqref{eq1}
then both $u$ and $\bar u$ possess a trace on $\CR(w)$ \(which is
contained in $\CR(u)$\) and the two traces are equal on this set.

In particular, if $\CR(w)=\CR(u)$ then \eqref{su=su*}, with $\unl u$
replaced by $\bar u$, holds and both $u$ and $\bar u$ possess a
trace on $\CR(u)$, the two traces being equal.

 \nind{\rm III.}\hskip 2mm Let $\nu$ denote the trace of
$u$ on $\CR(u)$. Then, for every compact set $F\sbs \CR(u)$,
\begin{equation}\label{s-meas2}
\int_{\Gw}\left(-u\Gd\eta+g(u)\eta\right)dx
\begin{cases}\geq-\int_{\Gw}\left(\BBK[\gn\chi\indx{F}]\Gd\eta\right)dx,&\text{$u$
supersolution,}\\
\leq-\int_{\Gw}\left(\BBK[\gn\chi\indx{F}]\Gd\eta\right)dx,&\text{$u$
subsolution}\end{cases}
\end{equation}
for every $\eta\in X(\Gw)$, $\eta\geq 0$, \sth $\supp \eta\cap\bdw\sbs F$.
\es

\nind \Proof Part  I. is a consequence of \rcor{ssol-2} I.

The first assertion in II. follows from \rcor{ssol-2} II. with
$A=\CR(w)$. The second assertion in II. is an immediate consequence
of the first.

By \rth{reg-trace}, $\unl u$ (resp. $\bar u$) satisfy
\eqref{meas2b}, where $\nu$ is the trace of $\unl u$ (resp. $\bar
u$) on $\CR(u)$. Since $\nu$ is also the trace of $u$ on $\CR(u)$ we
obtain statement III.
 \qed

 \bth{sing-trace} Assume that $g\in \CG$ satisfies the
Keller-Osserman condition.

\nind\(i\) Let $u$ be a positive solution of \eqref{eq1} and let
$\set{\Gw_n}$ be a \Lip exhaustion of $\Gw$.  If $y\in \CS(u)$ then,
for every nonnegative $Z\in C(\bar\Gw)$ \sth $Z(y)\neq0$
\begin{equation}\label{dynamic-S(u)}
  \lim \int_{\bdw_n}Zud\gw_n=\infty.
\end{equation}

\nind\(ii\) Let $u$ be a positive supersolution of \eqref{eq1} and
let $\set{\Gw_n}$ be a $C^1$ exhaustion of $\Gw$. If $y\in \CS(u)$
then \eqref{dynamic-S(u)} holds for every nonnegative $Z\in C(\bar\Gw)$ \sth
$Z(y)\neq0$. \es

The proof of satement (i) is essentially the same as for the
corresponding result in smooth domains \cite[Lemma 2.8]{MV4} and
therefore will be omitted. In fact the assumption that $g$ satisfies
the Keller-Osserman condition implies that the set of conditions II
in \cite[Lemma 2.8]{MV4} is satisfied. Here too, the Keller-Osserman
condition can be replaced by the weaker set of conditions II in the
same way as in \cite{MV4}.

Part (ii) is a consequence of \rth{ssol-3} and statement (i). \qed

\bdef{g-trace} Let $g\in \CG$. Let $u$ be a positive solution of
\eqref{eq1} with regular boundary set $\CR(u)$ and singular boundary
set $\CS(u)$. The Radon measure $\nu$ in $\CR(u)$ associated with
$u$  as in \rth{reg-trace} is called the {\em regular part of the
trace of $u$}. The couple $(\nu,\CS(u))$ is called the {\em boundary
trace} of $u$ on $\bdw$. This trace is also represented by the
(possibly unbounded) Borel measure $\bar\nu$ given by
\begin{equation}\label{barnu}
 \bar\nu(E)=\begin{cases} \nu(E), &\text{if $E\sbs
\CR(u)$}\\ \infty,&\text{otherwise.}\end{cases}
\end{equation}

The boundary trace of $u$ in the sense of this definition will be
denoted by $\tr{\bdw}u$.

Let
\begin{equation}\label{Vnu}
V_\nu:=\sup\set{u\indx{\nu\chi_F}: \;F\sbs \CR(u),\;F\text{
compact}}
\end{equation}
where $u\indx{\nu\chi_F}$ denotes the solution of \eqref{meas1} with
$\mu=\nu\chi_F$. Then $V_\nu$ is called the {\em semi-regular
component\/} of $u$. \es

\nind{\em Remark.} Let $\tau$ be a Radon measure on a relatively
open set  $A\sbs\bdw$. Suppose that for every compact set $F\sbs A$,
$u\indx{\tau\chi_F}$ is defined. If $V_\tau$ is defined as above, it
need not be a solution of \eqref{eq1} or even be finite.  However,
if $g$ satisfies the Keller--Osserman condition or if $u\indx{\tau\chi_F}$ is
dominated by a solution $w$, independent of $F$,  then $V_\tau$ is a solution.

\medskip
\bdef{removable} A compact set $F\sbs \bdw$ is removable relative to
\eqref{eq1} if the only non-negative solution $u\in C(\bar\Gw\sms
F)$ which vanishes on $\bar\Gw\sms F$ is the trivial solution $u=0$.
\es

\nind {\em Remark.} In the case of power nonlinearities in smooth
domains there exists a complete characterization of removable sets
(see \cite{MV3} and the references therein). In a later section we
shall derive such a characterization for a family of \Lip domains.

\medskip

\blemma{maxsol} Let $g\in\CG$ and assume that $g$ satisfies the
Keller-Osserman condition. Let $F\sbs\bdw$ be a compact set and
denote by $\CU_F$ the class of solutions $u$ of \eqref{eq1} which
satisfy the condition,
\begin{equation}\label{Fc-vanish}
   u\in C(\bar\Gw\sms F),\q u=0\txt{on $\bdw\sms F$.}
\end{equation}
Then there exists a function $U_F\in \CU_F$ \sth
$$u\leq U_F \forevery u\in \CU_F.$$

\nind Furthermore, $\CS(U_F)=:F'\sbs F$; $F'$ need not be equal to
$F$. \es

The proof is standard and will be omitted.

\bdef{maxsol} $U_F$ is called the {\em maximal solution} associated
with $F$. The set $F'=\CS(U_F)$ is  called the $g$-kernel of $F$ and
denoted by $k_g(F)$. \es

\nind{\em Note.} The situation $\CS(U_F)\subsetneq F$ occurs if and
only if there exists a closed set $F'\sbs F$ \sth $F\sms F'$ is a
non-empty removable set. In this case $U_F=U_{F'}$.

\blemma{UF-basics} Let $F_1, F_2$ be two compact subsets of $\bdw$.
Then,

\begin{equation}\label{Monotone1}
F_1\sbs F_2\Lra  U_{F_1}\leq U_{F_2}
\end{equation}
and
\begin{equation}\label{UF-ineq1}
   U_{F_1\cup F_2}\leq U_{F_1}+U_{F_2}.
\end{equation}

If $F$ is a compact subset of $\bdw$ and $\set{N_k}$ is a decreasing \seq of
relatively open neighborhoods of $F$ \sth $\bar N_{k+1}\sbs N_k$ and $\cap N_k=F$ then
\begin{equation}\label{UF-conv1}
  U_{\bar N_k}\to U_F
\end{equation}
uniformly in compact subsets of $\Gw$.
\es

\Proof The first statement is an immediate \cons of the definition of maximal solution.

Next we verify \eqref{UF-conv1}. By \eqref{Monotone1} the \seq $\set{U_{\bar N_k}}$ decreases
 and therefore it converges to a solution $U$. Clearly $U$ has trace zero outside $F$ so that $U\leq U_F$
 On the other hand, for every $k$, $U_{\bar N_k}\geq U_F$. Hence $U=U_F$

We turn to the verification of \eqref{UF-ineq1}. Let $u$ be a positive solution of \eqref{eq-q} which vanishes
on $\bdw\sms(F_1\cup F_2)$. We shall show that there exists solutions $u_1,u_2$ of \eqref{eq-q} \sth
\begin{equation}\label{u-ineq2}
  u_i=0 \txt{on $\bdw\sms F_i$,}  u\leq u_1+u_2.
\end{equation}
 First we prove this statement in the case where $F_1\cap F_2=\ems$. Let $E_1,E_2$ be $C^1$ domains \sth $\bar E_1\cap \bar E_2=\ems$ and $F_i\sbs E_i\cap \bdw$, (i=1,2).
Let $\set{\Gw_n}$ be a \Lip exhaustion of $\Gw$ and put $A_{n,i}=\bdw_n\cap E_i$, (i=1,2). Let $v_{n,i}$ be the solution of \eqref{eq-q} in $\Gw_n$  with boundary data $u\chi\indx{A_{n,i}}$ and $v_n$ be the solution in $\Gw_n$ with boundary data $u(1-\chi\indx{A_{n,1}\cup A_{n,2}})$. Then
$$u\leq v_n+v_{n,1}+v_{n,2}.$$
By taking a \sseq if necessary we may assume that the sequences $\set{v_n}$, $\set{v_{n,1}}$, $\set{v_{n,2}}$ converge.
Then $\lim v_{n,i}=U_i$ where $U_i$ vanishes on $\bdw\sms E_i$, (i=1,2). In addition, as the trace of $u$ on $\bdw\sms  (F_1\cup F_2)$ is zero, we have $\lim v_n=0$. Thus
$$u\leq U_1+U_2.$$
Now take   decreasing \seqs of $C^1$ domains $\set{E_{k,1}}$, $\set{E_{k,2}}$ \sth
$$\bar E_{k,1}\cap \bar E_{k,2}=\ems,  \q F_i\sbs E_{k,i}\cap\bdw, \q\bar E_{k,i}\cap\bdw\downarrow F_i\q i=1,2.$$
Construct $U_{k,i}$ corresponding to $E_{k,i}$ in the same way that $U_i$ corresponds to $E_i$. Then,
$$u\leq U_{k,1}+U_{k,2}$$
and, by \eqref{UF-conv1}, taking a \sseq if necessary,
$$u_i:=\lim_{k\tin}U_{k,i}=0 \txt{on $\bdw\sms F_i$},\q i=1,2. $$
This proves \eqref{u-ineq2} in the case where $F_1,F_2$ are disjoint.

In the general case, let $\set{N_j}$ be a decreasing \seq of relatively open neighborhoods of $F_1\cap F_2$ \sth

$$\bar N_{j+1}\sbs N_j,\q \cap N_j=F_1\cap F_2.$$
Put $F'_{j,2}=F_2\sms N_j$.
Let $\set{M_j}$ be a decreasing \seq of relatively open neighborhoods of $F_1$ \sth
$$\bar M_{j+1}\sbs M_j,\q \cap M_j=F_1,\q \bar M_j\cap F'_{j,2}=\ems.$$
Put $F'_{j,1}:=\bar M_j$.

Let $v_j$ be the largest solution dominated by $u$ and vanishing
on the complement of $F'_{j,1}\cup F'_{j,2}$:
$$\BAL
\bdw\sms(F'_{j,1}\cup F'_{j,2})&=\bdw\sms\big((F_1\cup F_2)\sms (N_j\sms \bar M_j)\big)\\
&=(\bdw\sms(F_1\cup F_2))\cup (N_j\sms \bar M_j).
\EAL$$
Furthermore, $(u- U_{\bar N_j\sms M_j})_+$ is a subsolution  which is dominated by $u$ and vanishes
on the complement of $F'_{j,1}\cup F'_{j,2}$.
Therefore $v_j$ satisfies
$$u\geq v_j\geq (u- U_{\bar N_j\sms M_j})_+,$$
which implies,
$$0\leq u-v_j\leq U_{\bar N_j\sms M_j}\leq U_{\bar N_j}.$$
By \eqref{UF-conv1}, $U_{\bar N_j}\downarrow U_{F_1\cap F_2}.$
Taking a converging \sseq $v_{j_i}\to v$ we obtain
$$0\leq u-v\leq U_{F_1\cap F_2}.$$

By the previous part of the proof there exist solutions $v_{j,1}$, $v_{j,2}$, whose boundary trace is supported in $F'_{j,1}$ and $F'_{j,2}$
respectively, \sth
$$v_j\leq v_{j,1}+v_{j,2}.$$
Taking a \sseq we may assume convergence of $\set{v_{j,1}}$ and $\set{v_{j,2}}$. Then  $u_i=\lim v_{j,i}$ has  boundary trace supported in $F_i$. Finally,
$$u\leq v+U_{F_1\cap F_2}\leq u_1+u_2 +U_{F_1\cap F_2}$$
and $\tr{\bdw}u_1$ is supported in $F_1$ while $\tr{\bdw}(u_2+U_{F_1\cap F_2})$ is supported in $F_2$. Since
$u-u_1$ is a subsolution dominated by the supersolution $u_2+U_{F_1\cap F_2}$ there exists a solution $w_2$ between them and we obtain
$$u\leq u_1+w_2$$
where $\tr{\bdw}w_2$ is supported in $F_2$.
\qed

The next theorem deals with some aspects of the generalized \bvp:
\begin{equation}\label{gbvp}\BAL
-\Gd u + g\circ u&=0,\q u\geq 0 \txt{in $\Gw$,}\\
\tr{\bdw}&=(\nu,F), \EAL\end{equation}

\nind where $F\sbs \bdw$ is a compact set and $\nu$ is a
(non-negative) Radon measure on $\bdw\sms F$.


\bth{reg+sing} Let $g\in\CG$ and assume that $g$ is convex and
satisfies the Keller-Osserman condition.

\nind{\sc Existence.}  The following set of conditions is  necessary
and sufficient  for existence of a solution $u$ of \eqref{gbvp}:

(i)   For every compact set $E\sbs \bdw\sms F$, the problem
\begin{equation}\label{meas-E}
  -\Gd u+g(u)=0 \txt{in $\Gw$,} u=\nu\chi\indx{E} \txt{on
  $\bdw$,}
\end{equation}
possesses a solution.

(ii) If $k_g(F)=F'$, then $F\sms F'\sbs\CS(V_\nu)$.

\smallskip

\nind When this holds,
\begin{equation}\label{u<V+U}
V_\nu\leq  u\leq V_\nu + U_F.
\end{equation}
Furthermore if $F$ is a removable set then \eqref{gbvp} possesses
exactly one solution.

\medskip
\nind{\sc Uniqueness.} Given a compact set $F\sbs \bdw$, assume that
\begin{equation}\label{uniq1}
  U_{E} \text{ is the unique solution with trace
 $(0,k_g(E))$}
\end{equation}
for every compact $E\sbs F$. Under this assumption:

\nind{\rm (a)} If $u$ is a solution of \eqref{gbvp} then
\begin{equation}\label{u>max(V,U)}
\max(V_\nu,U_F)\leq u\leq V_\nu + U_F.
\end{equation}

\nind{\rm (b)} Equation \eqref{eq-q} possesses at most one solution satisfying
\eqref{u>max(V,U)}.

\nind{\rm (c)} Condition \eqref{uniq1} is necessary and sufficient
in order that \eqref{gbvp} posses at most one solution.

\medskip

\nind{\sc Monotonicity.}

 \nind{\rm (d)} Let $u_1,u_2$
be two positive solutions of \eqref{eq1} with boundary traces
$(\nu_1,F_1)$ and $(\nu_2,F_2)$ respectively. Suppose that $F_1\sbs
F_2$ and that $\nu_1\leq \nu_2\chi\indx{F_1}=:\nu'_2$.
If \eqref{uniq1} holds for $F=F_2$ then $u_1\leq u_2$.

\es

\nind\Proof  First assume that there exists a solution $u$ of
\eqref{gbvp}. By \rth{reg-trace} condition (i) holds. \Consy $V_\nu$
is well defined by \eqref{Vnu}.

Since $V_\nu\leq u$ the function $w:=u-V_\nu$ is a subsolution of
\eqref{eq1}. Indeed, as $g$ is convex and $g(0)=0$ we have
\begin{equation}\label{g-convex}
  g(a)+g(b)\leq g(a+b) \forevery a,b\in \BBR_+.
\end{equation}
Therefore
$$0=-\Gd w +(g(u)-g(V_\nu)\geq -\Gd w +g(w).$$
 By \rth{ssol-1}, as $g$ satisfies
the Keller-Osserman condition, there exists a solution $\bar w$ of
\eqref{eq1} which is the smallest solution dominating $w$.

By \rth{ssol-3}, the traces of $w$ and $\bar w$ are equal on
$A=\CR(u)\sbs \CR(\bar w)$. Clearly the trace of $w$ on $\CR(u)$ is
zero. The definitions of $V_\nu$ and $\bar w$ imply,
\begin{equation}\label{ineqA}
 \max(V_\nu,\bar w)\leq u\leq V_\nu+\bar w.
\end{equation}
Therefore
$$\CS(\bar w)\cup\CS(V_\nu)=\CS(u).$$
In addition, as $\bar w$ has trace zero in $\bdw\sms F$, it follows,
 by the definition of the maximal function, that
 $$\bar w\leq U_F \txt{and \consy} \CS(\bar w)\sbs k_g(F).$$
These observations imply that condition (ii) must hold. Inequality
\eqref{u<V+U} follows from \eqref{ineqA} and this inequality implies
that if $F$ is a removable set then \eqref{gbvp} possesses exactly
one solution.

Now we assume that conditions (i) and (ii) hold and prove existence
of a solution. The function $V_\nu$ is well defined and $V_\nu +
U_F$ is a supersolution of \eqref{eq1} whose boundary trace is
$(\nu,F)$. Therefore, by \rth{ssol-3}, the largest solution
dominated by it has the same  boundary trace, i.e. solves
\eqref{gbvp}.

Next assume that condition \eqref{uniq1} is satisfied.
It is obvious that
\eqref{uniq1} is necessary for uniqueness. In addition, \eqref{uniq1}
implies that $U_F\leq u$  and \consy \eqref{u<V+U}
implies \eqref{u>max(V,U)}. It is also clear that (b) implies the sufficiency part of (c).

Therefore it remains to
prove statements (b) and (d).  Let $u$ be the smallest solution
dominating the subsolution $\max(V_\nu,U_F)$ and let $v$ be the
largest solution dominated by $V_\nu +U_F$.

To establish (b) we must show that $u=v$. By \eqref{u>max(V,U)} $ v-u\leq V_\nu$.
In addition the subsolution $v-u$ has trace zero on $\bdw\sms F$. Therefore
\begin{equation}\label{v-u}
    v-u\leq \min(V_{\nu}, U_F).
\end{equation}

 Let $\set{N_k}$
be a decreasing \seq of open sets converging to $F$ \sth
$N_{k+1}\Subset N_k$. Assuming for a moment that $\nu$ {\em is a
finite measure,} the trace of $V_\nu$ on $N_k$ is
$\nu_k:=\nu\chi\indx{N_k}$ and it tends to zero as $k\tin$.
Therefore, in this case,
$$\min(V_{\nu}, U_F)\leq V_{\nu_k}\to 0$$
and hence $u=v$. Of course this also implies uniqueness (statement (c))
in the case where $\nu$ is a finite measure.

In the general case we argue as follows. Let $v_k$ be the unique
solution with  boundary trace $(\nu'_k,\bar N_k)$ where
$\nu'_k=\nu(1-\chi\indx{\bar N_k})$. By taking a \sseq if necessary,
we may assume that $\set{v_k}$ converges to a solution $v'$. By
\eqref{u>max(V,U)},
$$\max(V_{\nu'_k}, U_{\bar N_k})\leq v_k\leq V_{\nu'_k}+ U_{\bar N_k}$$
and, by the previous part of the proof, $v_k$ is the largest solution dominated by
$V_{\nu'_k}+ U_{\bar N_k}$. We claim that if $w$ is a solution of \eqref{eq-q} then
\begin{equation}\label{Vnu'k}
  V_\nu\leq w\leq V_\nu+U_F\Lra w\leq V_{\nu'_k}+ U_{\bar N_k}.
\end{equation}
Indeed,
$$w\leq V_\nu+U_F\Lra w\leq V_{\nu'_k}+V_{\nu_k}+U_F\Lra w\leq V_{\nu'_k}+U_{\bar N_k}+U_F\Lra
w\leq V_{\nu'_k}+2U_{\bar N_k}.$$
Thus
$$0\leq w-V_{\nu'_k}\leq 2U_{\bar N_k}$$
which implies
$$w-V_{\nu'_k}\leq U_{\bar N_k},$$
 because any solution (or subsolution) dominated by $2U_{\bar N_k}$ is also dominated by $U_{\bar N_k}$.

Hence $v_k\geq v$ and \consy $v'\geq v$.

 By \eqref{UF-conv1} $U_{\bar N_{k}}\downarrow U_F$ and by definition $V_{\nu'_k}\uparrow V_\nu$. Therefore
$$\max(V_{\nu}, U_{F})\leq v'\leq V_{\nu}+ U_{F}.$$
Since $v$ is the largest solution dominated by $V_{\nu}+ U_{F}$ and $v\leq v'$ it follows that $v=v'$.


Let $u_k$ be the unique solution with boundary trace
$(\nu'_k,k_g(F))$. By \eqref{u>max(V,U)},
$$\max(V_{\nu'_k}, U_{k_g(F)})\leq u_k\leq V_{\nu'_k}+ U_{k_g(F)}.$$
Since $u_k\leq u$ and $\set{u_k}$ increases (because $\set{V_{\nu'_k}}$ increases)
it follows that $u'=\lim u_k\leq u$. Furthermore,
$$\max(V_{\nu}, U_{k_g(F)})\leq u'\leq V_{\nu}+ U_{k_g(F)}.$$

If \eqref{gbvp} possesses a solution then condition (ii) holds. Therefore for any solution $w$ of \eqref{eq-q}
$$ \max(V_{\nu}, U_{k_g(F)})\leq w\Lra
\max(V_{\nu}, U_{F})\leq w.$$
Hence $\max(V_{\nu}, U_{F})\leq u'$ and, as $u'\leq u$ we conclude that $u'=u$.

Finally, for every $\ge>0$,
$$(1-\ge)V_{\nu'_k}+ \ge U_{k_g(F)}\leq u_k$$
and \consy
$$\BAL
&v_k-u_k\leq V_{\nu'_k}+ U_{\bar N_k}-\big((1-\ge)V_{\nu'_k}+ \ge U_{k_g(F)})\big)=\\
&U_{\bar N_k}-(1-\ge) U_{k_g(F)}+\ge V_{\nu'_k}\leq U_{\overline{N_k\sms F}}+U_F-(1-\ge) U_{k_g(F)}+\ge V_{\nu'_k}\leq\\
&\ge(U_F +V_{\nu'_k})
\to \ge(U_{F}+V_{\nu}).
\EAL$$
This implies $u_k=v_k$ and hence $u=v$. This establishes statement (b) and hence the sufficiency in (c).

Finally we establish monotonicity.  Let $v_i$  be the unique solution of \eqref{eq-q} with boundary trace
$(\nu_i,F_i)$, (i=1,2). Then $v_i$ is the
largest solution dominated by
$V_{\nu_i} +U_{F_i}$ (i=1,2). The argument  used in proving \eqref{Vnu'k} yields
\begin{equation}\label{monotonicity}
  V_{\nu_1}\leq w\leq V_{\nu_1}+U_{F_1}\Lra w\leq V_{\nu_2}+ U_{F_2}.
\end{equation}
This implies $v_1\leq v_2$.
\qed

\mysection{Equation with power nonlinearity in a Lipschitz domain }
In this section we study the trace problem and the associated boundary value problem for equation
\begin{equation}\label{eq-q}
-\Gd u+\abs u^{q-1}u=0
\end{equation}
in a Lipschitz bounded domain  $\Gw$ and $q>1$. The main difference between the smooth cases and the Lipschitz case is the fact that the notion of critical exponent is pointwise.  If $G$ is any domain in $\BBR^N$ we denote
\begin{equation}\label{U(G)}
\CU(G):=\left\{\text{ the set of solutions (\ref{eq-q}) in }G\right\}.
\end{equation}
and $\CU_{+}(G)=\{u\in \CU(G):u\geq 0\text{ in } G\}$. Notice that any solution is at least $C^3$ in $G$ and any positive solution is $C^\infty$.
The next result is proved separately by Keller \cite{Ke} and Osserman \cite{Oss}.

\bprop{OK} Let $q>1$, $\Gw\subset\BBR^N$ be any domain and $u\geq\in C(\Gw)$ be a weak solution of
\begin{equation}\label{OKineq}
-\Gd u+Au^q\leq B\quad\text{in }\Gw.
\end{equation}
for some $A>0$ and $B\geq 0$. Then there exist $C_{i}(N,q)>0$ ($i=1,2$) such that
\begin{equation}\label{OKineq1}
u(x)\leq C_{1}\left(\myfrac{1}{\sqrt A\dist (x,\prt\Gw)}\right)^{2/(q-1)}
+C_{2}\left(\myfrac{B}{A}\right)^{1/q}\quad\forall x\in\Gw.
\end{equation}
\es
For a solution of (\ref{eq-q}) in $\Gw$ which vanishes on the boundary except at one point, we have a more precise estimate.

\bprop{OK2} Let $q>1$, $\Gw\subset\BBR^N$ be a bounded Lipschitz domain, $y\in\prt\Gw$ and  $u\in \CU_{+}(\Gw)$ is continuous in $\overline\Gw\setminus\{y\})$ and vanishes on $\prt\Gw\setminus\{y\}$. Then there exists $C_{3}(N,q,\Gw)>0$ and $\ga\in (0,1]$ such that
\begin{equation}\label{OKeq}
u(x)\leq C_{3}\left(\dist (x,\prt\Gw)\right)^{\ga}\abs{x-y}^{-2/(q-1)-\ga}
\quad\forall x\in\Gw.
\end{equation}
Furthermore $\ga=1$ if $\Gw$ is a $W^{2,s}$ domain with $s>N$.
\es
\Proof By translation we can assume that $y=0$. Let $\tilde u$ be the extension of $u_{+}$ by zero outside $\overline\Gw\setminus\{0\}$. Then it is a subsolution of (\ref{eq-q}) in $\BBR^N\setminus\{0\}$ (see \cite{GV} e.g.). Thus
$$\tilde u(x)\leq C_{1}|x|^{-2/(q-1)}\quad\forall x\neq 0,
$$
and, with the same estimate for $u_{-}$, we derive
\begin{equation}\label{OKeq1}\abs{u(x)}\leq C_{1}|x|^{-2/(q-1)}\quad\forall x\in \Gw.
\end{equation}
Next we define the transformation $T_{k}$ ($k>0$) by $T_{k}[u](x)=k^{-2/(q-1)}u(k^{-1}x)$, valid for any $x\in \Gw_{k}=k\Gw$. Then $u_{k}:=T_{k}[u]$ satisfies the same equation as $u$ in $\Gw_{k}$, is continuous in $\overline\Gw_{k}\setminus\{0\}$ and vanishes on $\prt\Gw_{k}\setminus\{0\}$. Then
$$u_{k}(x)\leq C_{1}|x|^{-2/(q-1)}\quad\forall x\in \Gw_{k},
$$
thus, by elliptic equation theory in uniformly Lipschitz domains, (which is the case if $k\geq 1$)
$$\norm{u_{k}}_{C^{\ga}(\Gw_{k}\cap (B_{7/4}\setminus B_{5/4}))}\leq C\norm{u_{k}}_{L^\infty(\Gw_{k}\cap (B_{2}\setminus B_{1}))}=C_{2}.
$$
This implies
$$|u(k^{-1}x')-u(k^{-1}z')|\leq C_{2}k^{-2/(q-1)-\ga}|x'-z'|^\ga
\quad\forall (x,z)\in \Gw_{k}\ti\Gw_{k}:5/4\leq |x'|,|z'|\leq 7/4.$$
Let $(x,z)$ in $\Gw\ti\Gw$ close enough to 0. First, if $5/7\leq |x|/|z|\leq 7/5$ there exists $k\geq 1$ such that
$5/4\leq |kx|,|kz|\leq 7/4$. Then
$$|u(x)-u(z)|\leq C_{3}|x|^{-2/(q-1)-\ga}|x-z|^{\ga}.
$$
If we take in particular $x$ such that $z=\rm{Proj}_{\prt\Gw}(x)$ satisfies the above restriction, we derive
$$u(x)\leq C_{3}|x|^{-2/(q-1)-\ga}\left(\dist (x,\prt\Gw)\right)^{\ga}.
$$
Because $\Gw$ is Lipschitz, it is easy to see that there exists $\gb\in (0,1/2)$ such that whenever
$\dist (x,\prt\Gw)=\abs{x-\rm{Proj}_{\prt\Gw}(x)}\leq\gb|x|$, there holds
$$5/7\leq |x|/\abs{\rm{Proj}_{\prt\Gw}(x)}\leq 7/5.
$$
Next we suppose $\abs{x-\rm{Proj}_{\prt\Gw}(x)}>\gb|x|$. Then, by the Keller-Osserman estimate,
$$u(x)\leq C|x|^{-2/q-1)-\ga}|x|^\ga\leq C\gb^{-\ga}|x|^{-2/q-1)-\ga}\abs{x-\rm{Proj}_{\prt\Gw}(x)}^\ga,
$$
which is (\ref{OKeq}). If we assume that $\prt\Gw$ is $W^{2,s}$, with $s>N$, then we can perform a change $W^{2,s}$ of coordinates near $0$ with transforms $\prt\Gw\cap B_{R}(0)$ into $\BBR^N_{+}\cap B_{R}(0)$ and the equation into
\begin{equation}\label{OKeq2}
-\sum_{i,j}\myfrac{\prt}{\prt x_{i}}\left(a_{ij}\myfrac{\prt \tilde u}{\prt x_{j}}\right)+|\tilde u|^{q-1}\tilde u=0,\quad\text{in }\BBR^N_{+}\cap B_{R}(0)\setminus\{0\},
\end{equation}
where the $a_{ij}$ are the partial derivatives of the coordinates and thus belong to $W^{1,s}(B_{R)}$. By developping, $\tilde u$ satisfies
$$
-\sum_{i,j}a_{ij}\myfrac{\prt^2 \tilde u}{\prt x_{i}\prt x_{j}}
-\sum_{j}b_{j}\myfrac{\prt \tilde u}{\prt x_{j}}+|\tilde u|^{q-1}\tilde u=0.
$$
Notice that, since $s>N$, the  $a_{ij}$ are continuous while the $b_{i}$ are in
$L^s$. The same regularity holds uniformly for the rescaled form of $\tilde u_{k}:=T_{k}[\tilde u]$. By the Agmon-Douglis-Nirenberg estimates $\tilde u_{k}$ belongs to $W^{2,s}$. Since $s>N$, $\tilde u$ satisfies an uniform $C^1$ estimates, which implies that we can take $\ga=1$.
\qeda

\subsection{Analysis in a cone}
\indent The removability question for solutions of (\ref{eq-q}) near
the vertex  of a cone has been studied in \cite{FV}, and we recall
this result below.

 If we look for separable solutions of (\ref{eq-q}) under the form
 $u(x)=u(r,\gs)=r^\gb\gw(\gs)$, where $(r,\gs)\in \BBR^+\ti S^{N-1}$ are the spherical coordinates,
 one finds immediately $\gb=-2/(q-1)$ and $\gw$ is a solution of
\begin{equation}\label{NEVP}
-\Gd'\gw-\gl_{_{N,q}}\gw+\abs\gw^{q-1}\gw=0
\end{equation}
on $S^{N-1}$ with
\begin{equation}\label{NEVP1}
\gl_{_{N,q}}=\myfrac{2}{q-1}\left(\myfrac{2q}{q-1}-N\right).
\end{equation}
Thus, a solution of (\ref{eq-q}) in the cone
$C_{_{S}}=\{(r,\gs):r>0,\gs\in S\subset S^{N-1}\}$, vanishing on
$\prt C_{_{S}}\setminus\{0\}$, has the form
$u(r,\gs)=r^{-2/(q-1)}\gw(\gs)$ if and only if $\gw$ is a solution
of (\ref{NEVP}) in $S$ which vanishes on $\prt S$. The next result
\cite[Prop 2.1]{FV} gives the  the structure of the set of positive
solutions of  (\ref{NEVP}).
\bprop{ex-uni}Let $\gl_{_{S}}$ be the first eigenvalue of the Laplace-Beltrami operator $-\Gd'$
in $W^{1,2}_{0}(S)$. Then \smallskip

\noindent (i) If $\gl_{_{S}}\geq \gl_{_{N,q}}$ there exists no
solution  to (\ref{NEVP}) vanishing on $\prt S$.\smallskip

\noindent (ii) If $\gl_{_{S}}< \gl_{_{N,q}}$ there exists a unique
positive solution $\gw=\gw_{_{S}}$ to (\ref{NEVP}) vanishing on
$\prt S$. Furthermore $S\subset S'\Longrightarrow \gw_{_{S}}\leq
\gw_{_{S'}}$. \es

The following is a consequence of \rprop{ex-uni}.
\bprop{remov}{\rm\cite{FV}} Assume $\Gw$ a bounded domain with a
purely conical part with vertex $0$, that is
$$\Gw\cap B_{r_{0}}(0)=C_{_{S}}\cap B_{r_{0}}(0)=\{x\in\cap B_{r_{0}}(0)\setminus\{0\}:x/\abs x\in S\}\cup\{0\}
$$
and that $\prt\Gw\setminus \{0\}$ is smooth. Then, if $\gl_{_{S}}\geq \gl_{_{N,q}}$, any solution $u\in \CU(\Gw)$ which is continuous in
$\overline\Gw\setminus\{0\}$ and vanishes on $\prt\Gw\setminus \{0\}$ is identically $0$.
\es

\noindent \Remark If $S\subset S^{N-1}$ is a domain and $\gl_{_{S}}$
the first eigenvalue of the Laplace-Beltrami operator $-\Gd'$ in
$W^{1,2}_{0}(S)$ we denote by  $\tl\ga_{_{S}}$ and $\ga_{_{S}}$ the
positive root and the absolute value of the negative root
respectively, of the equation

$$X^2+(N-2)X-\gl_{_{S}}=0.$$
Thus
\begin{equation}\label{alpha}\BAL
\tl\ga_{_{S}}&=\myfrac{1}{2}\left(2-N+\sqrt{(N-2)^2+4\gl_{_{S}}}\right),\\
\ga_{_{S}}&=\myfrac{1}{2}\left(N-2+\sqrt{(N-2)^2+4\gl_{_{S}}}\right).
\EAL\end{equation} It is straightforward that
$$\gl_{_{S}}\geq \gl_{_{N,q}}\Longleftrightarrow\ga_{_{S}}\geq \myfrac{2}{q-1},
$$
and, in case of equality, the exponent $q=q_{_{S}}$ satisfies $q_{_{S}}=1+2/\ga_{_{S}}$.

In subsection 6.2 we compute the Martin kernel $K$ and the first
eigenfunction $\gr$ of $-\Gd$  for cones with $k$-dimensional edge.
In particular, if $k=0$ and $C_S$ is the cone with vertex at the
origin and 'opening' $S\sbs S^{N-1}$, we have
\begin{equation}\label{K-ro-CS}
 K^{C_S}(x,0)=|x|^{-\ga_S}\gw_{_S}(\gs),\q
 \gr(x)=|x|^{\tl\ga_S}\gw_{_S}(\gs).
\end{equation}

Combining the  removability result with the admissibility condition
\rth{admissible}, we obtain the following.
 \bth{admiss} The problem
 \begin{equation}\label{tr-at-vertex}\BAL
   -&\Gd u+|u|^{q-1}u=0 \txt{in $C_S$,} \\
   &u\in C(\bar C_S\sms \{0\}),\q
    u=0 \txt{on $\prt C_S\sms\{0\}$}
\EAL \end{equation}
 possesses a non-trivial solution if and only if
 $$1<q<q\indx{S}=1+2/\ga_{_{S}}.$$

 Under this condition the following statements hold.

 \nind{\rm (a)} For every $k\neq0$ there exists a unique solution $v_k$ of
 \eqref{eq-q} with boundary trace $k\gd_0$. In addition we have
 \begin{equation}\label{vk/v1}
   v_k/v_1(x)\to k \q\text{uniformly as $ x\to 0$}.
 \end{equation}

 \nind{\rm (b)} Equation \eqref{eq-q} possesses a unique solution $U$ in $C_S$ \sth $\CS(U)=\{0\}$
 and its trace on $\prt C_S\sms\{0\}$ is zero. This solution satisfies
 \begin{equation}\label{sing-est1}
   |x|^{\frac{2}{q-1}}U(x)= U(x/|x|)=\gw_{S}(x/|x|)
 \end{equation}
 and
\begin{equation}\label{sing-est2}
   U=v_\infty:=\lim_{k\tin}v_k.
 \end{equation}
\es
\Proof (a) By \eqref{K-ro-CS},
$$\int_{C_{S}\cap B_{1}}K^q(x,0)\gr(x)\,dx
\leq C\myint{0}{1}r^{\tilde\ga_{_{S}}-q\ga_{_{S}}+N-1}dr<\infty,$$
since
$$\tilde\ga_{_{S}}-q\ga_{_{S}}+N-1=1-(q-1)\ga_{_{S}}>-1.$$
Thus $q$ is admissible for $C_S\cap B_1$ at $0$. By
\rth{admissible}, for every $k\in \BBR$, there exists a unique
solution of \eqref{eq-q} with boundary trace $k\gd_0$.

Observe that, for every $a,j>0$, $\tl v_j(x):=a^{2/(q-1)}v_j(ax)$ is
a solution of \eqref{eq-q} in $C_S$. This solution has boundary
trace $k\gd_0$ where $k=a^{2/(q-1)}j$. Because of uniqueness, $\tl
v_j=v_k$. Thus
\begin{equation}\label{vk--vj}
 v_k(x)=a^{2/(q-1)}v_j(ax),\q k=a^{2/(q-1)}j.
\end{equation}
This implies \eqref{vk/v1}.

\medskip
\nind(b) Let $w$ be a solution in $C_S$ \sth $\CS(w)=\{0\}$ and its
trace on $\prt C_S\sms\{0\}$ is zero. We claim that
\begin{equation}\label{w>v_infty}
  w\geq v_\infty:=\lim{k\tin}v_k.
\end{equation}
Indeed, for every $S'\Subset S$, $k>0$,
$$\int_{aS'}w\,d\gw_a\to\infty,\q \limsup\int_{aS'}v_kd\gw_a<\infty \txt{as $a\to 0$}$$
where $d\gw_a$ denotes the harmonic measure for a bounded Lipschitz domain $\Gw_a$ \sth $aS'\sbs \bdw_a$ and $\Gw_a\uparrow C_S$.
Therefore, using the classical Harnack inequality up to the boundary, $w/v_k\to\infty$ as $|x|\to 0$ in $C_{S'}$.
In addition, either by  Hopf's maximum principle (if $S$ is smooth) or by the boundary Harnack principle (if $S$ is merely Lipschitz), $$c^{-1}v_1\leq w\leq cv_1  \txt{in $C_{S\sms S'}$.}$$
 This inequality together
 with \eqref{vk--vj} yields,
 $$c^{-1}v_k\leq w\leq cv_k \txt{in $C_{S\sms S'}$}$$
 with $c$ independent of $k$. Therefore $c^{-1}v_k\leq w$ in $C_S$. If $1/c>k/cj>1$ then $\frac{k}{j} v_{j}\leq v_k\leq cw$ and \consy $v_j<w$. Here we used the fact that $\frac{k}{j} v_{j}$ is a subsolution with boundary trace $k\gd_0$.

 Let $U_{0}$ be the maximal solution with trace $0$ on $\prt C_S\sms \{0\}$ and singular boundary point at $0$. Then
 $$U_{0}(x)=a^{2/(q-1)}U_{0}(ax) \forevery a>0, \;x\in C_S,$$
 because $a^{2/(q-1)}U_{0}(ax)$ is again a solution which dominates every solution with trace $0$ on
 $\prt C_S\sms \{0\}$ and singular boundary point at $0$. Hence,
 \begin{equation}\label{selfsim}
 U_{0}(x)=|x|^{-2/(q-1)}U_{0}(x/|x|)=|x|^{-2/(q-1)}\gw_{S}(x/|x|).
 \end{equation}
The second equality follows from  the uniqueness part in
\rprop{ex-uni} since the function $x\to U_0(x/|x|)$ is continuous in
$\bar S$ and vanishes on $\prt S$.

Inequality \eqref{w>v_infty} implies that $v_\infty$ is the minimal
positive solution \sth $\CS(w)=\{0\}$ and its trace on $\prt
C_S\sms\{0\}$ is zero.
  Using this fact we prove in the same way that $v_\infty$ satisfies
  $$v_\infty(x)=|x|^{-2/(q-1)}v_\infty(x/|x|)=|x|^{-2/(q-1)}\gw_{S}(x/|x|).$$
   This implies \eqref{sing-est2} and the uniqueness in statement (b).
   \qed

In the next theorem we describe the precise asymptotic behavior of
 solutions in a conical domain with mass concentrated at the vertex.

\bth{cone1} Let $C_{_{S}}$ be a cone with vertex  $0$ and opening
$S\subset S^{N-1}$ and assume that $1<q<q_{_S}=1+2/\ga_{_S}$. Denote
by $\gf\indx{S}$ the first eigenfunction of $-\Gd'$ in
$W^{1,2}_{0}(S)$ normalized by $\max \gf_{_{S}}=1.$  Then the
function
$$\Gf_S=x^{-\ga_{_S}}\gf_{_{S}}(x/\abs x),$$
with $\ga\indx{S}$ as in \eqref{alpha}, is harmonic in $C_S$ and
vanishes on $\prt C_S\sms \{0\}$. Thus there exists $\gg>0$ \sth the
boundary trace of $\Gf_S$ is the measure $\gg\gd_0$. Put
$\Gf_1:=\rec{\gg}\Gf_S.$

Let $r_{0}>0$ and denote $\Gw_S=C_S\cap B_{r_{0}}(0)$. For every
$k\in \BBR$, let $u_k$ be the unique solution of \eqref{eq-q} in
$\Gw$ with boundary trace $k\gd_0$. Then
\begin{equation}\label{lim-k}
u_{k}(x)=k\Gf_1(x)(1+o(1))\quad\text{as }x\to 0.
\end{equation}

\nind\/If $v_k$ is the unique solution of \eqref{eq-q} in $C_S$ with
boundary trace $k\gd_0$ then
\begin{equation}\label{uk--vk}
u_k/v_k\to1 \txt{and} v_k/(k\Gf_1)\to 1 \quad\text{as }x\to 0.
\end{equation}

\nind\/The function $u_\infty=\lim_{k\tin} u_k$  is the unique
positive solution of \eqref{eq-q} in $\Gw_S$ which vanishes on $\prt
\Gw_S\sms \{0\}$ and is strongly singular at $0$ (i.e., $0$ belongs
to its singular set). Its asymptotic behavior at $0$ is given by,
\begin{equation}\label{lim-infi}
u_{\infty}(x)= |x|^{-\frac{2}{q-1}}\gw_{S}(x/|x|)(1+o(1))\quad
\txt{as}x\to0.
\end{equation}
\es
\Proof {\it Step 1: Construction of a fundamental solution. } Put
\begin{equation}\label{tilde-Gf}
\Gf(x)=\abs {x}^{-\ga_{_S}}\gf_{_{S}}(x/\abs x),\q  \tl\Gf(x)=\abs
{x}^{\tl\ga_{_S}}\gf_{_{S}}(x/\abs{x})
\end{equation}
with $\ga\indx{S}$, $\tl\ga\indx{S}$ as in \eqref{alpha}.
 Then $\Gf$ and $\tl\Gf$ are  harmonic
in $C_{_{S}}$, $\Gf$ vanishes on $\prt C_{_{S}}\setminus\{0\}$ and
$\tl\Gf$ vanishes on $\prt C_{_{S}}$. Furthermore, since
$q<1+2/\ga_{_S}$,
$$\int_{C_S\cap B_1(0)} \Gf^q\gr dx<\infty.$$
Therefore the boundary trace of $\Gf$ is a bounded measure
concentrated at the vertex of $C_S$, which means that the trace is
$\gg\gd_0$ for some $\gg>0$. (Here $\gd_0$ denotes the Dirac measure
on $\prt C_S$ concentrated at the origin.)

 The function
$$\Psi(x)=\rec{\gg}(\Gf(x)-r_0^{\tilde\ga_{_{S}}-\ga_{_{S}}}\tl\Gf(x))
$$
is harmonic and  positive in $\Gw_{_{S}}$ and  vanishes on $\prt
\Gw_{_{S}}\setminus\{0\}$. Its boundary trace is  $\gd_0$.

\medskip

\noindent{\it Step 2: Weakly singular behaviour. } By
\rth{admissible} , for any $k\geq 0$, there exists a unique function
$u_{k}\in L^q_{\gr}(\Gw_{_{S}})$ with trace $k\gd_0$ and by
\eqref{meas2}
 \begin{equation}\label{weak1}
 u_{k}(x)=k\Psi(x)-\BBG[\abs {u_{k}}^q].
\end{equation}

Since $\abs x^{\ga_{_{S}}}u_{k}$ is bounded, we set
$$v(t,\gs)=r^{\ga_{_{S}}}u_{k}(r,\gs),\quad t=-\ln r.
$$
Then $v$ satisfies
 \begin{equation}\label{weak2}v_{tt}+(2\ga_{_{S}}+2-N)v_{t}+\gl_{_{S}}v+\Gd'v-e^{(\ga_{_{S}}(q-1)-2)t}\abs v^{q-1}v=0
\end{equation}
in $D_{S,t_{0}}:=[t_{0},\infty)\ti S$ (with $t_{0}:=-\ln r_{0}$) and vanishes on $[t_{0},\infty)\ti \prt S$. Since $ 0\leq u_{k}(x)\leq k\Psi(x)$, $v$ is uniformly bounded, and, since $\ga_{_{S}}(q-1)-2<0$, $v(t,.)$ is uniformly bounded in
$C^\ga(\overline S)$ for some $\ga\in (0,1)$. Furthermore, $\nabla' v(t,.)$ (by definition $\nabla'$ is the covariant gradient on $S^{N-1}$) is bounded in $L^2(S)$, independently of $t$. Set
$$y(t)=\myint{S}{}v(t,\gs)\gf_{_{S}}dV(\gs),\quad F(t)=\myint{S}{}(\abs v^{q-1}v)(t,\gs)\gf_{_{S}}dV(\gs).$$
From (\ref{weak2}), it follows
$$\myfrac{d}{dt}\left(e^{(2\ga_{_{S}}+2-N)t}y'\right)=e^{((q+1)\ga_{_{S}}-N)t}F,
$$
where $dV$ is the volume measure on $S^{N-1}$. By (\ref{alpha}), $\gg:=2\ga_{_{S}}+2-N>0$, then
$$y'(t)=e^{-\gg(t-t_{0})}y'(t_{0})+e^{-\gg t}\myint{t_{0}}{t}
e^{((q+1)\ga_{_{S}}-N)s}F(s)ds,
$$
and
$$\abs {y'(t)}\leq c_{1}e^{-\gg(t-t_{0})}+c_{2}e^{(\ga_{_{S}}(q-1)-2)t}.
$$
This implies that there exists $k^*\in\BBR^+$ such that
 \begin{equation}\label{weak2'}\lim_{t\to\infty}y(t)=k^*.
  \end{equation}
 Next we use the fact that the following Hilbertian decomposition holds
$$L^2(S)=\oplus_{k=1}^\infty ker(-\Gd'-\gl_{k}I)
$$
where $\gl_{k}$ is the $k$-th eigenvalue of $-\Gd'$ in  $W^{1,2}_{0}(S)$ (and $\gl_{_{S}}=\gl_{1}$). Let $\tilde v$ and $\tilde F$ be the projections of $v$ and  $\abs v^{q-1}v$ onto $ker(-\Gd'-\gl_{_{S}}I)^{\perp}$.
Since
 \begin{equation}\label{weak3}\tilde v_{tt}+(2\ga_{_{S}}+2-N)\tilde v_{t}+\gl_{_{S}}\tilde v+\Gd'\tilde v-e^{(\ga_{_{S}}(q-1)-2)t}\tilde F=0
 \end{equation}
we obtain, by multiplying by $\tilde w$ and integrating on $S$,
$$V''+(2\ga_{_{S}}+2-N)V'-(\gl_{2}-\gl_{_{S}})V+e^{(\ga_{_{S}}(q-1)-2)t}\Gf\geq 0,
$$
where $V(t)=\norm{\tilde v(t,.)}_{L^2(S)}$ and $\Gf(t)=\norm{\tilde F(t,.)}_{L^2(S)}$. The associated o.d.e.
$$z''+(2\ga_{_{S}}+2-N)z'-(\gl_{2}-\gl_{_{S}})z+e^{(\ga_{_{S}}(q-1)-2)t}\Gf=0,
$$
admits solutions under the form
$$z(t)=a_{1}e^{-\gm_{1}t}+a_{2}e^{\gm_{2}t}+d(t)e^{(\ga_{_{S}}(q-1)-2)t}
$$
where $-\gm_{1}$ and $\gm_{2}$ are respectively the negative and the positive roots of
$$X^2+(2\ga_{_{S}}+2-N)X-(\gl_{2}-\gl_{_{S}})=0,$$
and $|d(t)|\leq c\Gf$ if $\ga_{_{S}}(q-1)-2\neq-\gm_{1}$, or  $|d(t)|\leq ct^{1}\Gf$  if $\ga_{_{S}}(q-1)-2=-\gm_{1}$. Applying the maximum principle to (\ref{weak3}), we derive
 \begin{equation}\label{weak4}
 \norm{\tilde v(t,.)}_{L^2(S)}\leq  \norm{\tilde v(t_{0},.)}_{L^2(S)}
 e^{-\gm_{1}(t-t_{0})}+d(t)e^{(\ga_{_{S}}(q-1)-2)t}\forevery t\geq t_{0}.
  \end{equation}
  By the standard elliptic regularity results in Lipschitz domains \cite{GT}, we obtain from (\ref{weak4}), for any $t>t_{0}+1$,
   \begin{equation}\label{weak5}
 \norm{\tilde v(t,.)}_{C^{\ga}(S)}\leq  c_{1}\norm{\tilde v}_{L^2((t-1,t+1)\ti S)}
 +c_{2}\norm{e^{(\ga_{_{S}}(q-1)-2)s}\tilde F}_{L^\infty((t-1,t+1)\ti S)},
  \end{equation}
  for some $\ga\in (0,1]$ depending of the regularity of $\prt S$. Thus
  \begin{equation}\label{weak6}\norm{\tilde v(t,.)}_{C^{\ga}(S)}\leq ce^{-\gm_{1}t}+
  c'te^{(\ga_{_{S}}(q-1)-2)t}.
  \end{equation}
  Combining (\ref{weak2'}) and (\ref{weak6}) we obtain that
    \begin{equation}\label{weak7}
    \abs x^{\ga_{_{S}}}u_{k}(x)-k^*\gf_{_{S}}(x/|x|)\to 0\quad \text{as }x\to 0
 \end{equation}
 in $C^\ga(S)$. Furthermore $0\leq k^*\leq k$.

 \bigskip

\noindent{\it Step 3: Identification of $k^*$. }

Let $\{\Gw_n\}$ be a \Lip exhaustion of $\Gw_S$ and denote by
$\gw_n$ (resp. $\gw$) the harmonic measure on $\bdw_n$ (resp.
$\bdw\indx{S}$). By \rprop{trace=limit}
$$\lim_{n\to\infty}\int_{\bdw_n}u_k\,d\gw_n=k.$$
On the other hand, by \eqref{weak7},
$$u_k/(k^*|x|^{-\ga\indx{S}}\gf\indx{S})\to 1 \txt{as}x\to0.$$
Hence
$$\BAL\lim_{n\to\infty}\int_{\bdw_n}u_k\,d\gw_n&=k^*\lim_{n\tin}
\int_{\bdw_n}|x|^{-\ga\indx{S}}\gf\indx{S}\,d\gw_n\\
&=k^*\gg\lim_{n\tin} \int_{\bdw_n}\Gf_1\,d\gw_n=k^*\gg. \EAL$$
Thus
\begin{equation}\label{k*gg}
   k=k^*\gg.
\end{equation}
This and \eqref{weak7} imply \eqref{lim-k}.

Further,
$$u_k\leq v_k\leq k\Gf_1$$
since $\Gf_1$ is harmonic in $C_S$. Therefore \eqref{lim-k} implies
\eqref{uk--vk}.

\bigskip

\noindent{\it Step 4: Study when $k\to\infty$. }  By \rth{admiss},
equation \eqref{eq-q} possesses a unique solution $U$ in $C_S$ \sth
$U=0$ on $\prt C_S\sms\{0\}$ and $U$ has strong singularity at the
vertex, i.e., $0\in \CS(U)$. By \eqref{sing-est1} and
\eqref{sing-est2} this solution satisfies
\begin{equation}\label{sing-est3}
   U=v_\infty:=\lim_{k\tin}v_k=|x|^{-\frac{2}{q-1}}\gw\indx{S}.
 \end{equation}

Let $V$ be the maximal solution in $\Gw_S$ vanishing on
$\bdw\indx{S}\sms\{0\}$. Its extension by zero to $C_S$ is a
subsolution and \consy, $V\leq U$.

Let $w$ be the unique solution of \eqref{eq-q} in $\Gw_S$ \sth $w=U$
on $\prt\Gw_S\cap B_{r_0}(0)$ and $w=0$ on the remaining part of the
boundary. Then $w<U$ so that $U-w$ is a subsolution of \eqref{eq-q}
in $\Gw_S$ which vanishes on $\bdw\indx{S}\sms \{0\}$. Therefore
$U-w\leq V$. Thus
\begin{equation}\label{UVw}
 U-w\leq V\leq U \txt{and} U/V\to 1 \txt{as} x\to0.
\end{equation}

\nind\/{\em Assertion 1.} If $u$ is a solution of \eqref{eq-q} in
$\Gw_S$ \sth
$$u=0 \txt{on} \bdw_{S}\sms \{0\} \txt{and} u/U\to 1\txt{as}x\to0$$

\nind\/then $u=V$.

\medskip
By \eqref{UVw} $u/V\to 1$ as $x\to0$. Therefore, by a standard
application of the maximum principle, $u=V$.

 Let $u$ be an arbitrary positive solution in $\Gw_S$ vanishing
on $\bdw\indx{S}\sms \{0\}$. Denote by $u^*$ its extension by zero
to $C_S$. Then $u^*$ is a subsolution and, by \rth{ssol-1}, there
exists a solution $\bar u$ of \eqref{eq-q} in $C_S$ which is the
smallest solution dominating $u^*$.
 The solution $\bar u$ can be obtained from $u^*$ as follows. Let $\{r_n\}$ be a
 sequence decreasing to zero, $r_1<r_0$, and denote
 $$D_n=C_S\sms B_{r_n}(0),\q h_n=u^*\lfloor\indx{\prt D_n}.$$
Let $w_n$ be the solution of \eqref{eq-q} in $D_n$ \sth $w_n=h_n$ on
the boundary. Then $\{w_n\}$ increases and
\begin{equation}\label{baru=lim}
\bar u=\lim w_n.
\end{equation}

If $u$ has strong singularity at the origin then, of course, the
same is true \wrto $\bar u$ and \consy, by \rth{admiss},
\begin{equation}\label{baru=U}
  \bar u=U.
\end{equation}
In the  the remaining part of the proof we assume only
\eqref{baru=U} and show that this implies $u=V$.

Let $z$ be the solution of \eqref{eq-q} in $\Gw_S$ \sth $z=U$ on
$\prt\Gw_S\cap\prt B_{r_0}$ and $0$ on $\prt\Gw_S\cap\prt C_S$. Then
$u+z$ is a supersolution  in $\Gw_S$. Let
$$\Gw_n=\Gw_S\sms B_{r_n}(0)=D_n\cap B_{r_0}(0).$$
The trace of $u+z$ on $\prt\Gw_n$ is given by
$$f_n=\begin{cases} U &\text{on }\prt\Gw_n\cap\prt B_{r_0}\\ h_n+z
&\text{on }\prt\Gw_n\sms \prt B_{r_0}.\end{cases}$$

\nind\/Since  $U=\bar u\geq u^*$ we have $f_n\geq h_n$. Therefore,
if $\tl w_n$ is the solution of \eqref{eq-q} in $\Gw_n$ \sth $\tl
w_n=f_n$ on the boundary then
$$ w_n\leq \tl w_n\leq u+z \txt{in}\Gw_n.$$
Hence, by \eqref{baru=lim},
$$U\leq u+z.$$

Since $z\to 0$ as $x\to 0$, it follows that

$$\limsup U/u\leq 1 \txt{as}x\to0.$$
Since $u<V$, \eqref{UVw} implies that
$$\liminf U/u\geq 1 \txt{as} x\to0.$$
Therefore $U/u\to1$ as $x\to0$ and \consy, by Assertion 1, $u=V$.
This  proves the uniqueness stated in the last part of the theorem
and \eqref{UVw} implies \eqref{lim-infi}. \qed

\bcor{uinfty} Suppose that $u$ is a positive solution of
\eqref{eq-q} in $\Gw_S$ which vanishes  on $\bdw_S\sms \{0\}$ and
\begin{equation}\label{sup=infty}
   \sup_{\Gw_S}|x|^{\ga\indx{S}}u=\infty.
\end{equation}
Then $u=u_\infty$. \es

\Proof Let $\bar u$ be as in \eqref{baru=lim}. Since $\bar u\geq u$
it follows that
$$ \sup_{\Gw_S}|x|^{\ga\indx{S}}\bar u=\infty.$$
By \rth{admiss} $\bar u=U$. The last part of the proof shows that
$u=u_\infty$. \qed
\medskip

As a consequence of \rth{cone1}  we obtain the classification of
positive solutions of (\ref{eq-q}) in conical domains with isolated
singularity located at the vertex. In the case of a half space  such
a classification was obtained in \cite{GV}.

\bth{cone2} Let $C_{_{S}}$ be as in \rth{cone1},
$\Gw_{s}=C_{_{S}}\cap B_{r_{0}}(0)$ for some $r_{0}>0$ and
$1<q<q_{_S}=1+2/\ga_{_S}$. If $u\in C(\bar\Gw_{s}\setminus\{0\})$ is
a  positive solution of (\ref{eq-q}) vanishing on $\prt C_{_{S}}\cap
B_{r_{0}}(0)\sms\{0\}$, the following alternative holds:

\smallskip
Either\\
 \noindent (i)  $\limsup_{x\to 0}\abs
x^{-\tilde\ga_{_{S}}}u(x)<\infty$ and thus $u\in C(\bar\Gw_{s})$.

\smallskip
or\\
 \noindent (ii)  there exist $k>0$ such that (\ref{lim-k})
holds

\smallskip
or\\
 \noindent (iii) (\ref{lim-infi}) holds.\es
\Proof Let  $u_{\ge}$ be the solution of (\ref{eq-q}) in $\Gw_{_{S,\ge}}=\Gw_{_{S}}\setminus B_{\ge}(0)$ with boundary data $u$ on $\Gw_{_{S,\ge}}\cap
\prt B_{\ge}(0)$ and zero on $\prt \Gw_{_{S,\ge}}\setminus \prt B_{\ge}(0)$. Then
$$0\leq u_{\ge}\leq u\leq u_{\ge}+Z(x)
\forevery x\in \Gw_{_{S,\ge}},
$$
where $Z$ is harmonic in $\Gw_{_{S}}$, vanishes on $\prt
\Gw_{_{S}}\setminus \prt B_{r_{0}}(0)$ and coincides with $u$ on
$C_{_{S}}\cap \prt B_{r_{0}}(0)$. Furthermore
$0<\ge<\ge'\Longrightarrow u_{\ge}\leq u_{\ge'}$ in
$\Gw_{_{S,\ge'}}$. Thus $u_{\ge}$ converges, as $\ge\to 0$, to a
solution $\tilde u$ of (\ref{eq-q}) which vanishes on
$\prt\Gw_{_{S}}\setminus\{0\}$ and satisfies
\begin{equation}\label{E1-0}
0\leq \tilde u(x)\leq u(x)\leq \tilde u(x)+Z(x)
\forevery x\in \Gw_{_{S}}.
\end{equation}
If
\begin{equation}\label{E2}\limsup_{x\to 0}\abs x^{\ga_{_{S}}}\tilde u(x)<\infty,
\end{equation}
it follows from \rth{cone1}-Step 2, that there exists $k^*\geq 0$ such that
\begin{equation}\label{E3}\tilde u(x)=k^*\abs x^{-\ga_{_{S}}}\gf_{_{S}}(x/|x|)(1+o(1))\quad\text{as }x\to 0.
\end{equation}
If $k^*>0$ then $u$ satisfies (ii). If $k^*=0$, it is
straightforward to see that, for any $\ge>0$,  $\tilde u(x)\leq
\ge\abs x^{-\ga_{_{S}}}$. Thus
\begin{equation}\label{E4}u(x)\leq Z(x)=c\abs x^{\tilde\ga_{_{S}}}\gf_{_{S}}(x/|x|)(1+o(1))\quad\text{as }x\to 0,
\end{equation}
by standard expansion of harmonic functions at $0$.

Finally, if
\begin{equation}\label{E5}\limsup_{x\to 0}\abs x^{\ga_{_{S}}}\tilde u(x)=\infty,
\end{equation}
then, by \rcor{uinfty}, $\tilde u=u_\infty$  and \consy, by
\rth{cone1}, $\tl u$ -- and therefore $u$ -- satisfies
\eqref{lim-infi}. \qeda

\medskip

\subsection{Analysis in a Lipschitz domain}
In a general Lipschitz bounded domain tangent planes have to be replaced by asymptotic cones, and  these  asymptotic cones can be inner or outer.
\bdef {critcone} Let $\Gw$ be a bounded Lipschitz domain and $y\in
\prt\Gw$. For $r>0$,  we denote by $\CC^{I}_{y,r}$ (resp.
$\CC^{O}_{y,r}$) the set of all open cones $C_{s,y}$ with vertex at
$y$ and smooth opening $S\subset\prt B_{1}(y)$ such that
$C_{s,y}\cap B_{r}(y)\subset\Gw$ (resp. $\Gw\cap B_{r}(y)\subset
C_{s,y}$). Further we denote
\begin{equation}\label{r-cones}
C^{I}_{y,r}:=\bigcup \left\{C_{S,y}:C_{S,y}\in
\CC^{I}_{y,r}\right\},\q C^{O}_{y,r}:=\bigcap
\left\{C_{S,y}:C_{S,y}\in \CC^{O}_{y,r}\right\}
\end {equation}
and 
\begin{equation}\label{lim-cone}
C^{I}_{y}:=\bigcup_{r>0}C^{I}_{y,r}, \q C^{O}_{y}:=\bigcap_{r>0}
C^{O}_{y,r}.
\end {equation}
The cone  $C^{I}_{y}$ (resp. $C^{O}_{y}$) is called the limiting
inner cone (resp. outer cone) at $y$. Finally we denote
\begin{equation}\label{Syr}\BAL S^{I}_{y,r}:=&C^{I}_{y,r}\cap \prt
B_{1}(y), &S^{O}_{y,r}:=&C^{O}_{y,r}\cap \prt B_{1}(y),\\
S^{I}_{y}:=&C^{I}_{y}\cap \prt B_{1}(y),&S^{O}_{y}:=&C^{O}_{y}\cap
\prt B_{1}(y).\EAL
\end{equation}
 \es

\Remark In this definition, we identify $\prt B_{1}(y)$ with the
manifold $S^{N-1}$.  Notice that the following monotonicity holds
\begin{equation}\label{mono}0<s<r\Longrightarrow\left\{\BA {l}C^{I}_{y,r}\subset C^{I}_{y,s}\\[2mm]
C^{O}_{y,s}\subset C^{I}_{y,r}. \EA\right.
\end {equation}

\bdef {critexp} If $C_S$ is a cone with vertex $y$ and opening $S$
and if $\gl_S$ is the first eigenvalue of $-\Gd'$ in
$W^{1,2}_{0}(S)$, we denote
\begin{equation}\label{exp1}
\ga_{_{S}}=\myfrac{1}{2}\left(N-2+\sqrt{(N-2)^2+4\gl_{_{S}}}\right),
\txt{and} q_{_{S}}=1+2/\ga_{_{S}}.
\end {equation}
Thus $q_{_{S}}$ is the critical value for the cone $C_S$ at its
vertex. \es

\noindent \Remark As $r\mapsto S^{I}_{y,r}$ is nondecreasing, it
follows that $r\mapsto \gl\indx{S^{I}_{y,r}}$ is nonincreasing and
\consy $r\mapsto q\indx{S^{I}_{y,r}}$ is nondecreasing. It is
classical that
\begin{equation}\label{limint-I}
\lim_{r\to 0}\gl_{_{S^{I}_{y,r}}}=\gl_{_{S^{I}_{y}}}.
\end{equation}

A similar observation holds with respect to $ S^{O}_{y,r}$ if we
interchange the terms  `nondecreasing' and `nonincreasing'. In
particular
\begin{equation}\label{limint-O}
\lim_{r\to 0}\gl_{_{S^{O}_{y,r}}}=\gl_{_{S^{O}_{y}}}.
\end{equation}
In view of \eqref{exp1} we conclude that,
\begin{equation}\label{lim-qry}
\lim_{r\to 0}q_{_{S^{I}_{y,r}}}=q_{_{S^{I}_{y}}}, \q \lim_{r\to
0}q_{_{S^{O}_{y,r}}}=q_{_{S^{O}_{y}}}.
\end{equation}

\medskip

We also need the following notation:

\bdef{q*y} Let $\Gw$ be a bounded Lipschitz domain. For every
compact set $E\sbs \bdw$ denote,
\begin{equation}\label{q*y}
q^*_{E}=\lim_{r\to 0}\inf\left\{q\indx{S^{I}_{z,r}}:z\in\prt\Gw,
\;\dist(z,E)<r\right\},
\end{equation}
If $E$ is a singleton, say $\{y\}$, we replace $q^*_E$ by $q^*_y$.
\es

\Remark For a cone $C_S$ with vertex $y$,   $q^*_y\leq q_{S}$.
However if $C_S$ is contained in a half space then $q^*_y= q_{S}$.
On the other hand, if $C_S$ strictly contains a half space then
$q^*_y< q_S$.

If $\Gw$ is the complement of a bounded convex domain then, for
every $y\in \bdw$,
\begin{equation}\label{concave-dom}
   q^*_y=(N+1)/(N-1)
\end{equation}
Indeed $q_{c,y}\geq (N+1)/(N-1)$. But for $\BBH_{N-1}$-a.e. point
$y\in\bdw$ there exists a tangent plane and \consy
$q_{c,y}=(N+1)/(N-1)$. This readily implies \eqref{concave-dom}.

Since $\Gw$ is \Lip, there exists $r_\Gw>0$ \sth,  for every $r\in
(0,r_\Gw)$ and every $z\in \bdw$, there exists a cone $C$ with
vertex at $z$ \sth $C\cap B_r(z)\sbs \bar\Gw$. Denote
$$a(r,y):=\inf\left\{q\indx{S^{I}_{z,r}}:z\in\prt\Gw\cap
B_{r}(y)\right\} \forevery r\in (0,r_\Gw),\;y\in \bdw.$$ Then,
\begin{equation}\label{max-min}\BAL
q^*_E:=&\lim_{r\to 0}\inf\{a(r,y):y\in E\}\\
\leq &\inf \,\{\lim _{r\to 0}a(r,y):y\in E\}=\inf\,\{ q^*_y:y\in
E\}. \EAL\end{equation}

Indeed, the monotonicity  of the function $r\mapsto
q_{_{S^{I}_{y,r}}}$ (for each fixed $y\in \bdw$) implies
\begin{equation}\label{max-min1}
q^*_y=\lim_{r\to 0}a(r,y)=\sup_{0<r<r_\Gw}a(r,y).
\end{equation}
As
$$q^*_E=\lim_{r\to0}\inf\{a(r,y):y\in E\}$$
inequality \eqref{max-min} follows immediately from
\eqref{max-min1}.

Finally we observe that, if $E$ is a compact subset of $\bdw$ then
\begin{equation}\label{q*Er}
 (E)_r:=\{z\in\bdw:\dist(z,E)\leq r\}\Lra  q^*\indx{{(E)}_r}\uparrow q^*_E \txt{as} r\downarrow0.
\end{equation}

In order to deal with boundary value problems in a general \Lip
domain $\Gw$ we must study the question of q-admissibility of
$\gd_y$, $y\in \bdw$. This question is addressed in the following:

\bth{admi} If $y\in\prt\Gw$ and
$1<q<q_{_{S^{I}_{y}}}:=1+2/\ga_{_{S^{I}_{y}}}$ then
\begin{equation}\label{adm0}
\myint{\Gw}{}K^q(x,y)\gr(x)dx<\infty.
\end {equation}
Furthermore, if $E$ is a compact subset of $\bdw$ and   $1<q<q^*_E$
then, there exists $M>0$ such that,
\begin{equation}\label{adm1}
\myint{\Gw}{}K^q(x,y)\gr(x)dx\leq M \forevery y\in E.
\end {equation}
\es
\Proof  We recall some sharp estimates of the Poisson kernel due to
Bogdan \cite{Bog}. Set $\gk=1/2(\sqrt{1+K^2})$, where $K$ is the
Lipschitz constant of the domain, seen locally as the graph of a
function from $\BBR^{N-1}$ into $\BBR$. Let $x_{0}\in \Gw$ and set
$\phi(x):=G(x,x_{0})$. Then there exists $c_{1}>0$ such that for any
$y\in\prt\Gw$ and $x\in\Gw$ satisfying $|x-y|\leq r_{0}$, there
holds
\begin{equation}\label{adm2}
c_{1}^{-1}\myfrac{\phi (x)}{\phi^2(\xi)}|x-y|^{2-N}
\leq K(x,y)\leq c_{1}\myfrac{\phi (x)}{\phi^2(\xi)}|x-y|^{2-N},
\end {equation}
for any $\xi$ such that $B_{\gk|x-y|}(\xi)\subset\Gw\cap B_{|x-y|}(y)$. This implies
\begin{equation}\label{adm3}c_{2}^{-1}\myfrac{\phi^{q+1} (x)}{\phi^{2q}(\xi)}|x-y|^{(2-N)q}
\leq K^q(x,y)\gr(x)\leq c_{2}\myfrac{\phi^{q+1} (x)}{\phi^{2q}(\xi)}|x-y|^{(2-N)q}
\end {equation}
for some $c_{2}$ since $\phi$ and $\gr$ are comparable in $B_{r_{0}}(y)$, uniformly with respect to $y$ (provided we have chosen $r_{0}\leq \dist(x_{0},\prt\Gw)/2$. Let $C_{s,y}$ be a smooth cone with vertex at $y$ and opening $S:=C_{s,y}\cap \prt B_{1}(y)$, such that  $\overline C_{s,y}\cap \prt B_{r_{0}}(y)\subset\Gw$. We can impose to the point $\xi$ in inequality (\ref{adm2})  to be such that  $\xi/|\xi|:=\Xi_{0}\in S$, or, equivalently, such that $|\xi-y|\leq \gg\dist (\xi,\prt\Gw)$ for some $\gg>1$ independent of $\xi$, $|x-y|$ and $y$. Then, by Carleson estimate \cite[Lemma 2.4]{Ba} and
Harnack inequality, there exists $c_{5}$ independent of $y$ such that there holds
\begin{equation}\label{adm3+2}\myfrac{\gf(\xi)}{\gf(x)}\geq c_{3}\end {equation}
for all $x\in
\Gw\cap B_{r_{0}}(y)$ and all $\xi $ as above. Consequently, (\ref{adm3}) yields to
\begin{equation}\label{adm2+1}
K^q(x,y)\gr(x)\leq c_{4}\phi^{1-q}(\xi)|x-y|^{(2-N)q}.
\end {equation}
There exists a separable harmonic function $v$ in $C_{s,y}$ under the form
$$v (z)=|z-y|^{\ga_{_{S}}+2-N}\phi_{_{S}}((z-y)/|z-y|)
$$
where $\phi_{_{S}}$ is the first eigenfunction of $-\Gd'$ in
$W^{1,2}_{0}(S)$ normalized by  $\max \phi_{_{S}}=1$, $\gl_{_{S}}$ the corresponding eigenvalue and $\ga_{_{S}}$ is given by (\ref{alpha}). By the maximum principle,
\begin{equation}\label{adm3+1}v (z)\leq c_{5}\phi(z)\quad\forall
z\in C_{_{S,y}}\cap B_{r_{0}}(y).
\end {equation}
Therefore there exists $c_{6}>0$ such that
\begin{equation}\label{adm3+3}\phi(\xi)\geq c_{6}\abs{\xi-y}^{\ga_{_{S}}+2-N}.
\end {equation}
Because  $|x-y|\geq \abs{\xi-y}\geq \gk\abs{x-y}/2$, from the choice of $\xi$, it follows
\begin{equation}\label{adm4}
K^q(x,y)\gr(x)\leq \myfrac{c_{7}}{|x-y|^{(q-1)\ga_{_{S}}+N-2}}\quad\forall x\in \Gw\cap B_{r_{0}}(y).
\end{equation}
Clearly, if we choose $q$ such that
$1<q<q_{_{S^{I}_{y}}}:=1+2/\ga_{_{S^{I}_{y}}}$,  then
$q<1+2/\ga_{_{S^{I}_{r,y}}}$ for some $r$ small enough and we can
take $C_{S,y}=C^{I}_{y,r}$. Thus (\ref{adm0}) follows.

We turn to the proof of \eqref{adm1}. To simplify the notation we
assume that $q<q^*_{\bdw}$. The argument is the same in the case
$q<q^*_E$.

If we assume $q<\lim_{r\to
0}\inf\{q_{_{S^{I}_{z,r}}}:z\in\prt\Gw\}$, then for $\ge>0$ small
enough, there exists $r_{\ge}>0$ such that
$$ 0<r\leq r_{\ge}\Longrightarrow 1<q<\inf\{q_{_{S^{I}_{z,r}}}:z\in\prt\Gw\}-\ge\quad \forall 0<r\leq r_{\ge}.$$
Notice that the shape of the cone may vary, but, since $\prt\Gw$ is
Lipschitz there exists a fixed relatively open subdomain
$S^*\subset\prt B_{1}$ such that for any $y\in\prt\Gw$, there exists
an isometry $\CR_{y}$ of $\BBR^{N}$ with the property that
$\CR_{y}(\overline S^*)\subset S^{I}_{y,r}$ for all $0<r\leq
r_{\ge}$. Here we use the fact that  $r\mapsto S^{I}_{y,r}$ is
increasing when $r$ decreases. If we take $\xi$ such that
$\xi/|\xi|=\Xi_{0}\in \CR_{y}( S^*)$,  then the constants in Bogdan
estimate (\ref{adm2}) and Carleson inequality (\ref{adm3+2}) are
independent of $y\in\prt\Gw$ if we replace $r_{0}$ by
$\inf\{r_{\ge},r_{0}\}$. Hereafter we shall assume that $r_{\ge}\leq
r_{0}$. Set
$$v_{S} (t)=|t-y|^{\ga_{_{S}}+2-N}\phi_{_{S}}((t-y)/|t-y|)
$$
with $S=S^{I}_{y,r_{\ge}}$. Then $v_{S}$ is well defined in the cone
$C_{S,y}$ with vertex $y$ and opening $S$. Let
$$\Gs_{cr_{\ge}}:=\{t\in\Gw:\dist (t,\prt\Gw)=cr_{\ge}\}.$$ Because
$\prt\Gw$ is Lipschitz, we can choose $0<c<1$ such that $C_{S,y}\cap
\Gs_{cr_{\ge}}\subset B_{r_{\ge}}(z)$. Then we can compare $v_{S}$
and $\phi$ on the set $\Gs_{cr_{\ge}}$. It follows by maximum
principle that estimate (\ref{adm3+1}) is still valid with a
constant may depend on $r_{\ge}$, but not on $y$. Because
$$\min_{\CR_{y}(S^*)}\phi_{_{S^{I}_{y,r_{\ge}}}}\geq c_{8}
$$
where $c_{8}$ is independent of $y$, (\ref{adm3+3}) holds under the
form
\begin{equation}\label{adm5}\phi(\xi)\geq c_{6}\abs{\xi-y}^{\ga_{_{S^{I}_{y,r_{\ge}}}}+2-N},
\end {equation}
where, we recall it, $\xi$ satisfies $\xi/|\xi|\in\CR_{y}(S^*)$,  and is associated to any $x\in B_{r_{\ge}}(y)\cap\Gw$ by the property that $B_{\gk|x-y|}(\xi)\subset B_{|x-y|}(y)\cap\Gw$, and thus $|x-y|\geq \abs{\xi-y}\geq \gk\abs{x-y}/2$. Then
(\ref{adm4}) holds uniformly with respect to $y$, with $r_{0}$ replaced by $r_{\ge}$. This implies (\ref{adm1}).\qeda
\medskip

The next proposition partially complements \rth{admi}.

\bprop{admi''} Let $y\in\prt\Gw$ and $q>q_{_{S^{O}_{y}}}$. Then any
solution of (\ref{eq-q}) in $\Gw$ which vanishes on
$\prt\Gw\setminus\{0\}$ is identically $0$. \es

\Remark This proposition implies that, if $q>q_{_{S^{O}_{y}}}$,
\begin{equation}\label{adm1bis}
\myint{\Gw}{}K^q(x,y)\gr(x)dx=\infty.
\end {equation}
Otherwise $\gd_y$ would be admissible.

\medskip
\Proof  We consider a local outer smooth cone with vertex at $y$, $C_{2}$, such that
$\overline\Gw\cap B_{r_{0}}(y)\setminus\{0\}\subset C_{2}\cap B_{r_{0}}(y):=C_{2,r_{0}}$. We denote by $S^*=C_{2}\cap \prt B_{1}(y)$ its opening.
For $\ge>0$ small enough, we consider the doubly truncated cone $C^\ge_{2,r_{0}}=\cap C_{2,r_{0}}\setminus B_{\ge}(y)\}$ and the solution $v:=v_{\ge}$ to
\begin{equation}\label{eq-cone1}\left\{\BA {l}
-\Gd v+ v^{q}=0\quad\text{in }C^\ge_{2,r_{0}}\\
\phantom{-,\Gd v+^{q}}
v=\infty\quad\text{on }\prt B_{\ge}(y)\cap C_{2}\\
\phantom{-,\Gd v+^{q}}
v=\infty\quad\text{on }\prt B_{r_{0}}(y)\cap C_{2}\\
\phantom{-,\Gd v+^{q}}
v=0\quad\text{on }\prt C_{2}\cap B_{r_{0}}(y)\setminus \overline B_{\ge}(y),
\EA\right.\end{equation}
where $q\geq q_{_{S^*}}:=1+2/\ga_{_{S^*}}$, and $\ga_{_{S^*}}$ is expressed by (\ref{alpha}) with $S$ replaced by $S^*$. Then $v_{\ge}$ dominates in $C^\ge_{2,r_{0}}\cap\Gw$ any positive solution $u$ of (\ref{eq-q}) in $\Gw$ which vanishes on
$\prt\Gw\setminus\{0\}$. Letting $\ge \to 0$, $v_{\ge}$ converges to $v_{0}$ which satisfies
\begin{equation}\label{eq-cone2}\left\{\BA {l}
-\Gd v+ v^{q}=0\quad\text{in }C_{2,r_{0}}\\
\phantom{-,\Gd v+^{q}}
v=\infty\quad\text{on }\prt B_{r_{0}}\cap C_{2}\\
\phantom{-,\Gd v+^{q}}
v=0\quad\text{on }\prt C_{2}\cap B_{r_{0}}(y).
\EA\right.\end{equation}
Furthermore $u\leq v_{0}$ in $B_{r_{0}}\cap \Gw$. Because $q_{_{S^*}}$ is the critical exponent in $C_{2}$ , the singularity at $0$ is removable, which implies that $v(x)\to 0$ when $x\to 0$ in $C_{2}$. Thus $u_{+}(x)\to 0$ when $x\to 0$ in $\Gw$. Thus $u_{+}=0$. But we can take any cone with vertex $y$ containing $\Gw$ locally in $B_{r}(y)$ for $r>0$. This implies that for any $q> q_{_{S^{O}_{y}}}$, any solution of (\ref{eq-q}) which vanishes on $\prt\Gw\setminus\{0\}$ is non-positive. In the same way it is non-negative.
\qeda \medskip

\bdef{crtI}   If $y\in\prt\Gw$  we say that an exponent $q\geq 1$ is: \medskip

\noindent (i) Admissible at $y$ if
$$\norm{K(.,y)}_{L^q_{\gr}(\Gw)}<\infty,$$
and we set
$$q_{1,y}=\sup\{q> 1: q\text{ admissible at } y\}.
$$
\noindent (ii) Acceptable at $y$ if there exists a solution of
(\ref{eq-q}) with boundary trace $\gd_{y}$, and we set
$$q_{2,y}=\sup\{q> 1: q\text{ acceptable at } y\}.
$$

\noindent (iii) Super-critical at $y$ if any solution of (\ref{eq-q}) which is continuous in $\Gw\setminus\{0\}$ and vanishes on
$\prt\Gw\setminus\{0\}$ is identically zero, and we set
$$q_{3,y}=\inf\{q> 1: q\text{ super-critical at } y\}.
$$
\es

\bprop {ord}Assume $\Gw$ is a bounded Lipschitz domain and
$y\in\prt\Gw$. Then
\begin{equation}\label{ord1}
q_{_{S^{I}_{y}}}\leq q_{1,y}\leq q_{2,y}\leq q_{3,y}\leq q_{_{S^{O}_{y}}}.
\end{equation}
If $1<q<q_{2,y}$ then, for any real $a$  there exists exactly one
solution of \eqref{eq-q}  with boundary trace $\gg\gd_{y}$.
 \es
\Proof It follows from \rth{admi} that $q_{_{S^{I}_{y}}}\leq
q_{1,y}$ and from \rprop{admi''} that  $q_{3,y}\leq
q_{_{S^{O}_{y}}}.$ It is clear from the definition and
\rth{admissible} that $q_{1,y}\leq q_{2,y}\leq q_{3,y}$. Thus
\eqref{ord1} holds.

Now assume that $q<q_{2,y}$ so that there exists a solution $u$ with
boundary trace $\gd_y$. By the maximum principle $u>0$ in $\Gw$. If
$a\in (0,1)$ then $au$ is a subsolution of \eqref{eq-q} with
boundary trace $a\gd_y$ and $au<u$. Therefore by \rcor{ssol-2} II,
the smallest solution dominating $au$ has boundary trace $a\gd_y$.
If $a>1$ then $au$ is a supersolution and the same conclusion
follows from \rcor{ssol-2} I. If $v_a$ is the (unique) solution of
\eqref{eq-q} with boundary trace $a\gd_y$ then $-v$ is the (unique)
solution  with boundary trace $-a\gd_y$. \qed

\bth{coin} Assume $y\in\prt\Gw$ is such that
$S^{O}_{y}=S^{I}_{y}=S$, let $\gl_{_{S}}$ be the first eigenvalue of
$-\Gd'$ in $W^{1,2}_{0}(S)$ and denote
\begin{equation}\label{qcy}
q_{_{c,y}}:=1+2/\ga_{_{S}}
\end{equation}
with $\ga\indx{S}$ as in \eqref{alpha}.
 Then
$q_{1,y}=q_{2,y}=q_{3,y}=q_{_{c,y}}$ and \\
{\rm (i)} if $1<q< q_{_{c,y}}$ then $\gd_y$ is admissible;\\
 {\rm (ii)} if $q>q_{_{c,y}}$ then the only solution of \eqref{eq-q}
in $\Gw$ vanishing on $\bdw\sms \{y\}$ is the trivial solution.\\
 {\rm (iii)} if $q=q_{_{c,y}}$ and $u$ is a solution of \eqref{eq-q}
 in $\Gw$ vanishing on $\bdw\sms\{y\}$ then
 \begin{equation}\label{uo(1)}
    u=o(1)|x-y|^{-\frac{2}{q-1}} \txt{as $x\to y$ in $\Gw$.}
\end{equation}

 \es

 \Remark We  know that, in the conical case, the conclusion of statement (ii) holds
 for
 $q=q_{c,y}$ as well. \Consy,  in a polyhedral domain $\Gw$, an
 isolated
 singularity at a point
 $y\in\bdw$ is removable if $q\geq q_c(y)$. We do not know if this
 holds in general \Lip domains.


\medskip
\Proof The above assertion, except for statement (iii), is an
immediate \cons of \rprop{ord}, \rdef {critexp} and the remark
following that definition.

It remains to prove (iii). We may assume that $u>0$. Otherwise we
observe that $|u|$ is a subsolution of \eqref{eq-q} and by
\rth{ssol-1}(ii) there exists a solution $v$ dominating it. It is
easy to verify that the smallest solution dominating $|u|$ vanishes
on $\bdw\sms \{y\}$.

For any $r>0$ let $u_r$ be the extension of $u$ by zero to
$D_r:=C_{S_r^O}\cap B_r(y)$. Thus $u_r$ is a subsolution in $D_r$,
$u_r\in C(\bar D_r\sms \{y\})$ and $u_r=0$  on $(\prt C_{S_r^O}\cap
B_r(y))\sms \{y\}$. The smallest solution above it, say $\tl u_r$ is
in $C(\bar D_r\sms \{y\})$ and $\tl u_r=0$  on $(\prt C_{S_r^O}\cap
B_r(y))\sms \{y\}$. By a standard argument this implies that there
exists a positive solution $\tl v_r$ in $D_r$ \sth $\tl v_r$
vanishes on $\prt D_r\sms \{y\}$ and
$$u_r\leq 2\tl v_r \txt{in} D_r.$$
 We extend
this solution by zero to the entire cone $C_{S_r^O}$, obtaining a
subsolution $\tl w_r$ and finally (again by \rth{ssol-1}(ii)) a
solution $w_r$ in $C_{S_r^O}$ which vanishes on $\prt C_{S_r^O}\sms
\{y\}$ and satisfies
$$u_r\leq 2w_r \txt{in} D_r.$$
Observe that $q_{_{S_r^O}}\downarrow q_{c,y}$  as $r\downarrow 0$.
If $q_{c,y}= q_{_{S_r^O}}$ for some $r>0$ then the existence of a
solution $w_r$ as above is impossible. Therefore we conclude that
$q_{c,y}< q_{_{S_r^O}}$ and therefore, by \rth{admiss}, there exists
a solution $v_{\infty,r}$ in $C_{S_r^O}$ \sth
$$v_{\infty,r}(x)=  |x-y|^{-\frac{2}{q-1}}\gw_{_{S_r^O}}((x-y)/|x-y|)
\forevery x\in C_{S_r^O}.$$
This solution is the maximal solution in
$C_{S_r^O}$ so that
$$w_r\leq v_{\infty,r} \txt{in} D_r.$$
But, since $q=q_{_{S^O}}$, it follows that $\gw_{_{S_r^O}}\to 0$ as
$r\to 0$. This implies \eqref{uo(1)}.

 \qeda

The next  result provides an important ingredient in the study of
general boundary value problems in \Lip domains.

\bth{sub} Assume that $q>1$, $\Gw$ is a bounded Lipschitz domain and
$u\in\CU_{+}(\Gw)$. If $y\in\CS(u)$ and $q<q^*_y$ then, for every
$k>0$,  the measure $k\gd_y$ is admissible and
\begin{equation}\label{u>uk}
  u\geq u_{k\gd_{y}} \forevery k\geq 0.
\end{equation}
\es

\Remark If $q>q^*_y$, \eqref{u>uk} need not hold. For instance,
consider the cone $C_S$ with vertex at the origin, \sth $S\sbs
S^{N-1}$ is a smooth domain and  $S^{N-1}\sms S$ is contained in an
open half space. Then $q_{c,0}>(N+1)/(N-1)$ while
$q_{c,x}=(N+1)/(N-1)$ for any $x\neq0$ on the boundary of the cone.
Thus $q^*(0)<q_{c,0}$. Suppose that $q\in (q^*_0, q_{c,0})$. Let $F$
be a closed subset of $\prt C_S$ \sth $0\in F$ but $0$ is a
$C_{2/q,q'}$-thin point of $F$. Let $u$ be the maximal solution in
$C_S$ vanishing on $\prt C_S\sms F$. Then $0\in \CS(u)$ but
\eqref{u>uk} does not hold for any $k>0$.

\medskip
\Proof Up to an isometry of $\BBR^N$, we can assume that $y=0$ and
represent $\prt\Gw$ near $0$ as the graph of  a Lipschitz function.
This can be done in the following way: we define the cylinder
$C'_{R}:=\{x=(x',x_{N}):x'\in B'_{R}\}$ where $B'_{R}$ is the
$(N-1)$-ball with radius $R$. We denote, for some $R>0$ and
$0<\gs<R$,
$$\prt\Gw\cap C'_{R}=\{x=(x',\eta(x')):x'\in B'_{R}\},
$$
and
$$\Gs_{\gd,\gs}=\{x=(x',\eta(x')+\gd):x'\in B'_{\gs}\},$$
 and assume that, if $0<\gd\leq R$,
$$\Gw^R_{\gd}=\{x=(x',x_{N}):x'\in B'_{R},\,\eta(x')<x_{N}<\eta(x')+R\}\subset\Gw.
$$
We can also assume that $\eta (0)=0$. Although the two harmonic measures in $\Gw$ and $\prt\Gw\cap C'_{R}$ differ, it follow by Dahlberg's result that there exists a constant $c>0$ such that, if $\gd<\gd_{0}\leq R/2$,
$$c^{-1}\gw_{\Gw}^{x_{0}}(E)\leq \gw_{\Gw^R_{\gd}}^{x_{0}}(E+\ge {\bf e}_{ N})\leq c\,\gw_{\Gw}^{x_{0}}(E),
$$
for any Borel set $E\subset \prt\Gw\cap C'_{\gd}$. Therefore, if we set
$$M_{\ge,\gs}=\myint{\Gs_{\ge,\gs}}{}u(x)d\gw^{x_{0}},(x),
$$
it follows that $\lim_{\ge\to 0}M_{\ge,\gs}=\infty$ since $0\in\CS(u)$. We can suppose that $\gs$ is small enough so that there exists
$\hat q\in (q,q^*_{y})$ and $M>0$ such that, for any $ p\in [1, \hat q]$
\begin{equation}\label{unif}
\myint{\Gw}{}K^p(x,z)\gr(x) dx\leq M\quad
\forall  z\in \prt\Gw\cap B_{\gs}.
\end{equation}
For fixed $k$ there exists $\ge=\ge(\gd)>0$ such that $M_{\ge,\gs}=k$.  There exists a uniform Lipschitz exhaustion $\{\Gw_{\ge}\}$ of $\Gw$ with the following properties: \smallskip

\noindent (i) $\Gw_{\ge}\cap C'_{R}\cap \{x=(x',x_{N}):a<x_{N}<b\}=\Gs_{\ge,R}$, for some fixed $a$ and $b$.\smallskip

\noindent (ii) The $\Gw_{\ge}$ and $\Gw$ have the same Lipschitz character $L$.\smallskip

\noindent It follows that the Poisson kernel
$K^{\Gw_{\ge}}$ in $\Gw_{\ge}$ respectively endows the same properties (\ref{unif}) as $K$ except $\Gw$ has to be replaced by $\Gw_{\ge}$, $\gr$ by $\gr_{\ge}:=\dist(.,\prt\Gw_{\ge}$ and $z$ has to belong to $\prt\Gw_{\ge}\cap B_{\gs}$. Next, we consider the solution $v=v_{\ge(\gs)}$ of
\begin{equation}\left\{\BA {l}
-\Gd v+v^q=0\quad\qquad\text{in }\Gw_{\ge}\\
\phantom{-\Gd +v^q}v=u\chi_{_{\Gs_{\ge,\gs}}}\quad\text{in }\prt\Gw_{\ge}
\EA\right.\end{equation}
By the maximum principle, $u\geq v$ in $\Gw_{\ge}$. Furthermore $v\leq \BBK^{\Gw_{\ge}}[u\chi_{_{\Gs_{\ge,\gs}}}]$. Let $\hat q=(q+\tilde q_{\gs})/2$ and $\gw\subset\Gw$ be a Borel subset. By  convexity
$$\myint{\gw}{}\left(\BBK^{\Gw_{\ge}}
[u\chi_{_{\Gs_{\ge,\gs}}}]\right)^{\hat q}\gr(x) dx\leq M\,M_{\ge,\gs}.
$$
Thus, by H\"older's inequality
$$\myint{\gw}{}\left(\BBK^{\Gw_{\ge}}
[u\chi_{_{\Gs_{\ge,\gs}}}]\right)^{ q}\gr(x) dx\leq
\left(\myint{\gw}{}\gr(x)dx\right)^{1-q/\hat q}\left(M\,M_{\ge,\gs}\right)^{q/\hat q}.
$$
By standard a priori estimates, $v_{\ge(\gs)}\to v_{0}$ (up to a subsequence) a.e. in $\Gw$, thus $v^q_{\ge(\gs)}\to v^q_{0}$. By Vitali's theorem and the uniform integrability of the
$\{v_{\ge(\gs)}\}$, $v_{\ge(\gs)}\to v_{0}$ in $L^q_{\gr}(\Gw)$. Because
$$v_{\ge(\gs)}+\BBG^{\Gw_{\ge}}[v^q_{\ge(\gs)}]
=\BBK^{\Gw_{\ge}}[u\chi_{_{\Gs_{\ge,\gs}}}]
$$
where $\BBG^{\Gw_{\ge}}$ is the Green operator in  $\Gw_{\ge}$,  and
$$\BBK^{\Gw_{\ge}}[u\chi_{_{\Gs_{\ge,\gs}}}]\to
M_{\ge,\gs}K(.,y)=kK(.,y)
$$
as $\gs\to 0$, it follows that $u\geq v_{0}$, and $v_{0}$ satisfies
$$v_{0}+\BBG^\Gw_{}[v^q_{0}]=kK(.,y).
$$
Then $v_{0}=u_{k\gd_{y}}$, which ends the proof.\qeda\medskip

\bcor{sub'} Let $\{y_{j}\}_{j=1}^n\subset\bdw$ be a set of points
\sth
\begin{equation}\label{y*CS}
 q<\inf \{q^*_{y_{j}}:j=1,...,n\}.
\end{equation}
Then, for any set of positive numbers $k_1,\cdots, k_n$, there
exists a unique solution $u_\mu$ of \eqref{eq-q} in $\Gw$ with
boundary trace $\mu=\sum_{j=1}^nk_j\gd_{y_j}$.

If $u\in \CU_{+}(\Gw)$ and $\{y_{j}\}_{j=1}^n\subset\CS(u)$ then
 $u\geq u_\mu$.
 \es

\Proof From \rth{sub}, $u\geq u_{k_{j}\gd_{y_{j}}}$ for any
$j=1,...,n$. Thus $u\geq \tilde u_{\{k\}} = \max(
u_{k_{j}\gd_{y_{j}}})$, which is a subsolution with boundary trace
$\sum_j k_j\gd_{y_j}$. But $\tilde v_{\{k\}}$, the solution with
boundary trace $\sum_{j}k_{j}\gd_{y_{j}}$ is the smallest solution
above $\tilde u_{\{k\}}$. Therefore the conclusion of the corollary
holds. \qeda\medskip

 As a \cons one obtains

 \bth{ex-sub} Let $E\sbs \bdw$ be a closed set and assume that $q<q^*_E$. Then,
for every $\mu\in \GTM(\Gw)$ \sth $\supp\mu\sbs E$ there exists a
(unique) solution $u_\mu$ of \eqref{eq-q} in $\Gw$ with boundary
trace $\mu$.

If $\{\mu_n\}$ is a \seq in $\GTM(\Gw)$ \sth $\supp\mu_n\sbs E$ and
$\mu_n\rightharpoonup\mu$ weak* then $u_{\mu_n}\to u_\mu$ locally
uniformly in $\Gw$.

If $u\in \CU_{+}(\Gw)$ and  $q<q^*_{_{\CS(u)}}$ then, for every
$\mu\in \GTM(\Gw)$ \sth $\supp\mu\sbs \CS(u)$,
\begin{equation}\label{umu<u}
 u_\mu\leq u.
\end{equation}
 \es

 \Proof Without loss of generality we assume that $\mu\geq)$. Let
$\{\gm_{n}\}$ be a \seq of measures on $\prt\Gw$ of the form
$$\gm_{n}=\sum_{j=1}^{k_{n}}a_{j,n}\gd_{y_{j,n}}
$$
where $y_{j,n}\in E$, $a_{j,n}>0$ and
$\sum_{j=1}^{k_{n}}a_{j,n}=\norm\gm$, such that
$\gm_{n}\rightharpoonup\gm$ weakly*. Passing to a \sseq if
necessary, $u_{\gm_{n}}\to v$ locally uniformly in $\Gw$. In order
to prove the first assertion it remains to show  that $v=u_{\gm}$.


If  $0<r$ is sufficiently small, there exists $\hat q_r\in
(q,q^*_{E})$ and $M_r>0$ such that, for any $p\in [1,\hat q_r]$ and
every $z\in \bdw$ such that $\dist(z,E)<r$, estimate (\ref{unif})
holds. It follows that the family of functions
$$\set{K(\cdot,z): z\in \bdw,\;\dist(z,E)<r}$$
is uniformly integrable in $L^q_\gr(\Gw)$ and \consy the family
$$\set{\BBK[\nu];\nu\in \GTM(\bdw),\;\norm{\nu}_{\GTM}\leq
1,\;\supp\nu\sbs \{z\in\bdw:\dist(z,E)<r\}}$$ is uniformly
integrable in $L^q_\gr(\Gw)$. By a standard argument (using Vitali's
convergence theorem) this implies that $v= u_\mu$. This proves the
first two assertions of the theorem.

The last assertion is an immediate consequence of the above together
with \rcor{sub'}. Indeed, if $E=\CS(u)$ then, by \rcor{sub'}, $u\geq
u_{\gm_{n}}$. Therefore $u\geq u_\mu$.

 \qeda\medskip


\bprop{critprop2} Let $y\in \prt\Gw$ and $1<q<q_{_{S^I_{y}}}$. Then there exists a maximal solution $u:=U_{y}$ of (\ref{eq-q}) such that $tr(U_{y})=(\{y\},0)$. It satisfies
\begin{equation}\label{mino1}
\liminf_{\tiny{\BA{l}x\to y\\
\frac{x-y}{|x-y|}\to\gs \EA}}|x-y|^{2/(q-1)}U_{y}(x)\geq \gw_{_{S^I_{y}}}(\gs),
\end {equation}
uniformly on any compact subset of $S^I_{y}$, where
$\gw_{_{S^I_{y}}}$ is the unique positive solution of
\begin{equation}\label{mino2}\left\{\BA {l}
-\Gd' \gw-\gl_{_{N,q}}\gw+|\gw|^{q-1}\gw=0\quad\text {in }
S^I_{y}\\[2mm]
\phantom{-\Gd' \gw-\gl_{_{N,q}}\gw+|\gw|^{q-1}}\gw=0 \quad\text {on
}\prt S^I_{y}, \EA\right.\end {equation} normalized by
$\gw(\gs_0)=1$ for some fixed $\gs_0\in S^I_y$.

For $r>0$ small enough, we denote by $\gw_{_{S^O_{y,r}}}$ the unique
positive solution of
\begin{equation}\label{maj1}\left\{\BA {l}
-\Gd' \gw-\gl_{_{N,q}}\gw+|\gw|^{q-1}\gw=0\quad\text {in }
S^O_{y,r}\\[2mm]
\phantom{-\Gd' \gw-\gl_{_{N,q}}\gw+|\gw|^{q-1}}\gw=0 \quad\text {on
}\prt S^O_{y,r}, \EA\right.\end {equation} normalized in the same
way. Then
\begin{equation}\label{maj2}
\limsup_{\tiny{\BA{l}x\to y\\
\frac{x-y}{|x-y|}\to\gs \EA}}|x-y|^{2/(q-1)}U_{y}(x)\leq
\gw_{_{S^O_{y,r}}}\left(\gs\right) .
\end{equation}
Finally, if $S^{O}_{y}=S^{I}_{y}=S$, then
\begin{equation}\label{maj2'}
\lim_{\tiny{\BA{l}x\to y\\
\frac{x-y}{|x-y|}\to\gs \EA}}|x-y|^{2/(q-1)}U_{y}(x)=
\gw_{_{S}}\left(\gs\right).
\end{equation}
\es
\Proof
We recall that $C^{I}_{y,r}$ (resp. $C^{O}_{y,r}$) is a r-inner cone (resp. r-outer cone) at $y$ with opening
$S^{I}_{y,r}\subset \prt B_{1}(y)$ (resp. $S^{O}_{y,r}\subset \prt B_{1}(y)$). This is well defined for a $r>0$ small enough so that $q<q_{_{S^I_{y,r}}}$. We denote by $\gw_{_{S^I_{y,r}}}$ the unique positive solution of
\begin{equation}\label{mino2r}\left\{\BA {l}
-\Gd' \gw-\gl_{_{N,q}}\gw+|\gw|^{q-1}\gw=0\quad\text {in }
S^I_{y,r}\\[2mm]
\phantom{-\Gd' \gw-\gl_{_{N,q}}\gw+|\gw|^{q-1}}\gw=0
\quad\text {on }\prt S^I_{y,r}.
\EA\right.\end {equation}
We construct  $U_{y}\in\CU_{+}(\Gw$), vanishing on $\prt\Gw\setminus\{y\}$ in the following way. For $0<\ge<r$, we denote by $v:=U_{y,\ge}$ the solution of
$$\left\{\BA {l}
-\Gd v+|v^{q-1}|v=0\quad\text{in }\Gw\setminus \overline B_{\ge}(y)
\\[2mm]\phantom{-\Gd v+|v^{q-1}|}
v=0\quad\text{in }\prt\Gw\setminus \overline B_{\ge}(y)
\\[2mm]\phantom{-\Gd v+|v^{q-1}|}
v=\infty\quad\text{in }\Gw\cap\prt B_{\ge}(y).
\EA\right.$$
Let $v:=V^{I}_{\ge}$ (resp. $v:=V^{O}_{\ge}$) be the solution of
$$\left\{\BA {l}
-\Gd v+|v^{q-1}|v=0\quad\text{in }C_{_{S^I_{y,r}}}\setminus
\overline B_{\ge}(y)\qquad  (\text{resp. }C_{_{S^O_{y,r}}}\setminus
\overline B_{\ge}(y))
\\[2mm]\phantom{-\Gd v+|v^{q-1}|}
v=0\quad\text{in }\prt C_{_{S^I_{y,r}}}\setminus \overline B_{\ge}(y)
\;\;\quad (\text{resp. }\prt C_{_{S^O_{y,r}}}\setminus \overline B_{\ge}(y))
\\[2mm]\phantom{-\Gd v+|v^{q-1}|}
v=\infty\quad\text{in }C_{_{S^I_{y,r}}}\cap\prt B_{\ge}(y)
\quad (\text{resp. }C_{_{S^O_{y,r}}}\cap\prt B_{\ge}(y)).
\EA\right.$$
Then there exist $m>0$ depending on $r$, but not on $\ge$, such that
\begin{equation}\label{s3}
V^{I}_{\ge}(x)-m\leq U_{y,\ge}(x)\leq V^{O}_{\ge}(x)+m
\end{equation}
for all $x\in C^I_{y,r}\setminus\{B_{\ge}(y)\}$ for the left-hand
side inequality, and $x\in\prt\Gw\cap
B_{r}(y)\setminus\{B_{\ge}(y)\}$ for the right-hand side one.  When
$\ge\to 0$,  $V^{I}_{\ge}$ converges to the explicit separable
solution $x\mapsto |x-y|^{-2/(q-1)}\gw_{_{S^I_{y,r}}}$ in
$C_{_{S^I_{y,r}}}$ (the positive cone with vertex generated by
$S^I_{y,r}$). Similarly $V^{O}_{\ge}$ converges to the explicit
separable solution $x\mapsto |x-y|^{-2/(q-1)}\gw_{_{S^O_{y,r}}}$ in
$C_{_{S^O_{y,r}}}$. Furthermore $\ge<\ge'\Longrightarrow
U_{y,\ge}\leq U_{y,\ge'}$. If $U_{y}=\lim_{\ge\to 0}\{U_{y,\ge}\}$,
there holds
\begin{equation}\label{s4}
|x-y|^{-2/(q-1)}\gw_{_{S^I_{y,r}}}(\frac{x-y}{|x-y|})-m\leq U_{y}(x)\leq
|x-y|^{-2/(q-1)}\gw_{_{S^I_{y,r}}}(\frac{x-y}{|x-y|})+m.
\end{equation}
These inequalities imply
\begin{equation}\label{mino3}
\liminf_{\tiny{\BA{l}x\to y\\
\frac{x-y}{|x-y|}\to\gs \EA}}|x-y|^{2/(q-1)}U_{y}(x)\geq \gw_{_{S^I_{y,r}}}(\gs),
\end {equation}
Inequality (\ref{maj2}) is obtained in a similar way. Since
$\lim_{r\to 0}\gw_{_{S^I_{y,r}}}=\gw_{_{S^I_{y}}}$ uniformly in
compact subsets of $S^I_y$ we also obtain (\ref{mino1}). If
$S^{O}_{y}=S^{I}_{y}=S$, then
$\gw_{_{S^I_{y}}}=\gw_{_{S^O_{y}}}=\gw_{_{S}}$, thus (\ref{maj2'})
holds. \qeda\medskip

\noindent \Remark Because $U_{y}$ is the maximal solution which
vanishes on $\prt\Gw\setminus\{y\}$, the function
$u_{\infty\gd_{y}}=\lim_{k\to \infty}u_{k\gd_{y}}$ also satisfies
inequality (\ref{maj2}). {\it We conjecture that $u_{\infty\gd_{y}}$
always satisfies estimate (\ref{mino1})}. This is true if the outer
and inner cone at $y$ are the same. In fact in that case we obtain a
much stronger result:


\bth{uniq} Assume $y\in\prt\Gw$ is such that $S^{O}_{y}=S^{I}_{y}=S$
and $q<q_{c,y}$.  Then $U_{y}=u_{\infty\gd_{y}}$. \es

\Proof Without loss of generality we can assume that $y=0$ and will
denote  $B_{r}=B_{r}(0)$ for $r>0$. Let $C^I_{r}$ (resp. $C^O_{r}$)
be a cone  with vertex $0$,  such that $\overline{C^I_{r}\cap
B_{r}}\setminus\{0\}\subset\Gw$ (resp. $\Gw\cap B_{r}\subset
C^O_{r}$). We recall that the characteristic exponents
$\ga_{_{S^I_{0}}}$ and $\ga_{_{S^O_{0}}}$ are defined according to
\rdef{critcone} and \rdef{critexp}. Since
$$\ga_{_{S^I_{0}}}=\lim_{r\to 0}\ga_{_{S^I_{0,r}}}=\lim_{r\to 0}\ga_{_{S^O_{0,r}}}=\ga_{_{S^O_{0}}}<2/(q-1),$$
we can choose $r$ such that
\begin{equation}\label {condi}
q\ga_{_{S^I_{0,r}}}-\ga_{_{S^O_{0,r}}}<2-(q-1)(\ga_{_{S^I_{0,r}}}-\ga_{_{S^O_{0,r}}}),
\end{equation}
and for simplicity, we set
$\ga_{_{S^I_{0,r}}}=\ga_{_{I}}$, $\ga_{_{S^O_{0,r}}}=\ga_{_{O}}$
and
$$\gg_{r}=\myfrac{q-1}{2+\ga_{_{O}}-q\ga_{_{I}}}.
$$
{\it Step 1.} We claim that there exists $c>0$ and $c^*>0$ such that, for any $m>0$
\begin{equation}\label {below1}
u_{m\gd_{0}}(x)\geq c^*m|x|^{-\ga_{_O}}\quad\forall x\in
B_{cm^{-\gg_{r}}}\cap C^I_{r}.
\end{equation}
Since $mK(.,0)$ is a super-solution for (\ref{eq-q}),
$$u_{m\gd}(x)\geq mK(x,0)-m^q\myint{\Gw}{}G(z,x)K^q(z,0)dz.
$$
If we assume that $x\in C^I_{r}\cap B_{r}$, then $\dist
(x,\prt\Gw)\geq \gth|x|$ for some $\gth>0$ since
$\overline{C^I_{r}\cap B_{r}}\setminus\{0\}\subset\Gw$.  Using
Bogdan's  estimate and Harnack inequality we derive
$$K(x,0)\geq c_{1}\myfrac{|x|^{2-N}}{G(x,x_{0})},
$$
for some fixed point $x_{0}$ in $\Gw$. But the Green function in
$\Gw\cap B_{r}$ is  dominated by the Green function in $C^O_{r}\cap
B_{r}$, thus $G(x,x_{0})\leq c_{2} |x|^{\tilde\ga_{_{O}}}$ where
$\tilde\ga_{_{O}}=2-N+\ga_{_{O}}$. This implies
\begin{equation}\label{below2}
K(x,0)\geq c_{3}|x|^{-\ga_{_{O}}}\quad\forall x\in C^I_{r}\cap B_{r}.
\end{equation}
Similarly (and it is a very rough estimate)
$$K(x,0)\leq c_{4}|x|^{-\ga_{_{I}}}\quad\forall x\in \Gw
$$
Because $G(x,z)\leq c_{5}|x-z|^{2-N}$, we obtain
$$\myint{\Gw}{}G(z,x)K^q(z,0)dz\leq c_{6}\myint{B_{R}}{}
|x-z|^{2-N}|z|^{-\ga_{_{I}}}dz.
$$
We write
$$\BA {l}\myint{B_{R}}{}
|x-z|^{2-N}|z|^{-q\ga_{_{I}}}dz=
\myint{B_{2|x|}}{}
|x-z|^{2-N}|z|^{-q\ga_{_{I}}}dz\\[4mm]
\phantom{\myint{B_{R}}{}
|x-z|^{2-N}|z|^{-q\ga_{_{I}}}dz-}
+\myint{B_{R}\setminus B_{2|x|}}{}
|x-z|^{2-N}|z|^{-q\ga_{_{I}}}dz.
\EA$$
But
$$\myint{B_{2|x|}}{}
|x-z|^{2-N}|z|^{-q\ga_{_{I}}}dz
=|x|^{2-q\ga_{_{I}}}
\myint{B_{2}(0)}{}|\xi-t|^{2-N}|t|^{-q\ga_{_{I}}}dt
$$
where $\xi=x/|x|$ is fixed. In the same way
$$\BA {l}\myint{B_{R}\setminus B_{2|x|}}{}
|x-z|^{2-N}|z|^{-q\ga_{_{I}}}dz
\leq \myint{B_{R}\setminus B_{2|x|}}{}
|z|^{2-N-q\ga_{_{I}}}
|x|^{2-q\ga_{_{I}}}dz\\[4mm]
\phantom{\myint{B_{R}\setminus B_{2|x|}}{}
|x-z|^{2-N}|z|^{-q\ga_{_{I}}}dz}
\leq |x|^{2-q\ga_{_{I}}}\myint{B_{R/|x|}\setminus B_{2}}{}
|t|^{2-N-q\ga_{_{I}}}dt\\[4mm]
\phantom{\myint{B_{R}\setminus B_{2|x|}}{}
|x-z|^{2-N}|z|^{-q\ga_{_{I}}}dz}
\leq c_{7} |x|^{2-q\ga_{_{I}}}\myint{2}{R/|x|}
s^{1-q\ga_{_{I}}}ds.
\EA$$
Thus
\begin{equation}\label {below3}\myint{B_{R}\setminus B_{2|x|}}{}
|x-z|^{2-N}|z|^{-q\ga_{_{I}}}dz\leq \left\{\BA{ll}c_{8}&\text{if }1-q\ga_{_{I}}>-1\\
c_{8}\abs{\ln |x|}&\text{if }1-q\ga_{_{I}}=-1\\
c_{8} |x|^{2-q\ga_{_{I}}}
&\text{if }1-q\ga_{_{I}}<-1.
\EA\right.\end{equation}
Combining (\ref{below2}) and (\ref{below3}) yields to (\ref{below1}).\smallskip

\noindent{\it Step 2.} There holds

\begin{equation}\label {below4}
u_{\infty\gd_{0}}(x)\geq\left( |x|^{-2/q-1)}-r^{-2/(q-1)}\right)\gw_{_{S^{I}_{r}}}(x/|x|)
\quad\forall x\in C^{I}_{r}\cap B_{r},
\end{equation}
where $\gw_{_{S^{I}_{r}}}$ is the unique positive solution of (\ref{mino2r}). For $\ell>0$, let $u^{I}_{\ell\gd_{0}}$ be the solution of
\begin{equation}\label{below5}\left\{\BA{l}
-\Gd u+u^q=0\quad\text{in }C^I_{r}\\
\phantom{-\Gd +u^q}
u=\ell\gd_{0}\quad\text{on }\prt C^I_{r}.
\EA\right.\end{equation}
By comparing $u^{I}_{\ell\gd_{0}}$ with the Martin kernel in $C^I_{r}$,
\begin{equation}\label{below6}
u^{I}_{\ell\gd_{0}}(x)\leq c_{10}\ell |x|^{-\ga_{_{I}}}\quad\forall x\in C^I_{r}.
\end{equation}
Because
\begin{equation}\label{below7}c_{10}\ell|x|^{-\ga_{_{I}}}\leq c^*m|x|^{-\ga_{_{O}}}\quad \forall
x\quad\text{s.t.}\;\abs x\geq c_{11}\left(\frac{\ell}{m}\right)^{(\ga_{_{I}}-\ga_{_{O}})^{-1}},\end{equation}
it follows
\begin{equation}\label{below8}
u_{m\gd_{0}}(x)\geq u^{I}_{\ell\gd}(x)
\quad \forall
x\quad\text{s.t.}\;c_{11}\left(\frac{\ell}{m}\right)^{(\ga_{_{I}}-\ga_{_{O}})^{-1}}\!\!\!\!\!\leq \abs x\leq c^*m^{-\gg_{r}}.
\end{equation}
Notice that (\ref{condi}) implies
$$\left(\frac{\ell}{m}\right)^{(\ga_{_{I}}-\ga_{0})^{-1}}
=o(m^{-\gg_{r}})\quad\text{as }\;m\to\infty.
$$
Since $u^{I}_{\ell\gd_{0}}(x)\leq |x|^{-2/(q-1)}\gw_{_{S^{I}_{r}}}(x/|x|)$, it follows, by the maximum principle, that
$$u_{m\gd_{0}}(x)\geq u^{I}_{\ell\gd_{0}}(x)-r^{-2/(q-1)}\gw_{_{S^{I}_{r}}}(x/|x|)$$
for every $x\in C^I_{r}\cap B_{r}$ \sth $\abs x\geq
c_{11}\left(\frac{\ell}{m}\right)^{(\ga_{_{I}}-\ga_{_{O}})^{-1}}$.
Letting successively $m\to\infty$ and $\ell\to\infty$ and using
$$\lim_{\ell\to\infty} u^{I}_{\ell\gd_{0}}(x)=|x|^{-2/(q-1)}\gw_{_{S^{I}_{r}}}(x/|x|)\quad\forall x\in C^I_{r},
$$
we obtain (\ref{below4}).\smallskip

\noindent{\it Step 3.} Let  $u\in\CU_{+}(\Gw)$, $u$ vanishing on $\prt\Gw\setminus\{0\}$. Because
$$u(x)\leq C_{N,q}|x|^{-2/(q-1)}
$$
and $ \overline{C^I_{r}\cap B_{r}}\setminus\{0\}\subset\Gw$, it is a classical consequence of Harnack inequality that, for any $x$ and $x'\in  \overline{C^I_{r}\cap B_{r/2}}$ such that  $2^{-1}|x|\leq |x'|\leq 2|x|$, $u$ satisfies
$$c_{12}^{-1}u(x')\leq u(x)\leq c_{12}u(x'),
$$
where $c_{12}>0$ depends on $N$, $q$ and $\min\left\{\dist (z,\prt\Gw)/|z|:z\in  \overline{C^I_{r}\cap B_{r}}\right\}$.
\smallskip

\noindent{\it Step 4.} There exists $c_{13}=c_{13}(q,\Gw)>0$ such that
\begin{equation}\label{below9}
U_{0}(x)\leq c_{13} u_{\infty\gd}(x)\quad\forall x\in\Gw.
\end{equation}
Because of (\ref{maj2}) and the fact that for $r>0$ and any compact subset $K\subset S^{I}_{0,r}$
$$1\leq \myfrac{\gw_{_{S^{O}_{0,r}}}(\gs)}{\gw_{_{S^{I}_{0,r}}}(\gs)}\leq M\quad\forall\gs\in K,
$$
where $M$ depends on $K$, there exists $c_{14}>0$ such that
$$1\leq \myfrac{U_{0}(x)}{u_{\infty\gd_{0}}(x)}\leq c_{14}\quad
\forall x\in B_{r}\;\text{ s.t. }\;x/|x|\in K.$$
Using Step 3, there also holds
\begin{equation}\label{below10}
c_{15}^{-1}\leq \min\left\{\myfrac{U_{0}(x')}{U_{0}(x)},\myfrac{u_{\infty\gd_{0}}(x')}{u_{\infty\gd_{0}}(x)}\right\}\leq
\max\left\{\myfrac{U_{0}(x')}{U_{0}(x)},\myfrac{u_{\infty\gd_{0}}(x')}{u_{\infty\gd_{0}}(x)}\right\}
\leq c_{15}\quad\forall x,x'\in B_{r/2},
\end{equation}
provided $x/|x|$ and $x'/|x'|\in K$ and $2^{-1}|x|\leq |x'|\leq 2|x|$.
For $0<s\leq r/2$,  set $\Gg_{s}=\Gw\cap\prt B_{s}$. There exists $n_{0}\in\BBN_{*}$ and $\gk\in (0,1/4)$, independent of $s$, such that for any $x\in\Gg_{s}$ such that $x/|x|\in K$, there exists at most $n_{0}$ points $a_{j}$ ($j=1,...j_{x}$) such that $a_{j}\in \Gg_{s}$,
$a_{1}\in \prt\Gw$, $\gk s\leq \dist (a_{j},\prt\Gw)\leq s$, $|a_{j}-a_{j+1}|\leq s/2$ for $j=1,...j_{x}$ and $a_{j_{x}}=x$. Using \rprop{BHI} and the remark hereafter,
$$
c^{-1}\myfrac{U_{0}(z)}{U_{0}(a_{1})}\leq \myfrac{u_{\infty\gd_{0}}(z)}{u_{\infty\gd_{0}}(a_{1})}\leq c\myfrac{U_{0}(z)}{U_{0}(a_{1})}\quad\forall z\in \Gg_{s}\cap B_{a_{0}}.
$$
Combining with (\ref{below10}) we derive
$$U_{0}(x)\leq cc_{15}^{n_{0}}u_{\infty\gd_{0}}(x)\quad\forall x\in \Gg_{s}.
$$
Because $cc_{15}^{n_{0}}u_{\infty\gd_{0}}$ is a super-solution of (\ref{eq-q})
(clearly $cc_{15}^{n_{0}}>1$),
$$U_{0}\leq cc_{15}^{n_{0}}u_{\infty\gd_{0}} \txt{in $\Gw\setminus \overline B_{s}$}\forevery s\in (0,r].$$
 Thus (\ref{below9}) follows with $c_{13}=cc_{15}^{n_{0}}$.\smallskip

\smallskip

\noindent{\it Step 5.} End of the proof. It is based upon an idea
introduced in \cite{MV1}. If we assume $U_{0}> u_{\infty\gd_{0}}$,
the convexity of $x\mapsto x^q$ implies that the function
$$v=u_{\infty\gd_{0}}-\myfrac{1}{2c_{13}}(U_{0}- u_{\infty\gd_{0}})
$$
is a super solution \sth
$$au_{\infty\gd_{0}}\leq v<u_{\infty\gd_{0}}$$
where  $a=\frac{1+c_{13}}{2c_{13}}<1$. Since $au_{\infty\gd_{0}}$ is
a subsolution, it follows that there exists a solution $w$ \sth
$$au_{\infty\gd_{0}}<w<v<u_{\infty\gd_{0}}.$$
But this is impossible because, for any  $a\in (0,1)$, the smallest
solution dominating $au_{\infty\gd_{0}}$ is $u_{\infty\gd_{0}}$.\qed


\medskip

The next result extends a theorem of Marcus and V\'eron \cite{MV1}.

\bth{wellposed} Assume that $\Gw$ is a bounded Lipschitz domain such
that $S^{^{O}}_{y}=S^{^{I}}_{y}=S_{y}$ for every $y\in \bdw$.
Further, assume that
$$1<q<q^*_{\prt\Gw}.$$
Then for any outer regular Borel measure $\bar\gn$ on $\prt\Gw$
there exists a unique solution $u$ of (\ref{eq-q}) such that
$tr_{_{\prt\Gw}}(u)=\bar\gn$. \es

\Proof We assume $\bar\gn\sim(\gn,F)$ in the sense of \rdef
{g-trace} where $F$ is a closed subset of $\prt\Gw$ and $\gn$ a
Radon measure on $\CR=\prt\Gw\setminus F$. We denote by $U_{F}$ the
maximal solution of (\ref{eq-q}) defined in \rlemma{maxsol}. Because
$q<q_{\prt\Gw}^*$, for any $y\in F$ there exists $u_{\infty\gd_{y}}$
(and actually $u_{\infty\gd_{y}}=U_{y}$ by \rth{uniq}). Then
$U_{F}\geq u_{\gd y}$ by \rlemma{UF-basics}, thus $\CS (U_{F})=F'=F$
with the notation of \rdef{maxsol}.  By \rth{ex-sub}, any Radon
measure is q-admissible thus for any compact subset $E\subset \CR$
there exist a unique solution $u_{\gn\chi_{E}}$ of (\ref{eq-q}) with
boundary trace $\gn\chi_{E}$. Therefore there exists a solution with
boundary trace $\bar\gn$ and, by \rth{reg+sing}, its uniqueness is
reduced to showing that $U_{F}$ is the unique solution with boundary
trace $(0,F)$. Assume $u_{F}$ is any solution with trace $(0,F)$. By
\rth{sub} and \rth{uniq}, there holds
\begin{equation}\label{WP0}
u_{F}(x)\geq u_{\infty\gd_{y}}(x)=U_{y}(x) \qquad\forall y\in
F,\;\forall x\in\Gw.
\end{equation}

\nind Next we prove:

\medskip
 \nind{\it Assertion.\hskip 2mm There exists $C>0$ depending on $F$, $\Gw$ and $q$ such that}
\begin{equation}\label{WP1}
U_{F}(x)\leq Cu_{F}(x)\qquad\forall x\in\Gw.
\end{equation}
There exists $r_{0}>0$ and a circular cone $C_{0}$ with vertex
$0$ and opening $S_{0}\subset\prt B_{1}$ such that for any $y\in \prt\Gw$ there exists an isometry $\CR_{y}$ of $\BBR^N$ such that
$\CR_{y}(\overline C_{0})\cap B_{r_{0}}(y)\subset \Gw\cup\{y\}$. We shall denote by $C_{1}$ a fixed sub-cone  of $C_{0}$ with vertex $0$ and opening $S_{1}\Subset S_{0}$. In order to simplify the geometry, we shall assume that both $C_{0}$ and $C_{0}$ are radially symmetric cones.
If $x\in \Gw$ is such that $\dist (x,\prt\Gw)\leq r_{0}/2$, either \smallskip

\noindent (i) there exists some $y\in\CS$ and an isometry $\CR_{y}$ such that $\CR_{y}(\overline C_{0})\cap B_{r_{0}}(y)\subset \Gw\cup\{y\}$ and $(x-y)/|x-y|\in S_{1} $,
\smallskip

\noindent (ii) or such a $y$ and $\CR_{y}$ does not exist.
\smallskip

\noindent In the first case, it follows from \rprop{critprop2} and \rth{uniq} that
\begin{equation}\label{WP2}
u_{F}(x)\geq c_{1}|x-y|^{-2/(q-1)}.
\end{equation}
Furthermore, the constant $c_{1}$ depends on $r$, $S$ $q$ and $\Gw$, but not on $u_{F}$. By (\ref{OKeq})
\begin{equation}\label{WP3}
U_{F}(x)\leq c_{2}\left(\dist (x,\prt\Gw)\right)^{-2/(q-1)}.
\end{equation}
Since in case (i), there holds $\dist (x,\prt\Gw)\geq c_{3}|x-y|$ for some $c_{3}>1$ depending on $S_{0}$ and $S_{1}$, it follows that (\ref{WP1}) holds with $c=c_{1}c_{2}^{2/(q-1)}/c_{3}$.\smallskip

 In case (ii), $x$ does not belong to any cone radially symmetric cones with opening $S_{1}$ and vertex at some $y\in \CS$. Therefore, there exists $c_{4}<1$ depending on $C_{1}$ such that
\begin{equation}\label{WP4}
\dist (x,\prt\Gw)\leq c_{4}\dist (x,\CS).
\end{equation}
We denote $r_{x}:=\dist (x,\CS)$. If
\begin{equation}\label{WP4'}
\dist (x,\prt\Gw)\leq \min \{c_{4},10^{-1}\}r_{x},
\end{equation}
there exists $\xi_{x}\in\prt\Gw$ such that $|x-\xi_{x}|\dist (x,\prt\Gw)$. Then $B_{9r_{x}/10}(\xi_{x})\subset B_{r_{x}}(x)$. We can apply \rprop{BHI} in $\Gw\cap B_{9r_{x}/10}(\xi_{x})$. Since $x\in B_{r_{x}/5}(\xi_{x})$, there holds
\begin{equation}\label{WP5}
c_{5}^{-1}\myfrac{u_{F}(z)}{U_{F}(z)}\leq
\myfrac{u_{F}(x)}{U_{F}(x)}\leq c_{5}\myfrac{u_{F}(z)}{U_{F}(z)}\qquad\forall z\in B_{r_{x}/5}(\xi_{x})\cap\Gw.
\end{equation}
We can take in particular $z$ such that $|z-\xi_{x}|=r_{x}/5$ and $\dist (z,\prt\Gw)=\max\{\dist (t,\prt\Gw):t\in  B_{r_{x}/5}(\xi_{x})\cap\Gw\}$. Since the distance from
$z$ to $\CS$ is comparable to $\dist (z,\prt\Gw)$, there exist $n_{0}\in\BBN_{*}$ depending on the geometry of $\Gw$ and $n_{0}$ points
$\{a_{j}\}$ with the properties that  $\dist (a_{j},\prt\Gw)\geq \dist (z,\prt\Gw)$,
$B_{r_{x}/10}(a_{j})\cap B_{r_{x}/10}(a_{j+1})\neq\emptyset$ for $j=1,...,n_{0}-1$, $a_{1}=z$ and $a_{n_{0}}$ have the property (i) above, that is  there exists some $y\in\CS$ and an isometry $\CR_{y}$ such that $\CR_{y}(\overline C_{0})\cap B_{r_{0}}(y)\subset \Gw\cup\{y\}$ and $(a_{n_{0}}-y)/|a_{n_{0}}-y|\in S_{1}$. By classical Harnack inequality (see \rth{uniq} Step 3), there holds
$$u_{F}(a_{j})\geq c_{6}u_{F}(a_{j+1})\quad\text{and }\;\;U_{F}(a_{j})\geq c^{-1}_{6}U_{F}(a_{j+1})
$$
for some $c_{6}>1$ depending on $N$, $q$ and $\Gw$ via the cone $C_{0}$. Therefore
\begin{equation}\label{WP6}
U_{F}(x)\leq c_{5}c_{6}^{2n_{0}}\myfrac{u_{F}(a_{n_{0}})}{U_{F}(a_{n_{0}})}u_{F}(x)\leq c_{7}u_{F}(x),
\end{equation}
which implies (\ref{WP1}) from case (i) applied to $a_{n_{0}}$. \smallskip

Finally, if (\ref{WP4}) holds, but also
\begin{equation}\label{WP4''}
\dist (x,\prt\Gw)\geq \min \{c_{4},10^{-1}\}r_{x},
\end{equation}
this means that $\dist (x,\prt\Gw)$ is comparable to $r_{x}$. Then we can perform the same construction as in the case (\ref{WP4'}) holds, except that we consider balls
$B_{\dist (x,\prt\Gw)/4}(a_{j)}$ in order to connect $x$ to a point $a_{n_{0}}$ satisfying (i). The number $n_{0}$ is always independent of $u_{F}$. Thus we derive again
estimate (\ref{WP1}) provided $\dist (x,\prt\Gw)\leq r_{0}/2$. In order to prove that this holds in whole $\Gw$, we consider some $0<r_{1}\leq r_{0}/2$ such that
$\Gw'_{r_{1}}:=\{x\in\Gw:\dist(x,\Gw)>r_{1}\}$ is connected. The function $v$ solution of
\begin{equation}\label{WP6'}\left\{\BA {l}
-\Gd v+v^q=0\quad\text{in }\Gw'_{r_{1}}\\
\phantom{-\Gd +v^q}
v=c_{1}u_{F}\quad\text{in }\prt\Gw'_{r_{1}}
\EA\right.\end{equation}
is larger that $U_{F}$ in $\Gw'_{r_{1}}$. Since $c_{1}u_{F}$ is a super solution, $v\leq c_{1}u_{F}$  in $\Gw'_{r_{1}}$. This implies that (\ref{WP1}) holds in $\Gw$.\smallskip

Inequality \eqref{WP1} implies uniqueness by the same argument as in
the proof of \rth{uniq}, Step 5. \qeda


\mysection{The Martin kernel and critical values for a cone.}

 \subsection{The geometric framework}
\setcounter{equation}{0}

 An N-dim polyhedra $P$ is the bounded domain bordered by a finite number of hyperplanes.
  Thus characteristic  elements of the boundary of  $P$ are the faces (subsets of an hyperplane),
  the vertex (intersection of $N$ hyperplanes) and a wide variety of N-k dimensional  edges,
  where k ranges from 2 to N. An N-k dimensional edge, i.e. an intersection of k hyperplanes,
  will be described by the characteristic angles of these hyperplanes.
  \medskip

 We recall that the spherical coordinates in $\BBR^N=\{x=(x_{1},...x_{N})\}$ are expressed by
 \begin{equation}\label{N-dim1}\left\{\BA {l}
 x_{1}=r\sin\gth_{N-1}\sin\gth_{N-2}...\sin\gth_{2}\sin\gth_{1}\\
 x_{2}=r\sin\gth_{N-1}\sin\gth_{N-2}...\sin\gth_{2}\cos\gth_{1}\\
 x_{3}=r\sin\gth_{N-1}\sin\gth_{N-2}...\cos\gth_{2}\\
 .\\
 .\\
 .\\
 x_{N-1}=r\sin\gth_{N-1}\cos\gth_{N-2},\\
 x_{N}=r\cos\gth_{N-1}
\EA\right. \end{equation}
where
$\gth_{1}\in [0,2\gp]$ and $\gth_{\ell}\in[0,\gp]$ for $\ell=2,3,...,N-1$ (the $\gth_{j}$ are the Euler
angles). Thus the ''angular'' component $\gs\in S^{N-1}$ of the spherical coordinates
$(r,\gs)$ of $x\in\BBR^{N}$ is denoted by $\gs=(\gth_{1},...,\gth_{N-1})$. \medskip

We consider an unbounded {\em non-degenerate k-dihedron}\/ defined as follows. Let
$k\in[2,N]\cap\BBN$ and  let $A$ be given by
 $$A=\begin{cases}(0,\ga_{1})\ti \prod_{j=2}^{k-1}(\ga_{j},\ga'_{j})
 &\text{if $k>2$}\\
(0,\ga_{1})
&\text{if $k=2$}\end{cases}
 $$
 where
 $$0<\ga_{1}< 2\gp,\q 0<\ga_{j}<\ga'_{j}<\gp \q j=2,...,k-1.$$
 We denote by $S_A$ the spherical domain
\begin{equation}\label{dom}
S_A=\{x\in \BBR^N:\abs x=1,\,\gs\in
A\ti\prod_{j=k}^{N-1}[0,\gp]\}\sbs S^{N-1}\}
\end{equation}
and by $D_A$ the corresponding  k-dihedron,
$$D_{A}=\{x=(r,\gs):r>0,\gs\in S_A\}.$$
The {\em edge}\/ of $D_{A}$ is the (N-k)-dimensional space
 \begin{equation}\label{edge}
d_{A}=\{x:x_{1}=x_{2}=...=x_{k}=0\}.
\end{equation}
 \subsection{Separable harmonic functions and the Martin kernel in a k-dihedron.}
In the system of spherical coordinates, the Laplacian takes the form
$$\Gd u=\prt_{rr}u+\myfrac{N-1}{r}\prt_{r}u+
\myfrac{1}{r^2}\Gd_{_{S^{N-1}}} u
$$
where the Laplace-Beltrami operator $\Gd_{_{S^{N-1}}}$ is expressed by induction by
 \begin{equation}\label{lapla}\BAL
\Gd_{_{S^{N-1}}} u=&\myfrac{1}{(\sin\gth_{N-1})^{{N-2}}}\myfrac{\prt}{\prt\gth_{N-1}}\left((\sin\gth_{N-1})^{{N-2}}\myfrac{\prt u}{\prt\gth_{N-1}}\right)\\
&+\myfrac{1}{(\sin\gth_{N-1})^{{2}}}\Gd_{_{S^{N-2}}} u.
\EAL \end{equation}
 and
  \begin{equation}\label{lapla1}
\Gd_{_{S^{1}}} u=\prt_{\gth_{1}\gth_{1}}u
 \end{equation}
If we compute the positive harmonic functions in the k-dihedron $D_{A}$ of the form
$$v(x)=v(r,\gs)=r^\gk\gw(\gs) \q\text{in }D_A,\q v=0\q\text{in }\prt D_{A}\setminus\{0\}.
$$
we find that
$\gk$ satisfies the algebraic equation
  \begin{equation}\label{kappa1}
\gk^2+(N-2)\gk-\gl_{A}=0
 \end{equation}
 where $\gl_{A}$ is the first eigenvalue of $-\Gd_{_{S^{N-1}}}$ in $W^{1,2}_{0}(S_A)$
 and $\gw$ is the corresponding normalized eigenfunction:
   \begin{equation}\label{kappa2-0}\left\{\BA {l}
\Gd_{_{S^{N-1}}}\gw+\gl_{A}\gw=0\quad\text{in }S_{A}\\
\phantom{\Gd_{_{S^{N-1}}}\gw+\gl_{A}}\gw=0\quad\text{on }\prt S_A.
\EA\right. \end{equation}
Thus
   \begin{equation}\label{kappa2}\BAL
   \gk_{+}&=\myfrac{1}{2}\left(2-N+\sqrt{(N-2)^2+4\gl_{A}}\right)\\
   \gk_{-}&=\myfrac{1}{2}\left(2-N-\sqrt{(N-2)^2+4\gl_{A}}\right).
\EAL\end{equation}
Since
 \begin{equation}\label{N-dim2}
 S^{N-1}=\left\{\gs\in \BBR^{N-1}\ti\BBR:\gs=(\gs_2\sin\gth_{N-1},\cos\gth_{N-1}),\;\gs_2\in S^{N-2}\right\},
 \end{equation}
we look for $\gw:=\gw^{\{1\}}$ of the form
$$\gw^{\{1\}}(\gs)=(\sin\gth_{N-1})^{\gk_{+}}\gw^{\{2\}}(\gs_2),
\quad\gth_{N-1}\in (0,\gp),\q\gs_2\in S^{N-2}.
$$
Here $S^{N-2}=S^{N-1}\cap\{x_N=0\}$ and we denote
$$S^{\{N-2\}}_A=S_A\cap\{x_{N}=0\},\q D^{\{N-2\}}_{A}:=D_A\cap\{x_{N}=0\}
\subset \BBR^{N-1}.$$
Then (\ref{kappa2})
jointly with relation (\ref{lapla}) implies
\begin{equation}\label{EV2}\left\{\BA {l}
\Gd_{_{S^{N-2}}}\gw^{\{2\}}+(\gl_{A}-\gk_{+})\gw^{\{2\}}=0\quad\text{on }S^{\{N-2\}}_A\\[2mm]
\phantom{\Gd_{_{S^{N-2}}}\gw^{\{2\}}+(\gl_{A}-\gk_{+})}
\gw^{\{2\}}=0\quad\text{on }\prt S^{\{N-2\}}_A.
\EA\right.\end{equation}
Since  we are interested in  $\gw^{\{2\}}$  positive, $\gl_{A}^{\{2\}}:=\gl_{A}-\gk_{+}$ must be
 the first eigenvalue of $-\Gd_{_{S^{N-2}}}$ in $W^{1,2}_{0}(S^{\{N-2\}}_A)$.

Next we look for positive harmonic functions  $\tl u$ in $D^{\{N-2\}}_{A}$ \sth
$$\tilde u(x_1,\ldots,x_{N-1})= r^{\gk'}\gw(\gs_2),\q  \tilde u=0 \text{ on }\prt D^{\{N-2\}}_{A}$$
 The algebraic equation which gives the exponents is
$$(\gk')^2+(N-3)\gk'-\gl_{A}^{\{2\}}=0.
$$
Denote by $\gk'_{+}$  the positive root of this equation.
By the definition of $\gl_{A}^{\{2\}}$,
$$\gk_+^2+(N-3)\gk_+-\gl_{A}^{\{2\}}=\gk_+^2+(N-2)\gk_+-\gl_{A}=0. $$
Therefore $\gk'_{+}=\gk_{+}$.
Accordingly, if $k\geq 3$, we set
$$\gw^{\{2\}}(\gs_2)=(\sin\gth_{N-2})^{\gk_{+}}\gw^{\{3\}}(\gs_{3}),
$$
an find that $\gw^{\{3\}}$ satisfies
\begin{equation}\label{EV3}\left\{\BA {l}
\Gd_{_{S^{N-3}}}\gw^{\{3\}}+(\gl_{A}-2\gk_{+})\gw^{\{3\}}=0\quad\text{in }S^{\{N-3\}}_{A}\\[2mm]
\phantom{\Gd_{_{S^{N-3}}}\gw^{\{3\}}+(\gl_{A}-2\gk_{+})}
\gw^{\{3\}}=0\quad\text{on }\prt S^{\{N-3\}}_{A},
\EA\right.\end{equation}
where
$$S^{\{N-3\}}_{A}=S_A\cap\{x_{N}=x_{N-1}=0\}.$$
Performing this reduction process (N-k) times, we obtain the following results.\smallskip

\noindent(i) If $k>2$ then
\begin{equation}\label{EV4}
\gw(\gs)=(\sin\gth_{N-1}\sin\gth_{N-2}...\sin\gth_{k})^{\gk_{+}}\gw^{\{N-k+1\}}(\gs_{N-k+1})
\end{equation}
where $$\gs_{N-k+1}\in S^{k-1}=S^{N-1}\cap\{x_N=,x_{N-1}=\cdots=x_{k+1}=0\},$$
and $\gw':=\gw^{\{N-k+1\}}$ satisfies
\begin{equation}\label{EV5-0}\left\{\BA {l}
\Gd_{_{S^{k-1}}}\gw'+(\gl_{A}-(N-k)\gk_{+})\gw'=0,\quad\text{in } S^{\{k-1\}}_{A}\\[2mm]
\phantom{\Gd_{_{S^{k-1}}}\gw'+(\gl_{A}-(N-k)\gk_{+})}\gw'=0,\quad\text{on
}\prt S^{\{k-1\}}_{A}, \EA\right.\end{equation}

$$S^{\{k-1\}}_{A}=S_A\cap\{x_{N}=x_{N-1}=...=x_{k+1}=0\}\approx A
$$
and $\gl_{A}-(N-k)\gk_{+}$ is the first eigenvalue of the problem.
In such a case, it is usually impossible to determine more explicitly $\gw^{\{N-k+1\}}$ and
$\gl_{A}-(N-k)\gk_{+}$, except for some very specific values of  $\ga_{j}$ and $\ga'_{j}$,
associated to consecutive zeros of generalized Legendre functions.\smallskip

\noindent(ii) If $k=2$ then
\begin{equation}\label{EV5-2}
\gw(\gs)=(\sin\gth_{N-1}\sin\gth_{N-2}...\sin\gth_{2})^{\gk_{+}}
\gw^{\{N-1\}}(\gth_{1})
\end{equation}
where $\gs_{N-1}\in S^{1}\approx \gth_{1}\in (0,2\gp)$, and $\gw^{\{N-1\}}$ satisfies
\begin{equation}\label{EV5-3}\left\{\BA {l}
\Gd_{_{S^{1}}}\gw^{\{N-1\}}+(\gl_{A}-(N-2)\gk_{+})\gw^{\{N-1\}}=0\quad\text{on }S^{\{1\}}_A\\[2mm]
\phantom{\Gd_{_{S^{1}}}\gw^{\{N-1\}}+(\gl_{A}-(N-2)\gk_{+})}\gw^{\{N-1\}}=0\quad\text{on }\prt S^{\{1\}}_A,
\EA\right.\end{equation}
with $\prt S^{\{1\}}_A\approx (0,\ga)$. In this case
\begin{equation}\label{EV6}
\gk_{+}=\myfrac{\gp}{\ga},\q \gw^{\{N-1\}}(\gth_{1})=\sin(\gp\gth_{1}/\ga),
\end{equation}
and
\begin{equation}\label{EV7}
\gl_{A}-(N-2)\gk_{+}=\myfrac{\gp^2}{\ga^2}\Longrightarrow
\gl_{A}=\myfrac{\gp^2}{\ga^2}+(N-2)\myfrac{\gp}{\ga}.
\end{equation}
Observe that $\rec{2}\leq\gk_+$ with equality holding only in the degenerate case $\ga=2\gp$ (which we exclude).
\medskip

In either case,  we find a positive harmonic function $v_{A}$ in
$D_{A}$, vanishing on $\prt D_{A}$, of the form
$$v_{A}(x)=\abs x^{\gk_{+}}\gw(x/\abs x)
$$
with $\gw$ as in \eqref{EV4} (for $k>2$) or \eqref{EV6} (for k=2).

Similarly we find a  positive harmonic function in $D_A$ vanishing  on $\prt D_A\sms \{0\}$,
singular at the origin, of the form
$$K'_{A}(x)=\abs x^{\gk_{-}}\gw(x/\abs x),\q \gk_{-}=2-N-\gk_{+}.$$
Because of the uniqueness of the kernel function (see {\bf A.2}) $K'_{A}(x)$ is, up to a
multiplicative constant $c_{A}$, the
Martin kernel of the Laplacian in $D_{A}$, with singularity at $0$.
  The Martin kernel, with singularity at a
point $z\in d_A$, is given by
\begin{equation}\label{k3}
K_{A}(x,z)=c_{_{A}}\myfrac{(\sin\gth_{N-1}\sin\gth_{N-2}...\sin\gth_{k})^{\gk_{+}}\gw^{\{N-k+1\}}(\gs_{N-k+1})}{|x-z|^{N-2+\gk_{+}}}
\end{equation}
for every $x\in D_A$. From (\ref{N-dim1})
$$\sin\gth_{N-1}\sin\gth_{N-2}...\sin\gth_{k}=|x-z|^{-1}{\sqrt{x_{1}^2+x_{2}^2+...+x_{k}^2}}.
$$
Therefore, if we write $x\in\BBR^N$ in the form $x=(x',x'')$,
$x'=(x_{1},...,x_{k})$, $x''=(x_{k+1},\cdots,x_N)$,
we obtain the formula,
\begin{equation}\label{k4}\BA{l}
K_{A}(x,z)=c_{_{A}}\myfrac{|x'|^{\gk_{+}}\gw^{\{N-k+1\}}(\gs_{N-k+1})}{|x-z|^{(N-2+2\gk_{+})}}\\[4mm]
\phantom{K_{A}(x,z)}
=c_{_{A}}\myfrac{|x'|^{\gk_{+}}\gw^{\{N-k+1\}}(\gs_{N-k+1})}{(|x'|^2+|x''-z|^2)^{(N-2+2\gk_{+})/2}}.
\EA\end{equation}
Therefore, the Poisson potential of a measure
$\gm\in\GTM (d_{A})$ is expressed by
\begin{equation}\label{k5}\BA{l}
\BBK[\gm](x)
\quad=c_{_{A}}|x'|^{\gk_{+}}\gw^{\{N-k+1\}}(\gs_{N-k+1})\\[2mm]
\phantom{----------}\ti\myint{\BBR^{N-k}}{}
\myfrac{ d\gm(z)}{(|x'|^2+|x''-z|^2)^{(N-2+2\gk_{+})/2}}.
\EA\end{equation}

\subsection{The admissibility condition}
Consider the boundary value problem
\begin{equation}\label{BVP}\left\{\BA {l}
-\Gd u+\abs u^{q-1}u=0\quad\text {in }D_A\\
\phantom{-\Gd u+\abs u^{q-1}}
u=\gm\in\GTM(\prt D_A).
\EA\right.\end{equation}
Let $r'=|x'|$, $r''=|x''|$,
\begin{equation}\label{GgR}
\Gg_{R}=\{x=(x',x''):r'\leq R,r''\leq R\}
\end{equation}
and let $\gr_R$ denote the first (positive) eigenfunction in
$D_{A,R}:=D_A\cap \Gg_R$.

By \rdef{admissible}, the admissibility condition for a measure
$\gm\in\GTM(d_{A}\cap\Gg_R)$ relative to \eqref{BVP} in $D_{A,R}$ is
\begin{equation}\label{admiss-R}
\int_{D_{A,R}}\BBK^R[|\gm|](x)^q\gr\ssub{R}(x)dx<\infty.
\end{equation}
where $\BBK^R$ is the Martin kernel of $-\Gd$ in $D_{A,R}$. Near $d_A$
this kernel behaves like the Martin kernel of the dihedron $D_{A}$.
Furthermore, the first eigenfunction $\gr\ssub{R}$ of $-\Gd$ in
$W^{1,2}_{0}(D_{A,R})$ behaves like the regular harmonic
function $v_{A}$. Therefore
\begin{equation}\label{grR}
 \gr\ssub{R}(x)\approx (r')^{\gk_{+}}\gw^{\{N-k+1\}}(\gs_{N-k+1})
\end{equation}
and the admissibility condition for a measure $\gm\in\frak M(d_{A})$
is
\begin{equation}\label{k6}
\int_{\Gg_{R}\cap D_{A}}\BBK[|\gm|](x)|^q\gr(x)dx<\infty \forevery
R>0,
\end{equation}
with $\Gg_R$ as in \eqref{GgR}. By \eqref{k5},
\begin{equation}\label{k7}
\BBK[|\gm|](x)\leq \
c_{_{A}}(r')^{\gk_{+}}\int_{\BBR^{N-k}}j(x',x''-z) d|\gm|(z)
\end{equation}
where
\begin{equation}\label{k7bis}
  j(x)=|x|^{-N+2-2\gk_+} \forevery x\in
  \BBR^N.
\end{equation}
 Therefore, using \eqref{grR}, condition \eqref{k6} becomes
\begin{equation}\label{k8}
\int_0^R\int_{|x''|<R} \Big(\int_{\BBR^{N-k}}
j(x',x''-z)d|\gm|(z)\Big)^q (r')^{(q+1)\gk_{+}+k-1}dx''dr'<\infty
\end{equation}
for every
$R>0$.
\subsection{The critical values.}
Relative to the equation
\begin{equation}\label{eqq}
   -\Gd u+|u|^{q-1}u=0
\end{equation}
there exist two thresholds of criticality associated with the edge
$d_A$.

 The first is the
value $q^*_c$ \sth, for $q^*_c\leq q$
 the whole edge $d_{A}$ is removable relative  to this  equation,
 but for $1<q<q^*_c$ there exist non-trivial solutions in $D_A$
 which vanish on $\prt D_A\sms d_A$.
The second $q_c<q^*_c$  corresponds to the removability of  points
on $d_{A}$. For $q\geq q_c$ points on $d_A$ are removable while for
$1<q<q_c$ there exist solutions with isolated point singularities on
$d_A$. In the next two propositions we determine these critical
values.
\bprop{admp0} Assume $q>1$, $1\leq k<N$. Then the condition
\begin{equation}\label {qk}
q<
q^*_c:=1+\myfrac{2-k+\sqrt{(k-2)^2+4\gl_{A}-4(N-k)\gk_{+}}}{\gl_{A}-(N-k)\gk_{+}}
\end{equation}
is necessary and sufficient for the existence of a non-trivial
solution $u$ of (\ref{eqq}) in $D_{A}$ which vanishes on $\prt
D_{A}\sms d_A$. Furthermore, when this condition holds, there exist
non-trivial positive bounded measures $\mu$ on $d_A$ \sth
$\BBK[\gm]\in L^q_{\gr}(\Gg_{R}\cap D_{A})$. \es

\Remark The statement remains true in the case $k=N$, which is the
case of the cone. In this case $q_c=q^*_c=1+(2/\ga_S)$ in the
notation of Section 5. However, in the present notation,
$\ga_S=-\gk_-$ and a straightforward computation yields:
\begin{equation}\label {qkN}
q_c=\frac{N+2+\sqrt{(N-2)^2+4\gl_A}}{N-2+\sqrt{(N-2)^2+4\gl_A}}.
\end{equation}

 \Proof Recall that $\gl_{A}-(N-k)\gk_{+}$ is the first
eigenvalue in $S_A^{\{k-1\}}$ (see \eqref{EV5-0} and the remarks
following it). Let $\gk'_+,\gk'_-$ be the two roots of the equation
$$X^2+(k-2)X-(\gl_{A}-(N-k)\gk_{+})=0,$$
i.e.
$$\gk'_\pm=\rec{2}\big(2-k\pm\sqrt{(k-2)^2+4(\gl_{A}-(N-k)\gk_{+}}\big).$$
Then, by \rth{admiss} and \rth{cone1}, if $1<q<1-(2/{\gk'_-})$ (note
that because of a change in notation the entity denoted by $\ga_S$
in subsection 5.1 is the same as  $-\gk'_-$ in the present section)
there exists a unique solution of (\ref{eqq}) in the cone
$C_{S_A^{k-1}}$ i.e. the cone with opening $S_A^{k-1}\sbs
S^{k-1}\sbs \BBR^k$
 with trace $a\gd_0$ (where $\gd_0$ denotes the Dirac measure at the vertex of the cone and
 $a>0$). By \rth{cone1} this solution satisfies
 \begin{equation}\label{lim-k'}
u_{a}(x)=a\abs x^{-\ga}\gf(x/\abs x)(1+o(1))\quad\text{as }x\to 0,
\end{equation}
where $\gf$ is the first positive eigenfunction of $-\Gd'$ in
$W^{1,2}_{0}(S_A^{k-1})$ normalized so that $u_1$ possesses trace
$\gd_0$.

The function $u$ given by
$$\tl u_a(x',x'')= u_a(x')\quad\forall (x',x'')\in D_{A}=C_{S_A^{k-1}}\ti \BBR^{N-k},
$$
is a nonzero solution of (\ref{eqq}) in $D_{A}$ which vanishes on
$\prt D_{A}\sms d_A$ and has bounded trace on $d_A$.

A simple calculation shows that $1-(2/{\gk_-})$ equals $q_c^*$ as
given in \eqref{qk}.

\smallskip

Next, assume that $q\geq q^*_c$ and let $u$ be a solution of
(\ref{eqq}) in $D_{A}$ which vanishes on $\prt D_{A}\sms d_A$.

Given $\ge>0$ let $v_\ge$ be the solution of \eqref{eqq} in
$D_A^{\{N-k-1\}}\setminus \{x'\in \BBR^k:|x'|\leq \ge\}$ \sth

$$v_\ge(x')=\begin{cases} 0, &\txt{if $x'\in \prt D_A^{\{N-k-1\}}$}\\
\infty, &\txt{if $|x'|=\ge$.}
\end{cases}$$
Given $R>0$ let $w\ssub{R}$ be the maximal solution in $\{x''\in
\BBR^{N-k}:|x''|<R\}$.

Then the function $u^*$ given by
$$u^*(x',x'')=v_\ge(x')+w\ssub{R}(x'')$$
is a supersolution of \eqref{eqq} in
$D_A\sms\{(x',x''):|x'|>\ge,\;|x''|<R\}$ and it dominates $u$ in
this domain. But $w\ssub{R}(x'')\to 0$ as $R\tin$ and, by \cite{FV},
$v_\ge(x')\to 0$ as $\ge\to 0$. Therefore $u_+=0$ and, by the same
token, $u_-=0$. \qed
\bprop{admpk} Let $A$ be defined as before. Then
\begin{equation}\label{admpk}
  \BBK[\gm]\in L^q_{\gr}(\Gg_{R}\cap D_{A})\forevery \mu\in \GTM(d_A),\forevery R>0
\end{equation}
 if and only if
\begin{equation}\label{q-critk}
1<q<q_{c}:=\myfrac{\gk_{+}+N}{\gk_{+}+N-2}.
\end{equation}
This statement is equivalent to the following:

Condition \eqref{q-critk} is necessary and sufficient in order that
the Dirac measure $\mu=\gd_P$, supported at a point $P\in d_A$, satisfy \eqref{admpk}.
\es \Proof It is sufficient to prove the result relative to the family of measures $\mu$ such that
 $\mu$ is positive,  has compact
support and  $\mu(d_A)=1$. Let $R>1$ be sufficiently large so
that the support of $\mu$ is contained in $\Gg_{R/2}$. The measure
$\mu$ can be approximated (in the sense of weak convergence of
measures) by a \seq $\set{\mu_n}$ of convex combinations of Dirac
measures supported in $d_A\cap \Gg_{R/2}$. For such a \seq
$\BBK[\mu_n]\to \BBK[\mu]$ pointwise and  $\set{\BBK[\mu_n]}$ is
uniformly bounded in $D_A\sms \Gg_{3R/4}$. Therefore it is
sufficient to prove the result when $\gm=\gd_{0}$.  In this case the
admissibility condition \eqref{k8}) is
$$\int_0^R\int_{|x''|<R} j(x)^q
(r')^{(q+1)\gk_{+}+k-1}dx''dr'<\infty,$$ i.e.,
$$\int_0^R\int_0^R
|x|^{q(2-N-2\gk_+)}(r')^{(q+1)\gk_{+}+k-1}(r'')^{N-k-1}dr''dr'<\infty.$$
Substituting $\tau:=r''/r'$ the condition becomes
$$\int_0^R\int_0^{R/r'}
(1+\tau^2)^{\myfrac{q}{2}(2-N-2\gk_+)}(r')^{q(2-N-\gk_+)+\gk_++N-1}\tau^{N-k-1}d\tau\,dr'<\infty.$$
This holds if and only if $q<(\gk_++N)/(\gk_++N-2)$. \qed
\medskip

\noindent\Remark It is interesting to notice that $k$ does not
appear explicitly in \eqref{q-critk}. Furthermore, we  observe that
\begin{equation}\label{q-critk2}
\myfrac{2}{q_{c}-1}\left(\myfrac{2q_{c}}{q_{c}-1}-N\right)=\gl_{A}\Longleftrightarrow \gk_{+}(\gk_{+}+N-2)=\gl_{A},
\end{equation}
which follows from \eqref{kappa2}. This implies that there does not exist a
nontrivial  solution of the nonlinear eigenvalue problem
\begin{equation}\label{nlne1}\BAL
 -\Gd_{_{S^N-1}}\psi-\myfrac{2}{q-1}\left(\myfrac{2q}{q-1}-N\right)\psi+
 |\psi|^{q-1}\psi&=0\q\text{in }S_{_{D_{A}}}\\[2mm]
\psi&=0\q\text{in }\prt S_{_{D_{A}}} \EAL
\end{equation}
which, in turn,  implies that there does not exists a nontrivial solution of
(\ref{eqq}) of the form $u(x)=u(r,\gs)=|x|^{-2/(q-1)}\psi(\gs)$, and
also no solution of this equation in $D_{A}$ which vanishes on $\prt
D_{A}\setminus\{0\}$. This is the classical ansatz for the
removability of isolated singularities in $d_{A}$. 

\section{The harmonic lifting of a Besov space $B^{-s,p}(d_A)$.}
Denote by $W^{\gs,p}(\BBR^\ell)$ ($\gs>0$, $1\leq p\leq\infty$)
the Sobolev spaces over $\BBR^\ell$. In order to use interpolation,
it is useful to introduce the Besov space $B^{\gs,p}(\BBR^{\ell})$
($\gs>0$). If $\gs$ is not an integer then
 \begin{equation}\label{B01}
B^{\gs,p}(\BBR^{\ell})=W^{\gs,p}(\BBR^{\ell}).
\end {equation}
If $\gs$ is an integer the space is defined as follows.
Put
$$\Gd_{x,y}f=f(x+y)+f(x-y)-2f(x).$$
Then
 \begin{equation}\label{B1}
B^{1,p}(\BBR^{\ell})=\left\{f\in L^p(\BBR^{\ell}):
\myfrac{\Gd_{x,y}f}{|y|^{1+\ell/p}}\in L^p(\BBR^{\ell}\ti
\BBR^{\ell})\right\},
\end {equation}
with norm
 \begin{equation}\label{B2}
\norm f_{B^{1,p}}=\norm f_{L^{p}}+ \left(\dint_{\BBR^{\ell}\ti
\BBR^{\ell}}
\myfrac{|\Gd_{x,y}f|^p}{|y|^{\ell+p}}dx\,dy\right)^{1/p},
\end {equation}
(with standard modification if $p=\infty$) and
 \begin{equation}\label{B2'}\BAL
B^{m,p}(\BBR^{\ell})=\Big\{&f\in
W^{m-1,p}(\BBR^{\ell}):\\
&D_x^{\ga}f\in B^{1,p}(\BBR^{\ell})\;\forall \ga\in \BBN^{\ell},\;|\ga|=m-1\Big\}
\EAL
\end{equation}
with norm
 \begin{equation}\label{B2''}\BAL
 \norm f_{B^{m,p}}=\norm f_{W^{m-1,p}}+
\left(\sum_{|\ga|=m-1}\dint_{\BBR^{\ell}\ti \BBR^{\ell}}
\myfrac{|D_x^{\ga}\Gd_{x,y}f|^p}{|y|^{\ell+p}}dx\,dy\right)^{1/p}.
\EAL \end {equation}

We recall that the following inclusions hold (\cite[p 155]{St})
  \begin{equation}\label{B3}\BA {l}
W^{m,p}(\BBR^{\ell})\subset B^{m,p}(\BBR^{\ell})
\quad\text{if }\,p\geq 2\\[2mm]
B^{m,p}(\BBR^{\ell})\subset W^{m,p}(\BBR^{\ell})
\quad\text{if }\,1\leq p\leq 2.
\EA\end {equation}
When $1<p<\infty$, the dual spaces of $W^{s,p}$ and $B^{m,p}$ are respectively denoted by $W^{-s,p'}$
and $B^{-m,p'}$.


The following is the main result of this section.
\bth{main1}Suppose that $q_c<q<q^*_c$ and let $A$ be defined as in
subsection {\bf 6.1}. Then there exist positive constants $c_1,c_2$, depending on $q,N,k,\gk_+$,
 such that for any $R>1$ and any $\gm\in\GTM_+(d_{A})$ with support in $B_{R/2}$:
\begin{equation}\label{L7}\BAL
&c_1\norm{\gm}^q_{B^{-s,q}(\BBR^{N-k})}\\
&\leq \int_{D_{A,R}}\BBK[|\gm|]^q(x)\gr(x)dx \leq
c_2(1+R)^{\gb}\norm{\gm}^q_{B^{-s,q}(\BBR^{N-k})},
\EAL\end{equation}
where $s=2-\frac{\gk_{+}+k}{q'}$, $\gb=(q+1)\gk_++k-1$ and $D_{A,R}=D_{A}\cap\Gg_R$. If $q=q_c$
the estimate
remains valid for measures $\mu$ \sth the diameter of $\supp \mu$ is sufficiently small
(depending on the parameters mentioned before).
 \es
\Remark When $q\geq 2$ the norms in the Besov space may be replaced
by the norms in the corresponding Sobolev spaces.\2
\indent Recall the admissibility condition for a measure
$\gm\in\GTM_{+}(d_{A})$:
$$\int_{D_{A,R}}\BBK[\gm]^q(x)\gr(x)dx <\infty \forevery R>0$$
and the equivalence (see \eqref{k6}--\eqref{k8})
\begin{align}\label{k8'}
&\int_{D_{A,R}}\BBK[\gm]^q(x)\gr(x)dx\approx J^{A,R}(\mu):=\\
&\int_0^R\int_{B''_R} \Big(\int_{\BBR^{N-k}}
\frac{d\gm(z)}{(\tau^2+|x''-z|^2)|)^{(N-2+2\gk_+)/2}}\Big)^q
\tau^{(q+1)\gk_{+}+k-1}dx''d\tau, \notag
\end{align}
where $x=(x',x'')\in
\BBR^k\ti\BBR^{N-k}$,  $\tau=|x'|$ and $B''_R=\set{x''\in\BBR^{N-k}: |x''|<R}$. We denote,
\begin{equation}\label{nu}
   \nu=N-2+2\gk_+.
\end{equation}

If $2\gk_+$ is an integer, it is natural to relate \eqref{k8'} to the Poison potential of $\mu$
in $\BBR_+^n=\BBR_+\ti\BBR_{n-1}$ where $n=N-2+2\gk_+$. We clarify this statement below.

Assuming that
$2\leq n+k-N$, denote
$$y=(y_1,\wtl y,y'')\in \BBR^n,\q \wtl y=(y_2,\cdots,y_{n+k-N}),\q
y''=(y_{n+k-N+1},\cdots,y_n).$$

 The Poisson kernel in $\BBR^n_{+}=\BBR_+\ti\BBR_{n-1}$ is given by
\begin{equation}\label{P1}
P_{n}(y)=\gg_{n}y_{1}|y|^{-n}\quad y_{1}>0,
\end{equation}
for some $\gg_{n}>0$,  and the Poisson potential of a bounded Borel measure $\gm$ with support in
$${\bf d}:=\{y=(0,y'')\in \BBR^n:\,y''\in \BBR_{N-k}\}$$
is
\begin{equation}\label{P2}\BBK_{n}[\gm](y)=
\gg_{n}y_{1}\int_{\BBR^{N-k}}\myfrac{d\gm(z)}{\left(y_1^2+|\wtl y|^2+|y''-z|^2\right)^{n/2}}\forevery y\in \BBR^n_{+}.
\end{equation}
In particular, for $\wtl y=0$,
\begin{equation}\label{P3}\BBK_{n}[\gm](y_{1},0,y'')=\gg_{n}y_{1}\myint{\BBR^{N-k}}{}\myfrac{d\gm(z)}{\left(y_{1}^2+|y''-z|^2\right)^{n/2}}.
\end{equation}
The  integral in \eqref{P3} is precisely the same as the inner
integral in \eqref{k8'}.

In fact, it will be shown that, if we set
\begin{equation}\label{def-n}
n:=\{\nu\}=\inf\set{m\in\BBN: m\geq\nu},
\end{equation}
this approach also works when $2\gk_+$ is not an integer.
We note that, for $n$ given by \eqref{def-n},
\begin{equation}\label{n-N+k}
n-N+k\geq 2,
\end{equation}
with equality only if $k=3$ and $\gk_+\leq 1/2$ or $k=2$ and $\gk_+\in(1/2,1]$.
Indeed, $$n-N+k=k+\{2\gk_+\}-2$$
and (as $\gk_+>0$) $\{2\gk_+\}\geq 1$. If $k=2$ then
$\gk_+>1/2$ and \consy $\{2\gk_+\}\geq 2$. These facts imply our assertion.

We also note that $\gk_+$ is strictly
increasing relative to $\gl_A$ and
\begin{equation}\label{gk+}
\gk_+\begin{cases}=1, &\text{if $D_A=\BBR^N_+$,}\\
<1, &\text{if $D_A\subsetneqq\BBR^N_+$,}\\
>1, &\text{if $D_A\supsetneqq\BBR^N_+$}.\end{cases}
\end{equation}
Finally we observe that  $\gg:=\gl_A-(N-k)\gk_+>0$ (see \eqref{EV5-0}) and, by
\eqref{kappa2} and \eqref{qk}:
\begin{equation}\label{q*c1}
 \gg=\gk_+^2+(k-2)\gk_+,\q q^*_c=1+ \frac{-(k-2)+\sqrt{(k-2)^2+4\gg}}{\gg}.
\end{equation}
Therefore $q^*_c$ is strictly decreasing relative to $\gg$ and \consy also relative to $\gk_+$.

The proof of the theorem is based on the following important result proved in \cite [1.14.4.]{Tri}
\bprop {repr}Let $1<q<\infty$ and $s>0$. Then for any bounded Radon measure $\gm$ in $\BBR^{n-1}$
there holds
\begin{equation}\label{P4}
I(\gm)=\int_
{\BBR_+^{n}}\abs{\BBK_{n}[\gm](y)}^qe^{-y_{1}}y_{1}^{sq-1}dy\approx
\norm{\gm}^q_{B^{-s,q}(\BBR^{n-1})}.
\end{equation}
\es

\smallskip

In the first part of the proof we derive inequalities comparing
$I(\mu)$ and $J^{A,R}(\mu)$. Actually, it is useful to consider  a
slightly more general expression than $I(\mu)$, namely:
\begin{equation}\label{Inu}
I_{\gn,\gs}^{m,j}(\gm):=\int_
{\BBR_+^{m+j}}\abs{\int_{\BBR^{m}}{}\myfrac{y_{1}d\gm(z)}{\left(y_1^2+
|\wtl y|^2+|y''-z|^2\right)^{\gn/2}}}^qe^{-y_{1}}y_{1}^{\gs q-1}dy,
\end{equation}
where $\nu$ is an arbitrary number \sth $\nu>m$, $j\geq 1$ and $\gs>0$. A point $y\in \BBR_+^{m+j}$
is written in the form $y=(y_1,\wtl y, y'')\in \BBR_+\ti \BBR^{j-1}\ti\BBR^m$. We assume that
 $\gm$  is supported in $\BBR^m$. Note that, 
 \begin{equation}\label{Imns}
    I(\mu)=\gg_{n}^qI_{n,s}^{m,j} \txt{where}
  m=N-k,\q j=n-m=n-N+k.
 \end{equation}
Put
\begin{equation}\label{Fnum}
 F_{\nu,m}[\mu](\tau):=\int_{\BBR^{m}}\abs{\int_{\BBR^{m}}
\myfrac{d\gm(z)}{(\tau^2+|y''-z|^2)^{\gn/2}}}^q dy''\forevery
\tau\in [0,\infty).
\end{equation}
With this notation, if $j\geq 2$ then
\begin{equation}\label{InuF}
I_{\gn,\gs}^{m,j}(\gm):=\int_0^\infty\int_
{\BBR^{j-1}} F_{\nu,m}[\mu](\sqrt{y_1^2+|\wtl y|^2}\,)e^{-y_{1}}y_{1}^{(\gs+1) q-1}d\wtl y\,dy_1
\end{equation}
and if $j=1$
\begin{equation}\label{InuF1}
I_{\gn,\gs}^{m,1}(\gm):=\int_0^\infty
 F_{\nu,m}[\mu](y_1)e^{-y_{1}}y_{1}^{(\gs+1) q-1}\,dy_1
\end{equation}


\blemma {repl}
Assume that $m<\gn$, $0<\gs$, $2\leq j$ and
$1<q<\infty$.
Then there exists a positive constant $c$, depending on $m,j,\nu,\gs,q$, such
that, for every bounded Borel measure $\gm$ with support in $\BBR^m$:

\begin{equation}\label{Ck2}
\rec{c}\int_0^\infty\,F_{\nu,m}[\mu](\tau)h_{\gs,j}(\tau)d\tau\leq
\,I_{\gn,\gs}^{m,j}(\gm)\leq
c\int_0^\infty\,F_{\nu,m}[\mu](\tau)h_{\gs,j}(\tau)d\tau,
\end{equation}
where $F_{\nu,m}$ is given by \eqref{Fnum} and, for every $\tau>0$,
\begin{equation}\label{hsj}
  h_{\gs,j}(\tau)=\begin{cases}\myfrac{\tau^{(\gs+1)q+j-2}}{(1+\tau)^{(\gs+1)q}}, &\text{if $j\geq2$,}\\[4mm]
  e^{-\tau}\tau^{(\gs+1)q-1},  &\text{if $j=1$.}\end{cases}
\end{equation}

\es
\Proof There is nothing to prove in the case $j=1$. Therefore we assume that $j\geq2$.

We use the notation $y=(y_{1},\wtl
y,y'')\in\BBR\ti\BBR^{j-1}\ti\BBR^{m}$. 
The integrand in
 \eqref{InuF} depends only on $y_1$ and $\gr:=|\wtl y|$. Therefore,
 $I_{\gn,\gs}^{m,j}$ can be written in the form
$$\BAL
I_{\gn,\gs}^{m,j}(\gm)
=c_{m,j}\int_0^\infty\int_0^\infty
F_{\nu,m}[\mu](\sqrt{y_1^2+\gr^2})e^{-y_1}y_1^{(\gs
+1)q-1}\,dy_1\gr^{j-2}d\gr.\EAL
$$
We substitute
 $y_1=(\tau^2-\gr^2)^{1/2}$, then  change the order of integration
 and finally  substitute $\gr=r\tau$. This yields,
$$\BAL &c_{m,j}^{-1}I_{\gn,\gs}^{m,j}(\gm)\\
&=\int_0^\infty\int_\gr^\infty
F_{\nu,m}[\mu](\tau)\gr^{j-2}e^{-\sqrt{\tau^2-\gr^2}}(\tau^2-\gr^2)^{(\gs +1)q/2-1}\tau\,d\tau\,d\gr\\
&=\int_0^\infty\int_0^\tau\,F_{\nu,m}[\mu](\tau)\gr^{j-2}e^{-\sqrt{\tau^2-\gr^2}}(\tau^2-\gr^2)^{(\gs
+1)q/2-1}\tau\,d\gr\,d\tau\\
&=\int_0^\infty\int_0^1\,F_{\nu,m}[\mu](\tau)\tau^{j-2+(\gs
+1)q}e^{-\tau\sqrt{1-r^2}}f(r) dr\,d\tau,
 \EAL$$
 where
 $$f(r)= r^{j-2}(1-r^2)^{(\gs +1)q/2-1}.$$
We denote
$$I_\gs^j(\tau)=\int_0^1\,e^{-\tau\sqrt{1-r^2}}f(r) dr,$$
so that
\begin{equation}\label{Isn}
 I_{\gn,\gs}^{m,j}(\gm)=c_{m,j}\int_0^\infty F_{\nu,m}[\mu](\tau)\tau^{j-2+(\gs
+1)q}I_\gs^j(\tau)d\tau.
\end{equation}

To complete the proof we estimate $I^j_\gs$. Since $j\geq2$,
$f\in L^1(0,1)$  and $I^j_\gs$ is continuous in
$[0,\infty)$ and positive everywhere. Hence, for every $\ga>0$,
there exists a positive constant $c_\ga=c_\ga(\gs)$ \sth
\begin{equation}\label{I*1}
\rec{c_\ga}\leq I_\gs^j\leq c_\ga \txt{in} [0,\ga).
\end{equation}

Next we estimate $I^j_\gs$ for large $\tau$. Since $j\geq2$,
$$I^j_\gs\leq 2^{(\gs +1)q/2-1}\int_0^1(1-r)^{(\gs +1)q/2-1}e^{-\tau\sqrt{1-r}}dr.$$
Substituting $r=1-t^2$ we obtain,
\begin{equation}\label{I*2'}
 I^j_\gs\leq 2^{(\gs +1)q/2}\int_0^1 t^{(\gs +1)q-1}e^{-t\tau}dt=c(\gs,q)\tau^{-(\gs +1)q}.
\end{equation}

On the other hand, if $\tau\geq 2$,
\begin{equation}\label{I*3}\BAL
I^j_\gs(\tau)&= \int_0^1(1-t^2)^{(j-3)/2}t^{(\gs+1)q-1} e^{-\tau t}dt\\
 &=\tau^{-(\gs+1)q} \int_0^\tau
 (1-(s/\tau)^2)^{(j-3)/2}s^{(\gs+1)q-1}e^{-s}ds\\
 &\geq \tau^{-(\gs+1)q}2^{-(j-3)}
 \int_0^1s^{(\gs+1)q-1}e^{-s}ds.
 \EAL\end{equation}
Combining \eqref{Isn} with \eqref{I*1}--\eqref{I*3} we obtain
\eqref{Ck2}. \qed

\par Next we derive an estimate in which integration over $\BBR^n_+=\BBR_+^j\ti\BBR^m$
is replaced by integration over a bounded domain, for measures supported
in  a fixed bounded subset of $\BBR^m$.

Let $B_R^j(0)$ and $B_R^m(0)$ denote the balls of radius $R$
centered at the origin, in $\BBR^j$ and $\BBR^m$ respectively. Denote
\begin{equation}\label{FnumR}
  F^R_{\nu,m}[\mu](\tau)=\int_{B^m_R}\abs{\int_{\BBR^{m}}
\myfrac{d\gm(z)}{(\tau^2+|y''-z|^2)^{\gn/2}}}^q dy''\forevery
\tau\in [0,\infty)
\end{equation}
and, if $j\geq 2$,
\begin{equation}\label{InuR}\BAL
I_{\gn,\gs}^{m,j}(\gm;R)=\int_{B^j_R\cap\{0<y_1\}}
F^R_{\nu,m}[\mu](\sqrt{y_1^2+|\wtl y|^2}\,)e^{-y_{1}}y_{1}^{\gs q-1}d\wtl y\,dy_1.
\EAL\end{equation}
where $(y_1,\wtl y)\in \BBR\ti\BBR^{j-1}$. If $j=1$ we denote,
\begin{equation}\label{InuR1}\BAL
I_{\gn,\gs}^{m,1}(\gm;R)=\int_0^R
F^R_{\nu,m}[\mu](y_1)e^{-y_{1}}y_{1}^{\gs q-1}\,dy_1.
\EAL\end{equation}

Similarly to \rlemma{repl} we obtain,
\blemma{replR} If $j\geq1$,
there exists a positive constant $c$ such
that, for any bounded Borel measure $\gm$ with support in $\BBR^m\cap B_R$

\begin{equation}\label{Ck2R}
c^{-1}\int_0^R\,F^R_{\nu,m}[\mu](\tau)h_{\gs,j}(\tau)d\tau\leq
\,I^{m,j}_{\gn,\gs}(\gm;R)\leq
c\int_0^R\,F^R_{\nu,m}[\mu](\tau)h_{\gs,j}(\tau)d\tau
\end{equation}
with $h_{\gs,j}$ as in \eqref{hsj}.
\es

\nind\Proof In the case $j=1$ there is nothing to prove .
Therefore we assume that $j\geq 2$.

From \eqref{InuR} we obtain,
$$\BAL
I^{m,j}_{\gn,\gs}(\gm;R)=c_{m,j}\int_0^R\int_0^{\sqrt{R^2-\gr^2}}
F^R_{\nu,m}[\mu](\sqrt{y_1^2+\gr^2})e^{-y_1}y_1^{(\gs
+1)q-1}dy_1\gr^{j-2}d\gr.\EAL
$$
Substituting $y_1=(\tau^2-\gr^2)^{1/2}$, then  changing the order of integration
 and finally  substituting $\gr=r\tau$ we obtain,

$$c_{m,j}^{-1}I^{m,j}_{\gn,\gs}(\gm;R)=\int_0^R\int_0^1\,F^R_{\nu,\mu}[\mu](\tau)\tau^{j-2+(\gs
+1)q}e^{-\tau\sqrt{1-r^2}}f(r) dr\,d\tau.$$
 where
 $$f(r)= r^{j-2}(1-r^2)^{(\gs +1)q/2-1}.$$
The remaining part of the proof is the same as for \rlemma{repl}.
\qed
\blemma{reduc}Let $1<q$, $0<\gs$ and assume that $m<\nu q$ and
$0\leq j-1<\nu$.
Then there exists a positive constant $\bar c$, depending on $j,m,q,\gs,\nu$, such
that, for every $R\geq1$ and every bounded Borel measure $\gm $ with support in $B_{R/2}(0)\cap
\BBR^m$,
\begin{equation}\label{R1}\BAL
\abs{\int_0^\infty\,F_{\nu,m}[\mu](\tau)h_{\gs,j}(\tau)d\tau-
\int_0^R\,F_{\nu,m}^R[\mu](\tau)h_{\gs,j}(\tau)d\tau}&\\[2mm]
\leq \bar c R^{(\gs +1-\nu)q+m+j-1}\norm{\gm}_{\GTM}^q&
\EAL\end{equation}
 with $h_{\gs,j}$ as in \eqref{hsj}. \es
\Proof
We estimate,
\begin{equation}\label{J1+J2}
\BAL
&\abs{\int_0^\infty\,F_{\nu,m}[\mu](\tau)h_{\gs,j}(\tau)d\tau-\int_0^R\,F_{\nu,m}^R[\mu](\tau)h_{\gs,j}(\tau)d\tau}\leq\\
&\int_R^\infty \abs{F_{\nu,m}[\mu]}(\tau)h_{\gs,j}(\tau)d\tau+
\int_0^R\,\abs{F_{\nu,m}[\mu]-F_{\nu,m}^R[\mu]}(\tau)h_{\gs,j}(\tau)d\tau.
\EAL
\end{equation}
For every $\tau>0$,
\begin{equation}\label{Fnm<}
\abs{F_{\nu,m}[\mu]}(\tau)\leq \tau^{-\nu q}\norm{\mu}_{\GTM}^q.
\end{equation}
Since $j-1<\nu q$, it follows that
\begin{equation}\label{J1-temp}
\BAL \int_R^\infty \abs{F_{\nu,m}[\mu]}(\tau)h_{\gs,j}(\tau)d\tau
&\leq \norm{\mu}_{\GTM}^q\int_R^\infty \tau^{-\nu q}h_{\gs,j}(\tau)d\tau\\
&\leq c(\gs,q)\norm{\mu}_{\GTM}^q \int_R^\infty \myfrac{\tau^{(\gs+1)q+j-2-\nu q}}{(1+\tau)^{(\gs+1)q}}d\tau\\
&\leq \frac{c(\gs,q)}{\nu q-j+1}
\norm{\mu}_{\GTM}^q R^{j-1-\nu q}.\EAL
\end{equation}
Since, by assumption, $\supp\mu\sbs B_{R/2}$, we have
\begin{equation}\label{J2-temp}
\BAL \int_0^R\,&\abs{F_{\nu,m}[\mu]-F_{\nu,m}^R[\mu]}(\tau)h_{\gs,j}(\tau)d\tau\\
&\leq\int_0^R\,\int_{|y''|>R}\abs{\int_{\BBR^{m}}
\myfrac{d\gm(z)}{(\tau^2+|y''-z|^2)^{\gn/2}}}^q dy''h_{\gs,j}(\tau)d\tau\\
&\leq \norm{\mu}_{\GTM}^q \int_0^R\int_{|\gz|>R/2}(|\tau^2+|\gz|^2)^{-\nu q/2}\,d\gz\,h_{\gs,j}d\tau\\
&\leq c(m,q)\norm{\mu}_{\GTM}^q \int_0^R\int_{R/2}^\infty(\tau^2+\gr^2)^{-\nu q/2}\gr^{m-1}\,d\gr\,h_{\gs,j}d\tau\\
&\leq c(m,q)\norm{\mu}_{\GTM}^q \int_0^R\,\tau^{m-\nu q}\int_{R/2\tau}^\infty(1+\eta^2)^{-\nu q/2}\eta^{m-1}\,d\eta\,h_{\gs,j}d\tau\\
&\leq \frac{c(m,q)}{\nu q-m}\norm{\mu}_{\GTM}^q R^{m-\nu q}\int_0^R\tau^{(\gs+1)q+j-2}\,d\tau\\
&\leq \frac{c(m,q)}{(\nu q-m)((\gs+1)q+j-1)}\norm{\mu}_{\GTM}^q R^{(\gs+1)q+j-1+m-\nu q}.\EAL
\end{equation}
Combining \eqref{J1+J2}--\eqref{J2-temp} we obtain
 \eqref{R1}.
\qed
\bcor{reduc}
For every $R>0$ put
\begin{equation}\label{Jns}
  J_{\nu,\gs}^{m,j}(\mu;R):=\int_0^RF^R_{\nu,m}[\mu](\tau) \tau^{(\gs+1)q+j-2}d\tau.
\end{equation}
Then
\begin{equation}\label{Est1-Jns}\BAL
\rec{c}I_{\gn,\gs}^{m,j}(\gm)-\bar c R^{\gb}\norm{\gm}_{\GTM}^q\leq
J_{\nu,\gs}^{m,j}(\mu;R)
\leq c R^{(\gs+1)q}I_{\gn,\gs}^{m,j}(\gm),\\
\gb=(\gs +1-\nu)q+j+m-1, \EAL
\end{equation}
for every $R>1$ and  every bounded Borel measure $\gm $ with support in $B^m_{R/2}(0):=B_{R/2}(0)\cap
\BBR^m$.
\es
\Proof This is an immediate consequence of \rlemma{reduc} and \rlemma{repl}.
\qed
\blemma{Est2-Jns} Let $m,j$ be positive integers \sth $j\geq 1$ and
let $1<q$, $0<\gs$. Put $n:=m+j$.

Then there exist positive constants $c,\bar c$, depending on $j,m,q,\gs$, such
that, for every $R>1$ and every measure $\gm\in \GTM_+(B^m_{R/2}(0))$,
\begin{equation}\label{Est2-Jns}\BAL
&\rec{c}\norm{\gm}^q_{B^{-\gs,q}(\BBR^{n-1})}-\bar c R^{q\left(\gs-\frac{n-1}{q'}\right)}\norm{\gm}_{\GTM}^q\leq
J_{n,\gs}^{m,j}(\mu;R)\\
&\leq c R^{(\gs+1)q}\norm{\gm}^q_{B^{-\gs,q}(\BBR^{n-1})}.\EAL
\end{equation}

If  $\gs<\frac{n-1}{q'}$,  there exists $R_0>1$ \sth, for all $R>R_0$
\begin{equation}\label{Est3-Jns}\BAL
\rec{2c}\norm{\gm}^q_{B^{-\gs,q}(\BBR^{n-1})}\leq
J_{n,\gs}^{m,j}(\mu;R).
\EAL
\end{equation}
If $\gs=\frac{n-1}{q'}$ then,  there exists $a>0$ \sth the inequality remains valid
for measures $\mu$ \sth $\mathrm{diam}(\supp\mu)\leq a$.

If, in addition, $\frac{j-1}{q'}<\gs$ then
\begin{equation}\label{Est4-Jns}\BAL
\rec{2c}\norm{\gm}^q_{B^{-s,q}(\BBR^{m})}\leq
J_{n,\gs}^{m,j}(\mu;R)
\leq c R^{(\gs+1)q}\norm{\gm}^q_{B^{-s,q}(\BBR^{m})},\EAL
\end{equation}
where $s:=\gs-\frac{j-1}{q'}$.
\es
\Remark Assume that $\mu\geq0$. Then:\\
 (i) If $\mu\in B^{-\gs,q}(\BBR^{n-1})$ and $\frac{j-1}{q'}\geq \gs$ then $\mu(\BBR^m)=0$.\\
(ii) If $\mu\in B^{-s,q}(\BBR^m)$ and $\gs> (n-1)/q'$ then $s> m/q'$ and therefore
$B^{s,q'}(\BBR^m)$ can be embedded in $C(\BBR^m)$.\2
\Proof Inequality \eqref{Est2-Jns} follows from \eqref{Est1-Jns} and \rprop{repr} (see also \eqref{Imns}).

For positive measures $\mu$,
$$\norm{\gm}_{\GTM}=\mu(\BBR^{n-1})\leq \norm{\gm}^q_{B^{-\gs,q}(\BBR^{n-1})}.$$
Therefore, if $\gs<\frac{n-1}{q'}$,  \eqref{Est2-Jns}   implies that
 there exists $R_0>1$ \sth
\eqref{Est3-Jns} holds for all $R>R_0$.

If $\gs=\frac{n-1}{q'}$ \eqref{Est2-Jns}   implies that
$$\rec{c}\norm{\gm}^q_{B^{-\gs,q}(\BBR^{n-1})}-\bar c \norm{\gm}_{\GTM}^q\leq
J_{n,\gs}^{m,j}(\mu;R).$$
But if $\mu$ is a positive bounded measure \sth $\mathrm{diam}(\supp\mu)\leq a$ then
$$\norm{\gm}_{\GTM}/\norm{\gm}^q_{B^{-\gs,q}(\BBR^{n-1})}\to 0 \txt{as $a\to 0$.}$$
The last inequality follows from the
imbedding theorem for Besov spaces according to which
there exists a
continuous trace operator $T:B^{\gs,q'}(\BBR^{n-1})\mapsto
B^{s,q'}(\BBR^{m})$ and a continuous lifting
$T':B^{s,q'}(\BBR^{m})\mapsto B^{\gs,q'}(\BBR^{n-1})$ where $s=\gs-\frac{n-m-1}{q'}$.
\qed

If $\nu\in\BBN$ and $\gs=s+\frac{\nu-m-1}{q'}$,
\begin{equation*}\BAL
 J_{\nu,\gs}^{m,\nu-m}(\mu;R)&= \int_0^RF^R_{\nu,m}[\mu](\tau) \tau^{(\gs+1)q+\nu-m-2}\,d\tau\\
& = \int_0^RF^R_{\nu,m}[\mu](\tau) \tau^{(s+\nu-m)q-1}\,d\tau.
\EAL\end{equation*}
However, if $\mu$ is positive, the expression
\begin{equation}\label{Mns}
M_{\nu,s}^{m}(\mu;R):=\int_0^RF^R_{\nu,m}[\mu](\tau) \tau^{(s+\nu-m) q-1}\,d\tau,
\end{equation}
is meaningful for any real $\nu>m$ and $s>0$. Furthermore, as shown  below, the results stated in \rlemma{Est2-Jns}
can be extended to this general case.
\bth{general-nu} Let $1<q$, $\nu\in \BBR$ and $m$ a positive
integer. Assume that $1\leq\nu-m$ and $0<s<m/q'$. Then there
exists a positive constant $c$ \sth, for every bounded positive
measure $\mu$ supported in $\BBR^m\cap B_{R/2}(0)$, $R>1$,
\begin{equation}\label{Est-Mns<}\BAL
\rec{c}\norm{\gm}^q_{B^{-s,q}(\BBR^{m})}\leq M_{\nu,s}^{m}(\mu;R)
\leq c R^{(s+\nu-m) q+1}\norm{\gm}^q_{B^{-s,q}(\BBR^{m})}.\EAL
\end{equation}
This also holds when $s=m/q'$, provided that the diameter of $\supp
\mu$ is sufficiently small. \es
\Proof If $\nu$ is an integer and $j:=\nu-m$ then this statement is part of \rlemma{Est2-Jns}.
Indeed the condition $s>0$ means that $\gs=s+\frac{j-1}{q'}>\frac{j-1}{q'}$ and the condition
$s< m/q' $ means that $\gs<\frac{n-1}{q'}$.

Therefore we assume that
$\nu\nin\BBN$. Let $n:=\{\nu\}$ and $\gth:=n-\nu$ so that $0<\gth<1$. 
Our assumptions imply that $1\leq n-m-1$ because (as $\nu$ is not an integer) $\nu-m>1$ and \consy
$n-m\geq 2$.

 If $a,b$ are positive numbers, put
\begin{equation*}\label{abt}
   A_\nu:=\frac{a^{(s+\nu-m)q-1}}{(a^2+b^2)^{\nu q/2}}.
\end{equation*}
Obviously $A_\nu$ decreases as $\nu$ increases. Therefore, $A_n\leq A_\nu\leq A_{n-1}$ which in turn implies,
$$M_{n,s}^{m}\leq M_{\nu,s}^{m}\leq M_{n-1,s}^{m}.$$
By \rlemma{Est2-Jns}, the assertions of the theorem are valid in the case that $\nu=n$ or $\nu=n-1$.
Therefore the previous
inequality implies that the assertions hold for any real $\nu$ subject to the conditions imposed.
\qed

By \eqref{k8'},
$$J^{A,R}=\int_0^R F_{\nu,m}^R(\tau)\tau^{(q+1)\gk_{+}+k-1}d\tau,$$
where $m=N-k$ and $\nu=N-2+2\gk_+$. 
\Consy, by \eqref{Jns},
$$J^{A,R}=M_{\nu,s}^m$$
where $s$ is determined by,
$$(s+\nu-m) q-1=(q+1)\gk_{+}+k-1,\q k=\nu-m+2-2\gk_+.$$
It follows that
$$sq=-(k-2+2\gk_+)q+(q+1)\gk_{+}+k=k(1-q)+2q-\gk_+(q-1)$$
and therefore
$$s=2-\frac{k+\gk_+}{q'}.$$

\medskip




\medskip

\noindent{\it Proof of \rth{main1}.}

Put
\begin{equation}\label{nusm}
 \nu:=N-2+2\gk_+,\q  s:=2-\frac{\gk_{+}+k}{q'},\q m:=N-k.    
\end{equation}

Recall that in the case $k=2$ we have $\gk_+>1/2$. Therefore
\begin{equation}\label{nu>m-1}
  \nu-m-1=k-3+2\gk_+>0.
\end{equation}
Furthermore,
$$(s+\nu-m) q-1=(q+1)\gk_{+}+k-1,\q k=\nu-m+2-2\gk_.$$
Thus
$$J^{A,R}=\int_0^R F_{\nu,m}^R(\tau)\tau^{(q+1)\gk_{+}+k-1}d\tau=M_{\nu,s}^m.$$

Next we show that $ 0<s\leq m/q'$.
More precisely we prove
\begin{equation}\label{0<s<m/q'}
 0<s\leq m/q' \;\iff\;q_c\leq q<q^*_c.
\end{equation}

Let $\mu$ be a bounded non-negative Borel measure in $B^{-s,q}(\BBR^m)$. If $s\leq 0$,
$B^{-s,q}(\BBR^m)\sbs L^q(\BBR^m)$. Therefore, in this case, every bounded Borel measure
on $\BBR^m$ is admissible i.e. satisfies \eqref{admpk}. \Consy, by \rprop{admpk}, $q<q_c$.
As we assume $q\geq q_c$ it follows that $s>0$.

If, $s>0$ and $sq'-m\geq0$ then $C_{s,q'}(K)=0$
for every compact subset of $\BBR^m$ and \consy $\mu(K)=0$ for any such set.
Conversely, if $sq'-m<0$ then there exist non-trivial positive bounded measures in $B^{-s,q}(\BBR^m)$.
Therefore, by \rprop{admp0},  $sq'<m$ if and only if $q<q^*_c$.

In conclusion, $0<s\leq m/q'$ and $\nu-m\geq 1$; therefore \rth{main1} is a consequence of
\rth{general-nu}.
\qed

\noindent\Remark Note that the critical exponent for the imbedding
of $B^{2-\frac{\gk_{+}+k}{q'},q'}(\BBR^{N-k})$ into $C(\BBR^{N-k})$
is again
$$ q=q_{c}=\myfrac{N+\gk_{+}}{N+\gk_{+}-2}.
$$



\mysection{ Supercritical equations in a polyhedral domain} In this
section $q$ is a real number larger than $1$ and $P$ an N-dim
polyhedral domain as described in subsection 6.1. Denote by
$\set{L_{k,j}:k=1,\dots,N,\;j=1,\dots, n_k}$ the family of faces,
edges and vertices of $P$. In this notation, $L_{1,j}$ denotes one
of the open faces of $P$; for $k=2,\dots,N-1$, $L_{k,j}$ denotes a
relatively open $N-k$-dimensonal edge and $L_{N,j}$ denotes a
vertex. For $1\leq k<N$, the $(N-k)$ dimensional space which
contains $L_{k,j}$ is denoted by $\BBR^{N-k}_j$. If $1<k<N$, the
cylinder of radius $r$ around the axis $\BBR^{N-k}_j$ will be
denoted by $\Gg^\infty_{k,j,r}$ and  the subset $A_{k,j}$ of
$S^{k-1}$ is defined by
$$\lim_{r\to0}\rec{r}(\prt\Gg^\infty_{k,j,r}\cap P)=L_{k,j}\ti A_{k,j}.$$
$A_{k,j}$ is the 'opening' of $P$ at the edge $L_{k,j}$.
For $k=N$ we replace in this definition the cylinder $\Gg^\infty_{N,j,r}$ by the ball $B_r(L_{N,j})$.
For $1<k\leq N$ and $A=A_{k,j}$ we use  $d_{A}$ as an alternative notation for $\BBR^{N-k}_j$ and denote by $D_A$ the k-dihedron  with edge $d_A$ and opening $A$ as in subsection 6.1 (with  $S_A$  defined as in \eqref{dom}). For $k=1$,
$D_A$ stands for the half space $\BBR^{N-1}_j\ti (0,\infty)$.

In what follows we denote by $\GTM_q^\Gw$ the set of bounded measures $\mu$ on the boundary of a \Lip domain $\Gw$ \sth
the \bvp
\begin{equation}\label{bvp8.1}
  -\Gd u+u^q=0 \txt{in $\Gw$,} u=\mu \txt{on $\prt \Gw$}
  \end{equation}
possesses a solution. A measure $\mu$ in this space is called a {\em q-good measure}.

The following statements can be proved in the same way as in the
case of smooth domains. For the proof in that case see \cite{MV2}.

\medskip
\nind{\bf I.} $\GTM_q^\Gw$ is a linear space and
$$\mu\in \GTM_q^\Gw\iff  |\mu|\in \GTM_q^\Gw.$$

\nind{\bf II.} If $\set{\mu_n}$ is an increasing \seq of measures in $\GTM_q^\Gw$ and $\mu:=\lim \mu_n$ is a finite measure then $\mu\in \GTM_q^\Gw$.

\bprop{edge condition} Let $\mu$ be a  bounded measure on $\prt P$.
($\mu$ may be a signed measure.)  For $i=1,\dots,N,\;j=1,\dots,
n_i$, we define the measure $\mu_{k,j}$ on $d_{A_{k,j}}$ by,
$$\mu_{k,j}=\mu\txt{on ${L_{k,j}}$,}\q \mu_{k,j}=0 \txt{on $d_{A_{k,j}}\sms {L_{k,j}}$}.$$

Then $\mu\in \GTM_q^P$, i.e., problem
\begin{equation}\label{Pbvp}
  -\Gd u+u^q=0 \txt{in $P$,} u=\mu \txt{on $\prt P$}
\end{equation}
possesses a solution, if and only if, $\mu_{k,j}$ is a q-good
measure relative to $D_{A_{k,j}}$ for all $(k,j)$ as above. \es

\Proof In view of statement {\bf I} above, it is sufficient to prove
the proposition in the case that $\mu$ is non-negative. This is
assumed hereafter. If $\mu\in \GTM_q^P$ then any measure $\nu$ on
$\prt P$ \sth $0\leq \nu\leq \mu$ is a q-good measure relative to
$P$. Therefore
$$\mu\in \GTM_q^P \Lra \mu'_{k,j}:=\mu\chi\indx{L_{k,j}}\in \GTM_q^P.$$
 Assume that $\mu\in \GTM_q^P$ and let $u_{k,j}$ be the solution of \eqref{Pbvp} when $\mu$ is replaced by $\mu'_{k,j}$. Denote by $u'_{k,j}$  the extension of $u_{k,j}$ by zero to the k-dihedron $D_{A_{k,j}}$. Then $u'_{k,j}$ is a subsolution of \eqref{eq-q} in
$D_{A_{k,j}}$ with boundary data $\mu_{k,j}$. In the present case
there always exists a supersolution, e.g. the maximal solution of
\eqref{eq-q} in $D_{A_{k,j}}$ vanishing outside $d_{A_{k,j}}\sms\bar
L_{k,j}$. Therefore there exists a solution $v_{k,j}$ of this
equation in $D_{A_{k,j}}$ with boundary data $\mu_{k,j}$, i.e.,
$\mu_{k,j}$ is q-good relative to $D_{A_{k,j}}$.

Next assume that $\mu\in \GTM(\prt P)$ and that $\mu_{k,j}$ is
q-good relative to $D_{A_{k,j}}$ for every $(k,j)$ as above. Let
$v_{k,j}$ be the solution of \eqref{eq-q} in $D_{A_{k,j}}$ with
boundary data $\mu_{k,j}$. Then $v_{k,j}$ is a supersolution of
problem \eqref{Pbvp} with $\mu$ replaced by $\mu'_{k,j}$ and \consy
there exists a solution $u_{k,j}$ of this problem. It follows that
$$w:=\max\set{u_{k,j}: k=1,\dots,N,\;j=1,\dots, n_k}$$
is a subsolution  while
$$\bar w:=\sum_{k=1,\dots,N,\;j=1,\dots, n_k}u_{k,j}$$
is a supersolution of \eqref{Pbvp}. \Consy there exists a solution of this problem, i.e., $\mu\in \GTM_q^P.$
\qed

\subsection{Removable singular sets and 'good measures', I}

\bprop{admi-conv} Let $A$ be a \Lip domain on $S^{k-1}$, $2\leq
k\leq N-1$, and let $D_A$ be the k-dihedron with opening $A$.
 Let $\gm\in\GTM(\prt D_{A})$ be a positive measure with compact support contained in $d_{A}$
 (= the edge of $D_A$). Assume that $\mu$ is q-good relative to $D_A$. Let $R>1$ be large enough so that
$\supp\mu\sbs B^{N-k}_R(0)$  and let $u$ be the solution of
\eqref{eq-q} in $D_A^R$ with trace $\mu$ on $d_A^R$ and trace zero
on $\prt D_A^R\sms d_A^R$. Then:

 \medskip
\nind{\rm (i)} For every non-negative $\eta\in
C_0^\infty(B_{3R/4}^{N-k}(0))$, 
\begin{equation}\label{inverse-ineq}\BAL
&\left( \int_{d_{A}^R}\eta^{q'} d\gm \right)\leq cM^{q'}
\int_{D_{A}^R}{}u^q\gr dx+\\
&\q+ cM^{q'}\left(\int_{D_{A}^R}{}u^q\gr dx\right)^{\rec{q}}
\left(1+M^{-1}\norm{\eta}_{L^{q'}(d_{A}^R)}\right). \EAL
\end{equation}
where $M=\norm{\eta}\indx{L^\infty}$ and $\gr$ is the first
eigenfunction of $-\Gd$ in $D_A^R$ normalized by $\gr(x_0)=1$ at
some point $x_0\in D_A^R$. The constant $c$ depends only on $N,q,k,
x_0,\gl_1,R$ where $\gl_1$ is the first eigenvalue.

\smallskip
 \nind{\rm (ii)} For any compact set
 $E\subset d_{A}$,

\begin{equation}\label{BVP-A2}
C^{N-k}_{s,q}(E)=0\Longrightarrow \mu (E)=0, \q s=2-\frac{\gk_{+}+k}{q'},
\end{equation}
where $C^{N-k}_{s,q}$ denotes the Bessel capacity with the indicated
indices in $\BBR^{N-k}$. \es

\Remark If we replace $D_A^R$ by $D_A^{\tl R,  R}=D_A\cap B_{\tl
R}^k(0)\cap B_{ R}^{N-k}(0)$,  $\tl R>1$,  then the constant $c$ in
(i) depends on $\tl R$ but {\em not on $ R$.}

\medskip
\Proof  We identify $d_A$ with $\BBR^{N-k}$ and  use the notation
$$x=(x',x'')\in \BBR^k\ti\BBR^{N-k},\q  y=|x'|.$$

 Let $\eta\in C_{0}^\infty(\BBR^{N-k})$ and let $R$ be large enough so that
$\supp\eta\sbs B^{N-k}_{R/2}(0)$.
 Let $w=w_R(t,x'')$   be the solution of the following problem in $\BBR_+\ti B^{N-k}_{R}(0)$:
  \begin{equation}\label{heat1}\BAL
 \prt_{t}w-\Gd_{x''} w&=0&&\text{in }\BBR^+\ti B_R^{N-k}(0),\\
 w(0,x'')&=\eta (x'')&&\text{in }B_R^{N-k},\\
 w(t,x'')&=0&&\text{on }\prt B^{N-k}_{R}(0).
 \EAL\end{equation}
Thus $w_R(t,\cdot)=S_R(t)[\eta]$ where $S_R(t)$ is the semi-group operator corresponding to the above problem. Denote,
   \begin{equation}\label{heat2}
    H_R[\eta](x',x'')=w_R(|x'|^2,x'')=S_R(y^2)[\eta](x''),\q y:=|x'|.
   \end{equation}

We assume, as we may, that $R>1$. Let $\gr^R$ be the first eigenfunction of $-\Gd_{x''}$ in the ball $B^{N-k}_{R}(0)$ normalized by $\gr^R(0)=1$ and let $\gr_A$ be the first eigenfunction of $-\Gd_{x'}$ in $C_{A}$
(where $C_{A}$ denotes the cone with opening $A$ in $\BBR^k$) normalized so that $\gr_A(x'_0)=1$ at some point $x'_0\in S_A$. Then $\gr^R\gr_{A}$ is the first eigenfunction of $-\Gd$ in $\set{x\in D_{A}:|x''|<R}$. Note that $\gr^R\leq 1$ and $\gr^R\to 1$  as $R\tin$ in $C^2(I)$ for any bounded set $I\sbs \BBR^{N-k}$.

Let $h\in C^\infty(\BBR)$ be a monotone decreasing function \sth $h(t)=1$ for $t<1/2$ and $h(t)=0$ for $t>3/4$. Put
$$\psi_R(x')=h(|x'|/R)$$
and
\begin{equation}\label{gz1}
   \gz_R:=\gr_{A}\psi_R H_R[\eta]^{q'}.
\end{equation}
If $\gr_A^R$ is the first eigenfunction (normalized at $x_0$) of $D_A^R:=D_A\cap \Gg_R$ ($\Gg_R$
as in \eqref{GgR}) then
\begin{equation}\label{rhoAR}
  \gr_A\psi_R\leq c\gr_A^R
\end{equation}
 and $\gr^R\gr_A^R$ is the first eigenfunction in $D_A^R$.

Hereafter we shall drop the index $R$  in $\gz_R, H_R, w_R$ but keep it in the other notations in order
 to avoid confusion.

We shall verify that $\gz\in D_A^R$. To this purpose we compute,
  \begin{equation}\label{heat2'}\BAL
\Gd\gz=&-\gl_{1}(\gr_A\psi_R) H[\eta]^{q'}+(\gr_A\psi_R)\Gd H[\eta]^{q'}+2\nabla(\gr_A\psi_R)\cdot\nabla H[\eta]^{q'}\\
=&-\gl_{1}\gz+q'(\gr_A\psi_R) (H[\eta])^{q'-1}\Gd H[\eta]\\
&+q(q'-1)(\gr_A\psi_R) (H[\eta])^{q'-2}| \nabla H[\eta]|^2\\
&+2q'(H[\eta])^{q'-1}\nabla(\gr_A\psi_R)\cdot\nabla H[\eta].
 \EAL\end{equation}
In addition,
$$\BA {l}\nabla H[\eta]=\nabla_{x'} H[\eta]+\nabla_{x''} H[\eta]
 =\prt_{y}H[\eta]\myfrac{x'}{y}+\nabla_{x''} H[\eta]\\[2mm]\phantom{\nabla H[\eta]}
 =2y\prt _{t}w(y^2,x'')\myfrac{x'}{y}+\nabla_{x''} H[\eta](x',x'')
 \EA$$
and \consy (recall that $y$ stands for $|x'|$),
$$\BAL &\nabla H[\eta]\cdot\nabla(\gr_A\psi_R)\\
&=2\prt_{t}w(y^2,x'')x'\cdot\Big (\psi_R\big(|x'|^{\kappa_+-1}(\kappa_+\frac{x'}{y}\gw_k(x'/y)+ |x'|\nabla \gw_k(x'/y))\big)+\gr_A\nabla \psi_R\Big)\\
&=2\kappa_+\prt_{t}w(y^2,x'')
|x'|^{\kappa_+}\gw_k(x'/y)=2\prt_{t}w(y^2,x'')(\kappa_+\gr_A \psi_R+
\gr_A x'\cdot\nabla\psi_R). \EAL$$ Since $w=w_R$ vanishes for
$|x''|=R$ and $\eta=0$ in a neighborhood of this sphere,
$|\prt_{t}w(y^2,x'')|\leq c\gr^R$. As $\psi_R$ vanishes for
$|x'|>3R/4$ we have $\gr_A  \nabla\psi_R\leq c\gr_A^R$. Therefore
$$
\abs{\nabla H[\eta]\cdot\nabla\gr_A}\leq c \gr^R\gr_A^R
$$
and, in view of \eqref{heat2'},
  \begin{equation}\label{Gdz}
  \abs{\Gd\gz}\leq c \gr^R\gr_A^R.
\end{equation}
  Thus $\gz\in X(D_{A}^R)$ and \consy
  \begin{equation}\label{int-eq-8.1}
\int_{D_{A}^R}\left(-u\Gd\gz+u^q\gz\right)dx=-\int_{D_{A}^R}\BBK[\gm]\Gd\gz dx.
\end{equation}
Since $q(q'-1)\gr_A (H[\eta])^{q'-2}| \nabla H[\eta]|^2\geq 0$, we have
  \begin{equation}\label{heat3}\BAL
 & \abs{\int_{D_{A}^R}u\Gd\gz dx}\\
&\leq  \int_{D_{A}^R}u\left(\gl_{1}\gz+q' (H[\eta])^{q'-1}\left(\gr|\Gd H[\eta]|+2|\nabla\gr.\nabla H[\eta]| \right)\right) dx\\
&\leq  \int_{D_{A}^R}{}u\left(\gl_{1}\gz + q' \gz^{1/q}\left(\gr^{1/q'}|\Gd H[\eta]|+2\gr^{-1/q}|\nabla\gr.\nabla H[\eta]| \right)\right) dx\\
&\leq  \left(\int_{D_{A}^R}{}u^q\gz dx\right)^{\rec{q}}
\left(\gl_{1}\left(\myint{D_{A}^R}{}\gz dx\right)^{\rec{q'}}+q'\norm{L[\eta]}_{L^{q'}(D_{A}^R)}\right) \EAL\end{equation}
where
  \begin{equation}\label{heat4}
  L[\eta]=\gr^{1/q'}|\Gd H[\eta]|+2\gr^{-1/q}|\nabla\gr.\nabla H[\eta]|.
  \end{equation}
By \rprop{IBP}
    \begin{equation}\label{heat4'}-\myint{D_{A}^R}{}\BBK[\gm]\Gd \gz dx=\myint{d_{A}^R}{}\eta^{q'}d\gm.
 \end{equation}
  Therefore
    \begin{equation}\label{heat5}\BAL
&\left( \int_{d_{A}^R}\eta^{q'} d\gm \right)\leq
\int_{D_{A}^R}{}u^q\gz dx+\\
&\q+ \left(\int_{D_{A}^R}{}u^q\gz dx\right)^{\rec{q}}
\left(\gl_{1}\left(\myint{D_{A}^R}{}\gz
dx\right)^{\rec{q'}}+q'\norm{L[\eta]}_{L^{q'}(D_{A}^R)}\right). \EAL
\end{equation}

 Next we prove that
\begin{equation}\label{heat6}
\norm{L[\eta]}_{L^{q'}(D_{A}^R)}\leq
C\norm\eta_{W^{s,q'}(\BBR^{N-k})}
\end{equation}
starting with the estimate of the first term on the right hand side
of \eqref{heat4}.

 $$\BA {l}\Gd H[\eta]=\Gd_{x'} H[\eta]+\Gd_{x''} H[\eta]
 =\prt^2_{y}H[\eta]+\myfrac{k-1}{y}\prt_{y}H[\eta]+\Gd_{x''} H[\eta]\\[2mm]\phantom{\nabla H[\eta]}
 =2y^2\prt_{tt}w(y^2,x'')+k\prt_{t}w(y^2,x'')+\Gd_{x''} H[\eta]\\[2mm]\phantom{\nabla H[\eta]}
 =2y^2\prt_{tt}w(y^2,x'')+(k+1)\prt_{t}w(y^2,x'').
 \EA$$
 Then
  $$\BA {l}\myint{\BBR^N}{}\gr \abs{\Gd H[\eta]}^{q'}dx\leq
 c\myint{0}{1}\myint{\BBR^{N-k}}{} \abs{\prt_{tt}w(y^2,x'')}^{q'}dx''y^{\gk_{+}+2q'+k-1}dy\\[4mm]
 \phantom{\myint{\BBR^N}{}\gr \abs{\Gd H[\eta]}^{q'}dx}
 +
 c\myint{0}{1}\myint{\BBR^{N-k}}{} \abs{\prt_{t}w(y^2,x'')}^{q'}dx''y^{\gk_{+}+k-1}dy\\[4mm]
 \phantom{\myint{\BBR^N}{}\gr \abs{\Gd H[\eta]}^{q'}dx}
 \leq c\myint{0}{1}\myint{\BBR^{N-k}}{} \abs{\prt_{tt}w(t,x'')}^{q'}dx''t^{(\gk_{+}+k)/2+q'}\myfrac{dt}{t}
 \\[4mm]
 \phantom{\myint{\BBR^N}{}\gr \abs{\Gd H[\eta]}^{q'}dx}
 +
c\myint{0}{1}\myint{\BBR^{N-k}}{} \abs{\prt_{t}w(t,x'')}^{q'}dx''t^{(\gk_{+}+k)/2}\myfrac{dt}{t}
\\[4mm]
 \phantom{\myint{\BBR^N}{}\gr \abs{\Gd H[\eta]}^{q'}dx}
\leq c\myint{0}{1}\norm{t^{2-(1-\frac{\gk_{+}+k}{2q'}))}\myfrac{d^2S(t)[\eta]}{dt^2}}^{q'}_{L^{q'}(\BBR^{N-k})}\myfrac{dt}{t}\\[4mm]
 \phantom{\myint{\BBR^N}{}\gr \abs{\Gd H[\eta]}^{q'}dx}
+ c\myint{0}{1}\norm{t^{1-(1-\frac{\gk_{+}+k}{2q'})}\myfrac{dS(t)[\eta]}{dt}}^{q'}_{L^{q'}(\BBR^{N-k})}\myfrac{dt}{t}.
 \EA$$
Put $\gb=\frac{\gk_{+}+k}{2q'}$ and note that
$0<\gb=\rec{2}(2-s)<1$. By standard interpolation theory,
$$\BAL
&\myint{0}{1}\norm{t^{1-(1-\gb)}\myfrac{dS(t)[\eta]}{dt}}^{q'}_{L^{q'}(\BBR^{N-k})}\myfrac{dt}{t}\\
&\approx\norm {\eta}^{q'}_{\left[W^{2,q'},L^{q'}\right]_{1-\gb,q'}}\approx
\norm\eta^{q'}_{W^{2(1-\gb),q'}
(\BBR^{N-k})},
\EAL$$
and
$$\BAL
&\myint{0}{1}\norm{t^{2-(1-\gb))}\myfrac{d^2S(t)[\eta]}{dt^2}}^{q'}_{L^{q'}(\BBR^{N-k})}\myfrac{dt}{t}\\
&\approx\norm
{\eta}^{q'}_{\left[W^{4,q'},L^{q'}\right]_{\frac{1}{2}(1-\gb),q'}}\approx
\norm\eta^{q'}_{W^{2(1-\gb),q'} (\BBR^{N-k})}. \EAL$$

The second term on the right hand side of \eqref{heat4} is estimated
in a similar way:
 $$\BAL &\int_{\BBR^N}\gr^{-q'/q}\abs{\nabla H[\eta]\cdot\nabla\gr}^{q'}dx\leq c\myint{0}{1}\int_{\BBR^{N-k}}
 \abs{\prt _{t}w(y^2,x'')}^{q'}dx'y^{\gk_{+}+k-1}dy\\
  &\leq c\myint{0}{1}\int_{\BBR^{N-k}}
 \abs{\prt _{t}w(t,x'')}^{q'}dx't^{\frac{\gk_{+}+k}{2}}\myfrac{dt}{t}\\
&\leq c\myint{0}{1}\norm{t^{1-(\frac{1}{2}-\gb)}\myfrac{d S(t)[\eta]}{dt}}_{L^{q'}(\BBR^{N-k})}^{q'}\myfrac{dt}{t}\\
&\approx\norm{\eta}^{q'}_{W^{2(1-\gb),q'}(\BBR^{N-k})}.
\EAL$$

This proves  \eqref{heat6}. Further, \eqref{heat5} and \eqref{heat6}
imply \eqref{inverse-ineq}.

We turn to the proof of part (ii). Let  $E$ be a closed subset of
$B^{N-k}_{R/2}(0)$ such that $C^{N-k}_{s,q'}(E)=0$. Then there
exists a \seq $\{\eta_{n}\}$ in $C^\infty_{0}(d_{A})$ such that
$0\leq\eta_{n}\leq 1$, $\eta_{n}=1$ in a neighborhood of $E$ (which
may depend on $n$), $\supp\eta_n\sbs B^{N-k}_{3R/4}(0)$ and
$\norm{\eta_{n}}_{W^{s,q'}}\to 0$. Then, by \eqref{heat6},
$$\norm{L[\eta_{n}]}_{L^{q'}(D_{A}^R)}\to 0.$$
Furthermore
$$\norm{w}\indx{L^{q'}((0,R)\ti B_R^{N-k}(0))}\leq c\norm{\eta_n}\indx{L^{q'}(B_R^{N-k}(0))}$$
 and \consy
 $$H[\eta_n]\to 0
 \txt{in}  L^{q'}(D_A^R).$$
 (Here we use the fact that $k\geq2$.) In addition
 $$0\leq H[\eta_n]\leq
 1,\q H[\eta_n]\leq c(R-|x'|)$$
 with a constant $c$ independent of $n$. Hence (see
 \eqref{rhoAR})
 $$\gz_{n,R}:=\gr_{A}\psi_R H[\eta_n]^{q'}\leq \gr^R\gr_{A}\psi_R H[\eta_n]^{q'-1}\leq
 \gr^R\gr_{A}^R H[\eta_n]^{q'-1}.$$
As $u^q\gr^R\gr_A^R\in L^1(D_A^R)$ we obtain,
$$\lim_{n\to\infty}\myint{D_{A}}{}u^q\gz_{n}dx=0.
$$
This fact and \eqref{heat5} imply that
$$\int_{d_{A}^R}\eta_n^{q'} d\gm\to 0.$$
As $\eta_n=1$ on a \ngh of $E$ in $\BBR^{N-k}$ it follows that
$\mu(E)=0$. \qed

\bprop{edge-bvp} Let $D_A$ be a k-dihedron, $1\leq k< N$. Let $k_+$
be as in \eqref{kappa2} and let $q^*_c$ and $q_c$ be as in
\rprop{admp0} and \rprop{admpk} respectively. Assume that $q_c\leq
q<q^*_c$. A measure $\mu\in \GTM(\prt D_A)$, with compact support
contained in $d_A$, is q-good relative to $D_A$ if and only if $\mu$
vanishes on every Borel set $E\sbs d_A$ \sth $C_{s,q'}(E)=0$, where
$s=2-\frac{k+\kappa_+}{q'}$. \es

\Remark We shall use the notation $\mu\prec C_{s,q'}$ to say that
$\mu$ vanishes on any Borel set $E\sbs(d_A)$ \sth $C_{s,q'}(E)=0$.

In the case $k=N$: $D_A=C_A$ (= the cone with vertex $0$ and opening
$A$ in $\BBR^k$) and $q_c=q^*_c$. By \rth{admiss},
$q_c=1-\frac{2}{\kappa_-}=\frac{N+\kappa_+}{N+\kappa_+ -2}$. (Note
the difference in notation; the entity denoted by $\kappa_-$ in
section 6 and in the present section is denoted by $-\ga_S$ in
subsection 5.1. See \eqref{alpha} and \eqref{kappa2}.) If $1<q<q_c$
then, again by \rth{admiss}, there exist solutions for every measure
$\mu=k\gd_0$ on $\prt C_A$.

In the case $k=1$, $q^*_c=\infty$, $\kappa_+=1$ and
$q_c=\frac{N+1}{N-1}$. Thus $s=2/q$ and the statement of the theorem
is well known (see \cite{MV3}).

\medskip
\Proof In view of the last remark, it remains to deal only with
$2\leq k\leq N-1$. We shall identify $d_A$ with $\BBR^{N-k}$.

It is sufficient to prove the result for positive measures because
$\mu\prec C_{s,q'}$ if and only if $|\mu|\prec C_{s,q'}$. In
addition, if $|\mu|$ is a q-good measure then $\mu$ is a q-good
measure.

First we show that if $\mu$ is non-negative and q-good then
$\mu\prec C_{s,q'}$. If $E$ is a Borel subset of $\bdw$ then
$\mu\chi\indx{E}$ is q-good. If $E$ is compact and $C_{s,q'}(E)=0$
then, by \rprop{admi-conv}, $E$ is a removable set. This means that
the only solution of \eqref{bvp8.1} \sth $\mu(\bdw\sms E)=0$ is the
zero solution. This implies that $\mu\chi\indx{E}=0$, i.e.,
$\mu(E)=0$. If $C_{s,q'}(E)=0$ but $E$ is not compact then
$\mu(E')=0$ for every compact set $E'\sbs E$. Therefore, we conclude
again that $\mu(E)=0$.

Next, assume that $\mu$ is a positive measure in $\GTM(\prt D_A)$
supported in a compact subset of $\BBR^{N-k}$.

If $\mu\in B^{-s,q}(\BBR^{N-k})$ then, by \rth{main1} and
\rth{admissible}, $\mu$ is q-good relative to $D_A\cap \Gg_{k,R}$,
for every $R>0$. (As before $\Gg_{k,R}$ is the cylinder with radius
$R$ around the 'axis' $\BBR^{N-k}$.) This implies that $\mu$ is
q-good relative to $D_A$.

If $\mu\prec C_{s,q'}$ then, by a theorem of Feyel and de la
Pradelle \cite{FD} (see also \cite{BP}), there exists a \seq
$\{\mu_n\}\sbs (B^{-s,q}(\BBR^{N-k}))_+$ \sth $\mu_n\uparrow\mu$. As
$\mu_k$ is q-good,  it follows that  $\mu$ is q-good. \qed

\bth{poly-good} Let $P$ be an $N$-dimensional polyhedron as
described in \rprop{edge condition}.  Let $\mu$ be a bounded measure
on $\prt P$, \(may be a signed measure\). Let
$k=1,\dots,N,\;j=1,\dots, n_k$, and let ${L_{k,j}}$ and $A_{k,j}$ be
defined as at the beginning of section 8. Further, put
\begin{equation}\label{skj}
 s(k,j)=2-\frac{k+(\kappa_+)_{k,j}}{q'},
\end{equation}
where $(\kappa_+)_{k,j}$ is defined as in \eqref{kappa2} with
$A=A_{k,j}$.
Then $\mu\in \GTM_q^P$, i.e., $\mu$ is a good measure for
\eqref{eq-q} relative to $P$, if and only if, for every pair $(k,j)$
as above and every Borel set $E\sbs L_{k,j}$:\\
 If $1\leq k<N$ then
\begin{equation}\label{poly-good1}\BAL
 (q_c)_{k,j}\leq q<(q^*_c)_{k,j},\; C^{N-k}_{s(k,j),q'}(E)=0 &\Lra
 \mu(E)=0 \\
 q\geq (q_c^*)_{k,j} &\Lra
 \mu(L_{N,j})=0
\EAL\end{equation}
and if $k=N$, i.e., $L$ is a vertex,
\begin{equation}\label{poly-good2}
q\geq
(q_c)_{k,j}=\frac{N+2+\sqrt{(N-2)^2+4\gl_A}}{N-2+\sqrt{(N-2)^2+4\gl_A}}\Lra
\mu(L)=0.
\end{equation}
Here $(q_c^*)_{k,j}$ and $(q_c)_{k,j}$ are defined  as in \eqref{qk}
and \eqref{q-critk}respectively,  with $A=A_{k,j}$.

If $1<q<(q_c)_{k,j}$ then there is no restriction on
$\mu\chi\indx{L_{k,j}}$. \es

\Proof This is an immediate consequence of \rprop{edge condition}
and \rprop{edge-bvp} (see also the Remark following it). In the case
$k=N$, $L_{N,j}$ is a vertex and the condition says merely that for
$q\geq q_c(L_{N,j}$, $\mu$ does not charge the vertex. \qed

\subsection{Removable singular sets II.}

\bprop{uqro} Let $A$ be a \Lip domain on $S^{k-1}$, $2\leq k\leq
N-1$, and let $D_A$ be the k-dihedron with opening $A$. Let $u$ be a
positive solution of \eqref{eq-q} in $D_A^R$, for some $R>0$.
Suppose that $F=\CS(u)\sbs d_A^R$ and let $Q$ be an open \ngh of $F$
\sth $\bar Q\sbs d_A^R$. \(Recall that $d_A^R=d_A\cap B^{N-k}_R(0)$
is an open subset of $d_A$.\) Let $\mu$ be the trace of $u$ on
$\CR(u)$.

Let $\eta\in W_0^{s,q'}(d_A^R)$ \sth
\begin{equation}\label{eta-cond}
 0\leq\eta \leq 1,\q \eta=0 \txt{on $Q$.}
\end{equation}
Employing the notation in the proof of \rprop{admi-conv}, put

\begin{equation}\label{gz-eta}
  \gz:=\gr_{A}\psi_RH_R[\eta]^{q'}.
\end{equation}
 Then
\begin{equation}\label{uqro1}
  \int_{D_A^R} u^q \gz\,dx\leq
  c(1+\norm{\eta}\indx{W^{s,q'}(d_A)})^{q'}+\mu(d_A^R\sms Q)^q,
\end{equation}
$c$ independent of $u$ and $\eta$.

\es

\Proof First we prove \eqref{uqro1} for $\eta\in C_0^\infty(d_A^R)$.
Let $\gs_0$ be a point in $A$ and let $\set{A_n}$ be a \Lip
exhaustion of $A$. 
If $0<\ge<\dist(\prt A,\prt A_n)=\bar \ge_n$ then

$$\ge\gs_0+C_{A_n}\sbs C_A.$$
Denote
$$D_A^{R', R''}=D_A\cap[|x'|<R']\cap[|x''|<R''].$$
Pick a \seq $\{\ge_n\}$ decreasing to zero \sth
$0<\ge_n<\min(\bar\ge_n/2^n, R/8)$. Let $u_n$ be the function given
by
$$u_n(x'x'')=u(x'+\ge_n\gs_0,x'') \forevery x\in D_{A_n}^{R_n,R},
\q R_n=R-\ge_n.$$ Then $u_n$ is a solution of \eqref{eq-q} in
$D_{A_n}^{R_n,R}$ belonging to $C^2(\bar D_{A_n}^{R_n,R})$ and we
denote its boundary trace by $h_n$. Let
$$\gz_n:=\gr_{A_n}\psi_RH_R[\eta]^{q'},$$
with $\psi_R$ and $H_R[\eta]$ as in the proof of \rprop{admi-conv}.
By \rprop{IBP}
    \begin{equation}\label{heat4''}-\int_{D_{A_n}^{R_n,R}}\BBP[h_n]\Gd \gz_n dx=\int_{B_R^{N-k}(0)}\eta^{q'}h_nd\gw_n
 \end{equation}
 where $\gw_n$ is the harmonic measure on $d_{A_n}^R$ relative to $D_{A_n}^{R_n,R}$. (Note that $d_{A_n}^R=d_{A}^R$
 and we may identify it with  $B_R^{N-k}(0)$.) Hence
  \begin{equation}\label{eq-8.2}
\int_{D_{A_n}^{R_n,R}}\left(-u_n\Gd\gz_n+u_n^q\gz_n\right)dx=-\int_{B_R^{N-k}(0)}\eta^{q'}h_n
\,d\gw_n.
\end{equation}
Further,
$$\int_{B_R^{N-k}(0)}\eta^{q'}h_n \,d\gw_n\to
\int_{B_R^{N-k}(0)}\eta^{q'}d\mu\leq \mu(d_A^R\sms Q),$$ because
$\eta=0$ in $Q$. By \eqref{heat3}, \eqref{heat6} we obtain,
  \begin{equation}\label{ineq-8.2}\BAL
&\abs{\int_{D_{A_n}^{R_n,R}}u_n\Gd\gz_n\,dx}\leq\\
c\Big(&\int_{D_{A_n}^{R_n,R}}u_n^q\gz_n dx\Big)^{\rec{q}}
\Big(\Big(\int_{D_{A_n}^{R_n,R}}\gz_n
dx\Big)^{\rec{q'}}+\norm{\eta}_{W^{s,q'}(B^{N-k}_R(0))}\Big).
\EAL\end{equation}

\nind From the definition of $\gz_n$ it follows that
$$\int_{D_{A_n}^{R_n,R}}\gz_n\,dx\leq \int_{D_{A_n}^{R_n,R}}\gr_n\,dx\to \int_{D_{A}^{R}}\gr\,dx,$$
where $\gr$ (resp. $\gr_n$)  is the first eigenfunction of $-\Gd$ in
$D_A^R$ \(resp. $D_{A_n}^{R_n,R}$\) normalized by $1$ at some
$x_0\in D_{A_1}^{R_1,R}$. Therefore, by \eqref{eq-8.2},

$$\BAL&\int_{D_{A_n}^{R_n,R}}u_n^q\gz_ndx\leq c\Big(&\int_{D_{A_n}^{R_n,R}}u_n^q\gz_n dx\Big)^{\rec{q}}
\big(1+\norm{\eta}_{W^{s,q'}(B^{N-k}_R(0))}\big) + \mu(d_A^R\sms
Q).\EAL$$
This implies
\begin{equation}\label{ineq-8.3'}
   \int_{D_{A_n}^{R_n,R}}u_n^q\gz_ndx\leq
   c\big(1+\norm{\eta}_{W^{s,q'}(B^{N-k}_R(0))}\big)^{q'}+\mu(d_A^R\sms
   Q)^q.
\end{equation}
To verify this fact, put
$$\BAL m=\Big(\int_{D_{A_n}^{R_n,R}}u_n^q\gz_ndx\Big)^{1/q},\; b=\mu(d_A^R\sms
Q),\; a=c\big(1+\norm{\eta}_{W^{s,q'}(B^{N-k}_R(0))}\big)\EAL$$ so
that \eqref{ineq-8.3'} becomes
$$m^q-am-b\leq 0.$$
If $b\leq m$ then
$$m^{q-1}-a-1\leq 0.$$
Therefore,
$$m\leq (a+1)^{\rec{q-1}}+b$$
which implies \eqref{ineq-8.3'}. Finally, by  the lemma of Fatou we
obtain \eqref{uqro1} for $\eta\in C_0^\infty$. By continuity we
obtain the inequality for any $\eta\in W^{s,q'}_0$ satisfying
\eqref{eta-cond}. \qed

\bth{removable} Let $A$ be a \Lip domain on $S^{k-1}$, $2\leq k\leq
N-1$, and let $D_A$ be the k-dihedron with opening $A$.
  Let $E$ be a compact subset of $d_A^R$ and let $u$ be a non-negative solution
 of \eqref{eq-q} in $D_A^R$ (for some $R>0$) \sth $u$ vanishes on
 $\prt D_A^R\sms E$. Then
\begin{equation}\label{REM1}
C^{N-k}_{s,q'}(E)=0,\q s=2-\frac{\gk_{+}+k}{q'} \Longrightarrow u=0,
\end{equation}
where $C^{N-k}_{s,q'}$ denotes the Bessel capacity with the
indicated indices in $\BBR^{N-k}$. \es

\Proof By \rprop{admi-conv}, \eqref{REM1} holds under the additional
assumption
\begin{equation}\label{REM2'}
\int_{D_A^R}u^q\gr_R\gr_A^R dx<\infty.
\end{equation}
 Indeed, by \rprop{reg-sol},
\eqref{REM2'} implies that the solution $u$ possesses a boundary
trace $\mu$ on $\prt D_A^R$. By assumption, $\mu(\prt D_A^R\sms
E)=0$. Therefore, by \rprop{edge-bvp}, the fact that
$C^{N-k}_{s,q'}(E)=0$ implies that $\mu(E)=0$. Thus $\mu=0$ and
hence $u=0$.

We show that, under the conditions of the theorem, if
$C^{N-k}_{s,q'}(E)=0$ then \eqref{REM2'} holds.

By \rprop{uqro}, for every $\eta\in W_0^{s,q'}(d_A^R)$ \sth $0\leq
\eta\leq 1$ and $\eta=0 \txt{in a \ngh of $E$,}$
\begin{equation}\label{ineq-8.3}
   \int_{D_{A}^{R}}u^q\gz\,dx\leq
   c\big(1+\norm{\eta}_{W^{s,q'}(B^{N-k}_R(0))}\big)^{q'},
\end{equation}
for $\gz$ as in \eqref{gz-eta}. (Here we use the assumption that
$u=0$ on $\prt D_A^R\sms E$.)

 Let $a>0$ be sufficiently small so that $E\sbs
B_{(1-4a)R}^{N-k}(0)$. Pick a \seq $\{\phi_n\}$  in
$C_0^\infty(\BBR^{N-k})$ \sth, for each $n$, there exists a \ngh
$Q_n$ of $E$,
   $\bar Q_n\sbs B^{N-k}_{(1-3a)R}(0)$ and
\begin{equation}\label{phin-1}\BAL
   &0\leq \phi_n\leq 1 \txt{everywhere,} \phi_n=1 \txt{in $Q_n$,}\\
   &\tl\phi_n:=\phi_n\chi\indx{[|x''|<(1-2a)R]}\in
   C_0^\infty(\BBR^{N-k}),\\
&\big\|\tl\phi_n\big\|_{W^{s,q'}(\BBR^{N-k})}\to 0\txt{as} n\tin\\
&\eta_n:=(1-\phi_n)\lfloor\indx{[|x''|<R]}\in
   C_0^\infty(d_A^R),\\ & \eta_n=0 \;\text{in  }[(1-a)R<|x''|<R].
 \EAL\end{equation}

\nind  Such a \seq exists because $C^{N-k}_{s,q'}(E)=0$. Applying
\eqref{ineq-8.3} to $\eta_n$ we obtain,
\begin{equation}\label{ineq-8.4}
 \sup  \int_{D_{A}^{R}}u^q\gz_n\,dx\leq c<\infty,
\end{equation}
where $\gz_n=\gr_A\psi_RH_R^{q'}[\eta_n]$ (see \eqref{gz-eta}). By
taking a \sseq we may assume that $\{\eta_n\}$ converges (say to
$\eta$) in $L^{q'}(B_{R}^{N-k}(0))$ and \consy $H[\eta_n]\to
H[\eta]$ in the sense that
$$H_R[\eta_n](x',\cdot)=w_{n,R}(y^2,\cdot)\to
w_{R}(y^2,\cdot)=H_R[\eta](x',\cdot) \txt{in $L^{q'}$}$$ uniformly
\wrto $y=|x'|$.  It follows that
\begin{equation}\label{ineq-8.5}
\int_{D_{A}^{R}}u^q\gz\,dx<\infty,\q  \gz=\gr_A\psi_RH_R^{q'}[\eta].
\end{equation}

As $\tl\phi_n\to 0$ in $W^{s,q'}(\BBR^{N-k})$ it follows that
$\phi_n\to 0$ and hence $\eta_n\to 1$ a.e. in
$B_{(1-2a)R}^{N-k}(0)$. Thus $\eta=1$ in this ball, $\eta=0$ in
$[(1-a)R<|x''|<R]$ and $0\leq \eta\leq 1$ everywhere.

 \Consy, given
$\gd>0$, there exists an $N$-dimensional \ngh $O$ of $d_A\cap
B_{(1-2a)R}^{N-k}(0)$ \sth
$$1-\gd<H_R[\eta]<1 \txt{and} 1-\gd<\psi_R/\gr_A^R<1 \txt{in}O.$$
Therefore \eqref{ineq-8.5} implies that
\begin{equation}\label{ineq-8.6}
   \int_{D_{A}^{(1-3a)R}}u^q\gr^R\gr_A^R\,dx\leq c<\infty.
\end{equation}
Recall that the trace of $u$ on  $\prt D_A^R\sms d_A^{(1-4a)R}$ is
zero. Therefore $u$ is bounded in $D_A^R\sms D_{A}^{(1-3a)R}$. This
fact and \eqref{ineq-8.6} imply \eqref{REM2'}. \qed

\bdef{ro-capacity} Let $\Gw$ be a bounded \Lip domain. Denote by
$\gr$ the first eigenfunction of $-\Gd$ in $\Gw$ normalized by
$\gr(x_0)=1$ for a fixed point $x_0\in\Gw$.

For every compact set $K\sbs \bdw$ we define
$$ M_{\gr, q}(K)=\set{\mu\in \GTM(\bdw):\mu\geq 0,\;\mu(\bdw\sms
K)=0,\; \BBK[\mu]\in L^1_\gr(\Gw)}$$
and

$$\tl C_{\gr,q'}(K)=\sup\set{\mu(K)^q:\;\mu\in M_{\gr, q}(K),\;\int_\Gw
\BBK[\mu]^{q}\gr\,dx=1}.$$

Finally we denote by $C_{\gr,q'}$ the outer measure generated by the
above functional. \es

The following statement is verified by standard arguments:

\blemma{cgrq} For every compact $K\sbs \bdw$, $C_{\gr,q'}(K)=\tl
C_{\gr,q'}(K)$. Thus $C_{\gr,q'}$ is a capacity and,
\begin{equation}\label{cgrq}
    C_{\gr,q'}(K)=0 \iff M_{\gr, q}(K)=\{0\}.
\end{equation}
\es

\bth{removable set} Let $\Gw$ be a bounded polyhedron in $\BBR^N$. A
compact set $K\sbs \bdw$ is removable if and only if
\begin{equation}\label{remov1}
  C_{s(k,j),q'}(K\cap L_{k,j})=0,
\end{equation}
for $k=1,\cdot,N$ $j=1,\cdots, n_k$, where $s(k,j)$ is defined as in
\eqref{skj}. This condition is equivalent to
\begin{equation}\label{remov1'}
  C_{\gr,q'}(K)=0.
\end{equation}
A measure $\mu\in\GTM(\bdw)$ is q-good if and only if it does not
charge sets with $C_{\gr,q'}$-capacity zero.
 \es

\Proof The first assertion is an immediate consequence of
\rprop{edge condition} and \rth{removable}. The second assertion
follows from the fact that
$$C_{\gr,q'}(K\cap L_{k,j})=C_{s(k,j),q'}(K\cap L_{k,j}).$$
The third assertion follows from \rth{poly-good} and the previous
statement. \qed


\mysection{Appendix--Boundary Harnack inequality}
In this section we prove the following

\bprop {BHI} Assume $\Gw$ is a bounded Lipschitz domain,
$A\subset\prt\Gw$ is relatively open and $q>1$. Let $(r_0,\gl_0)$ be
the \Lip characteristic of $\Gw$ (see subsection 2.1).

 Let $u_{i}\in C(\Gw\cup A)$, $i=1,2$, be positive solutions of
$$-\Gd u+u^q=0\quad\text {in }\;\Gw,
$$
such that  $u_{i}=0$ on $A$. Put $S=\bdw\sms A$ and
$d(x,S)=\dist(x,S)$. Let $y\in A$  and let
$$r:=\min(r_0/8,\rec{4}d(y,S)$$
so that
 $$\prt (B_{4r}(y)\cap\Gw)=(\overline
B_{4r}(y)\cap\prt\Gw)\cup (\prt B_{4r}(y)\cap\Gw).$$ Then
\begin{equation}\label{pot3}
c^{-1}\frac{u_{1}(z')}{u_{1}(z)}\leq \frac{u_{2}(z')}{u_{2}(z)}\leq
c\frac{u_{1}(z')}{u_{1}(z)} \forevery z,z'\in\,B_{r}(y)\cap\Gw,
\end{equation}
 where the constant $c>0$ depends only on $N,q$ and the \Lip characteristic of $\Gw$.
\es

\medskip
\Proof  Without loss of generality we assume that $y=0$.

Let $b=d(0,S)$ and put
$$\tl u_i(x)=b^{-\frac{2}{q-1}}u_1(x/b), \q i=1,2.$$
Then $\tl u_i$ has the same properties as $u_i$ when $\Gw$ is
replaced by $\Gw^b=\rec{b}\Gw$, $S$ by $ S^b=\rec{b}S$ and $r$ by
$\gd=r/b$. Of course $d(0,S^b)=1$ so that
$$\gd=\min(r_0/(8b),1/4).$$
The functions $\tl u_i$ satisfy the equation
$$-\Gd\tl u_i+\tl u_i^q=0 \txt{in} B_{4\gd}(0)\cap \Gw^b$$
and $\tl u_i=0$ on $B_{4\gd}(0)\cap \bdw^b$. Therefore, by the
Keller--Osserman estimate,
$$\tl u_i\leq c(N,q)\gd^{-2/(q-1)} \txt{in} \bar B_{3\gd}(0)\cap\Gw^b.$$

If $a(x)=\tl u_{1}^{q-1}$ then $\tl u_{1}$ satisfies
$$-\Gd\tl u_1+a(x)\tl u_1=0\quad\text {in }\;(\rec{b}\Gw)\cap B_{1}(0),
$$
and $a(\cdot)$ is bounded in $\bar B_{3\gd}(0)$.

Let $w$ be the solution of
$$\left\{\BAL
-\Gd w+a(x)w&=0\quad &&\text{in }\;B_{3\gd}(0)\cap\Gw^b\\
w&=0 &&\text{on }\;\overline B_{3\gd}(0)\cap\rec{b}\prt\Gw\\
w&=\tl u_{2}&&\text{on }\;\prt B_{3\gd}(0)\cap\bdw^b. \EAL\right.$$

\nind By applying the boundary Harnack principle in
$B_{3\gd}(0)\cap\Gw^b$ (using the slightly more general form derived
in \cite[Theorem 2.1]{Ba})  we obtain
\begin{equation}\label{pot0}
c^{-1}\frac{\tl u_{1}(\gz')}{\tl u_{1}(\gz)}\leq
\frac{w(\gz')}{w(\gz)}\leq c\frac{\tl u_{1}(\gz')}{\tl u_{1}(\gz)}
\forevery \gz,\gz'\in B_{2\gd}(0)\cap\Gw^b,
\end{equation}
where the constant $c$ depends only on the \Lip characteristic of
$\Gw^b$ (which is $(r_0/b,\gl_0b)$ and therefore 'better' then that
of $\Gw$ when $b\leq 1$). Since $w\leq \tl u_2$ the above inequality
implies,
$$\frac{w(\gz')}{\tl u_2(\gz)}\leq c\frac{\tl u_{1}(\gz')}{\tl
u_{1}(\gz)}\txt{and} c^{-1}\frac{\tl u_{1}(\gz')}{\tl
u_{1}(\gz)}\leq \frac{\tl u_{2}(\gz')}{w(\gz)}$$ which in turn
implies
$$\frac{w(\gz')}{w(\gz)}\leq c\frac{\tl u_{2}(\gz')}{\tl u_{2}(\gz)}
\forevery \gz,\gz'\in B_{2\gd}(0)\cap\Gw^b$$ and therefore
$$\frac{\tl u_{1}(\gz')}{\tl u_{1}(\gz)}\leq c^2 \frac{\tl u_{2}(\gz')}{\tl u_{2}(\gz)}
\forevery \gz,\gz'\in B_{2\gd}(0)\cap\Gw^b.$$ Switching the roles of
$\tl u_1$ and $\tl u_2$ we obtain,
$$\frac{\tl u_{2}(\gz')}{\tl u_{2}(\gz)}\leq c^2 \frac{\tl u_{1}(\gz')}{\tl u_{1}(\gz)}
\forevery \gz,\gz'\in B_{2\gd}(0)\cap\Gw^b.$$ This completes the
proof. \qed



\end {document}